\documentclass[12pt]{book}
\usepackage{amsfonts}
\usepackage{epsfig}

\title{The Phase Transition in Five Point Energy Minimization}
\author{Richard Evan Schwartz}

\newtheorem{theorem}{Theorem}[section]

\newtheorem{lemma}[theorem]{Lemma}

\newtheorem{corollary}[theorem]{Corollary}

\def\startproof{{\bf {\medskip}{\noindent}Proof: }}

\def\endproof{$\spadesuit$  \newline}

\DeclareFontFamily{U}{rcjhbltx}{}
\DeclareFontShape{U}{rcjhbltx}{m}{n}{<->rcjhbltx}{}
\DeclareSymbolFont{hebrewletters}{U}{rcjhbltx}{m}{n}

\let\aleph\relax\let\beth\relax
\let\gimel\relax\let\daleth\relax

\DeclareMathSymbol{\aleph}{\mathord}{hebrewletters}{39}
\DeclareMathSymbol{\beth}{\mathord}{hebrewletters}{98}
\DeclareMathSymbol{\gimel}{\mathord}{hebrewletters}{103}
\DeclareMathSymbol{\daleth}{\mathord}{hebrewletters}{100}

\DeclareMathSymbol{\lamed}{\mathord}{hebrewletters}{108}
\DeclareMathSymbol{\mem}{\mathord}{hebrewletters}{109}
\DeclareMathSymbol{\ayin}{\mathord}{hebrewletters}{96}
\DeclareMathSymbol{\tsadi}{\mathord}{hebrewletters}{118}
\DeclareMathSymbol{\qof}{\mathord}{hebrewletters}{114}
\DeclareMathSymbol{\shin}{\mathord}{hebrewletters}{152}

\def\D{\mbox{\boldmath{$D$}}}% 
\def\N{\mbox{\boldmath{$N$}}}% 
\def\Q{\mbox{\boldmath{$Q$}}}% 
\def\R{\mbox{\boldmath{$R$}}}% 
\def\T{\mbox{$TBP$}}% 
\def\Z{\mbox{\boldmath{$Z$}}}% 
\def\#{\sharp}%

\begin{document}
\maketitle

\tableofcontents
\newpage

\section*{Preface}

This monograph is a rigorous computer assisted
analysis of the extent to
which the triangular bi-pyramid (TBP) is the
minimizer, amongst $5$-point configurations
on the sphere, with respect to a power law
potential.  Our main result settles a long-standing
open question about this, but perhaps some people
will object to both the length of the solution
and to the computer-assisted nature.  On the
positive side, the proof divides neatly into
$4$ modular pieces, and there is a clear idea
behind the attack at each step.  

Define the Riesz $s$-potential 
$R_s(r)={\rm sign\/}(s)\ r^{-s}.$
We prove
there exists a computable number 
$\shin=15.0480773927797...$
such that the TBP is the unique
minimizer with respect to
$R_s$ if and only if 
$s \in (-2,0) \cup (0,\shin)$.
The Hebrew letter $\shin$ is pronounced ``shin''.
Our result does not say what happens for
power law potentials much beyond $\shin$,
though we do prove that in the tiny
interval
$(\shin,15+25/512]$
 the minimizer is a pyramid with square
base. Thus, we prove the existence of the
 phase transition that was conjectured to
exist in the 1977 paper [{\bf MKS\/}] by
T. W. Melnyk, O. Knop, and W. R. Smith.

This is still a draft of the monograph.
The mathematics is all written down and the
computer programs have all been run, but I
am not sure if the programs are in their
final form. While the programs are now
pretty well documented and organized, I
might still like to do more in this
direction. 
Meanwhile, I am happy to answer any questions about the code.

I thank Henry Cohn, 
John Hughes, Abhinav Kumar,
Curtis McMullen, Ed Saff, and Sergei Tabachnikov,
for discussions related to this monograph.
I would specially like to thank Alexander Tumanov for
his great idea about this problem, and Jill Pipher
for her enthusiasm and encouragement.
I thank my wife Brienne Brown for suggesting the name
of the phase-transition constant.  I thank my
family for their forbearance while I worked on
this project like Ahab going after the whale.

I thank I.C.E.R.M. for facilitating
this project. My interest in this problem
was rekindled during discussions about point configuration
problems around the time of the I.C.E.R.M. Science
Advisory Board meeting in Nov. 2015.  I thank
the National Science Foundation for their
continued support, currently in the form of
the Research Grant DMS-1204471.  Finally, I
thank the Simons Foundation for their support,
in the form of a Simons Sabbatical Fellowship.

\newpage

\part{Introduction and Outline}

\chapter{Introduction}

\section{The Energy Minimization Problem}

Let $S^2$ denote the unit sphere in $\R^3$ and
let $P=\{p_1,...,p_n\}$ be a finite list of 
distinct points on $S^2$.  
Given some function $f: (0,2] \to \R$ we can
compute the total $f$-potential
\begin{equation}
{\cal E\/}_f(P)=\sum_{i<j} f(\|p_i-p_j\|).
\end{equation}
For fixed $f$ and $n$,
one can ask which configuration(s)
minimize ${\cal E\/}_f(P)$. 

For this problem, the energy functional $f=R_s$, where
\begin{equation}
R_s(r)={\rm sign\/}(s) r^{-s},
\end{equation}
is a natural one to consider.
When $s>0$, this is called
the {\it Riesz potential\/}.
When $s<0$ this is called the
{\it Fejes-Toth potential\/}.
The case $s=1$ is specially called
the {\it Coulomb potential\/} or the
{\it electrostatic potential\/}.
This case of the energy minimization
problem is known as {\it Thomson's problem\/}.
See [{\bf Th\/}].
The case of $s=-1$, in which one tries
to maximize the sum of the distances, is
known as {\it Polya's problem\/}. 

There is a large literature on the energy 
minimization problem. See [{\bf F\"o\/}] and
[{\bf C\/}] for some early local results.
See [{\bf MKS\/}] for a definitive numerical study
on the minimizers of the Riesz potential for $n$ 
relatively small. The website
[{\bf CCD\/}] has a compilation of experimental
results which stretches all the way up to 
about $n=1000$. 
The paper [{\bf SK\/}] gives a nice survey of results,
with an emphasis on the case when $n$ is large. 
See also [{\bf RSZ\/}].
The paper [{\bf BBCGKS\/}] gives a survey
of results, both theoretical and experimental, about
highly symmetric configurations in higher dimensions.

When $n=2,3$ the problem is fairly trivial.
In [{\bf KY\/}] it is shown that when
$n=4,6,12$, the most symmetric
configurations -- i.e. vertices of the relevant
Platonic solids -- are the unique 
minimizers for all $R_s$
with $s \in (-2,\infty)-\{0\}$.  See 
[{\bf A\/}] for just the case $n=12$ and
see [{\bf Y\/}] for an earler, partial result
in the case $n=4,6$. The result in
[{\bf KY\/}] is contained in the much more
general and powerful result
[{\bf CK\/}, Theorem 1.2] concerning 
the so-called sharp configurations.

The case $n=5$ has been notoriously intractable.
There is a general feeling that 
for a wide range of energy choices, and in particular
for the Riesz potentials,  the
global minimizer
is either the triangular bi-pyramid or else some
pyramid with square base.
The {\it triangular bi-pyramid\/} (TBP) is
the configuration of $5$ points with one
point at the north pole, one at the south
pole, and three arranged in an equilateral
triangle around the equator.  The pyramids
with square base come in a $1$-parameter
family.  Following Tumanov [{\bf T\/}] we
call them {\it four pyramids\/} (FPs).

[{\bf HS\/}] has a rigorous computer-assisted
proof that the TBP is the unique minimizer
for $R_{-1}$ (Polya's problem) and my paper [{\bf S1\/}] has a
rigorous computer-assisted proof that the TBP
is the unique minimizer for $R_1$ (Thomson's problem)
and $R_2$.
The paper [{\bf DLT\/}] gives a traditional
proof that the TBP is the unique minimizer
for the logarithmic potential.

In [{\bf BHS\/}, Theorem 7] it is shown that,
as $s \to \infty$, any sequence of
$5$-point minimizers w.r.t. $R_s$ must converge (up to
rotations) to the FP having one point at the north pole
and the other $4$ points on the equator.
In particular, the TBP is not a minimizer
w.r.t $R_s$ when $s$ is sufficiently large.
In 1977, T. W. Melnyk, O. Knop, and W. R. Smith, [{\bf MKS\/}]
conjectured the existence of the phase transition
constant, around $s=15.04808$, at which point the
TBP ceases to be the minimizer w.r.t. $R_s$.

Define
\begin{equation}
G_k(r)=(4-r^2)^k, \hskip 30 pt
k=1,2,3,...
\end{equation}
In [{\bf T\/}], A. Tumanov gives a traditional
proof of the following result.
\begin{theorem}[Tumanov]
\label{tumanov2X}
Let $f=a_1G_1+a_2G_2$ with $a_1,a_2>0$.  The
TBP is the unique global minimizer with respect to $f$.
Moreover, a critical point of $f$ must be the TBP
or an FP.
\end{theorem}

As an immediate corollary, the TBP is a minimizer
for $G_1$ and $G_2$.
Tumanov points out that these potentials do not have
an obvious geometric interpretation, but they are
amenable to a traditional analysis.
Tumanov also mentions
that his result might be a step towards proving that the TBP
minimizes a range of power law potentials.  He
makes the following obserbation.
\newline
\newline
{\bf Tumanov's Observation:\/}
If the TBP is the unique 
minimizer for $G_2$, $G_3$ and $G_5$,
then the TBP is the unique minimizer
for $R_s$ provided that $s \in (-2,2]-\{0\}$.
\newline

Tumanov does not offer a proof for his
observation, but we will prove related results
in this monograph.
Tumanov also remarks that it seems hard to
establish that the TBP is the minimizer w.r.t  
$G_3$ and $G_5$.
The family of potentials $\{G_k\}$ behaves somewhat
like the Riesz potentials. For instance, the
TBP is not the global minimizer w.r.t. $G_7$.
I checked that the TBP is not the global
minimizer for $G_8,...,G_{100}$, and
I am sure that this fact remains true
for $k>100$ as well.

\section{Main Results}

The idea of this monograph is to combine the
divide-and-conquer approach in [{\bf S1\/}] with
(elaborations of) Tumanov's observation (and
other tricks) to 
rigorously establish the phase transition 
explicitly conjectured in [{\bf MKS\/}].
We will show that the TBP is the unique
global minimizer for
$$G_3, \hskip 10 pt
G_4, \hskip 10 pt
G_5, \hskip 10 pt
G_6, \hskip 10 pt
G_5-25 G_1, \hskip 10 pt
G_{10}+28 G_5+102 G_2,
$$
and that most configurations (in a sense
made precise in the next chapter) have higher
$G_{10}+13 G_5+68G_2$ energy than the TBP.
We will then use these facts to establish the following two
results.

\begin{theorem}
\label{aux2}
Let $s \in (-2,0)$ be arbitrary.
Amongst all $5$-point configurations
on the sphere, the TBP is
the unique global minimizer w.r.t. $R_s$.
\end{theorem}

\begin{theorem}
\label{aux1}
Let $s \in (0,15+25/512]$ be arbitrary.
Amongst all $5$-point configurations
on the sphere, either the TBP or
some FP is
the unique global minimizer w.r.t. to
$R_s$.  
\end{theorem}

We chose $15+25/512=15.04809...$ because this is a dyadic
rational that is very slightly larger than the
cutoff in [{\bf MKS\/}].
The above two results combine with some
elementary analysis to
establish our main result.

\begin{theorem}[Main]
\label{aux3}
There exists a computable number 
$$\shin=15.0480773927797...$$
such that the TBP is the unique global
minimizer w.r.t.
the $R_s$ energy if and only if 
$s \in (-2,0) \cup (0,\shin)$.
\end{theorem}

Our results say nothing about what happens
beyond $s=15+25/512$ (though 
our Symmetrization Lemma in \S \ref{symm0}
perhaps suggests a path onward; see
the discussion in \S \ref{finalYY}.) Even so,
the Main Theorem is the definitive result about
the extent that the TBP minimizes power law
potentials.  In particular, it establishes the
existence of the phase transition constant which
has been conjectured to exist for many decades.

The constant
$\shin$ is probably some transcendental number that
is unrelated to all other known constants of nature
and mathematics.  So, saying that an ideal computer can
compute $\shin$ to arbitrary precision
is the best we can do.  In \S \ref{cc} we will
sketch an algorithm which computes $\shin$
to arbitrary precision.

\section{Ideas in the Proof}

In the next chapter we will give a detailed
outline of the proofs of the results above.
  Here we make
some general comments about the strategy.

We are interested in searching through the
moduli space of all $5$-point configurations
and eliminating those which have higher
energy than the TBP.   We win if we eliminate
everything but the TBP.   Assuming that the
TBP really is the minimizer for a given 
energy function $f$, it is probably the case
that most configurations are not even
close to the TBP in terms of $f$-energy.  So, in principle, 
one can eliminate most of the configuration
space just by crude calculations.
If this elimination procedure works well,
then what is left is just a small neighborhood
$B$ of the TBP.  The TBP
is a critical point for ${\cal E\/}_f$, by symmetry,
and with luck the Hessian of
${\cal E\/}_f$ is positive definite throughout $B$.
In this case, we can say that the TBP must
be the unique global minimizer.

Most of our proof involves implementing
such an elimination scheme.  We use a 
divide-and-conquer approach. We
normalize so that $(0,0,1)$ is one of
the points of the configuration, and then
we use stereographic projection to
move this point to $\infty$ and the
other $4$ points to $\R^2$.
Stereographic projection is the map
\begin{equation}
\label{stereo}
\Sigma(x,y,z)=\bigg(\frac{x}{1-z}, \frac{y}{1-z}\bigg).
\end{equation}
We call these points {\it stereographic coordinates\/}
for the configuration.  Stereographic coordinatss
give the moduli space a natural flat
structure, making divide-and-conquer
algorithms easy to manage.
Our basic object is a {\it block\/}, a
rectangular solid in
the moduli space.  The main mathematical
feature of the paper is a result which
gives a lower bound on the energy of any
configuration in a block based on the
potentials of the configurations corresponding
to the vertices, and an error term.   The
Energy Theorem (Theorem \ref{ENERGY}) is our
main result along these lines.  The Energy
Theorem is related to the same
kind of result we had in [{\bf S1\/}]
but the result here is cleaner because
we are dealing with polynomial expressions.
We state the Energy Theorem specifically
for the energy potentials $G_k$ mentioned
above, but in principle it should work
much more broadly.

It is hard to over-emphasize the importance
of having an efficient error term in the
Energy Theorem.  This makes the difference
between a feasible calculation and one which
would outlast the universe.  I worked quite
hard on making the Energy Theorem simple,
efficient, and pretty sharp.
Another nice feature of the Energy Theorem is
that every quantity involved in it is a rational
function of the vertices of the block.  Thus,
at least in principle, we could run all our
computer programs using exact integer arithmetic.
Such integer calculations seemed too slow, however.
So, I implemented the calculations using
interval arithmetic.  Since the terms in the
Energy Theorem are all rational, our
calculations only involve the
operations plus, minus, times, divide.  
Such operations are contolled by the IEEE
standards, so that interval arithmetic is
easy to implement.  In
contrast, if we had to rigorously
bound expressions involving the power
function  \footnote{Towards the end
of the monograph, in \S \ref{cluge},
we are forced to confront the power function
directly for a very limited calculation,
but we will reduce everything
to the basic arithmetic operations in this case.}
we would be in deeper and murkier 
waters computationally. 

So far we have discussed one energy function at a time,
but we are interested in a $1$-parameter
family of power laws and we can only run
our program finitely many times.  This is
where Tumanov's observation comes in.
For instance, if we just run our program
on $G_3$ and $G_5$ (which we do) we pick
up all the $R_s$ functions for
$s \in (-2,0) \cup (0,2]$.  After a lot
of experimenting I found variants of
Tumanov's result which cover the
rest of the exponent range $[0,15+25/512]$.
Our main result along these lines, which
is really a compilation of results, is
the Forcing Lemma. See \S \ref{forcing}.

Since $\shin<15+25/512$ something must go wrong
with our approach above.  What happens is
the TBP is not quite the minimizer with
respect to one of the magic energy functions,
\begin{equation}
G_{10}^{\sharp}=G_{10}+13 G_5+68 G_2.
\end{equation}
However
we can still eliminate most of the configuration
space.  What remains is a small region {\bf SMALL\/}
of mystery
configurations. See \ref{bigsmall}.

  We deal with the configurations in
{\bf SMALL\/} using a symmetrization construction
that, so to speak, is neutral with respect
to the TBP/FP dichotomy.  The symmetrization
construction, which only decreases energy in a small
but critical range, somehow makes the TBP and the FPs
work together rather than as competitors.
We explain the construction in
\S \ref{symm0} and elaborate on the
philosophical reason why it works.

\section{Computer Code}

The computer code involved in this paper is
publicly available by download from my website,
\newline
\newline
{\bf http://www.math.brown.edu/$\sim$res/Java/TBP.tar\/}
\newline
\newline
All the algorithms are implemented in Java.  We also
have a small amount of Mathematica code needed for
the analysis of the Hessians done in \S \ref{local},
and also for the handful of evaluations done in
\S \ref{finalXX}. In both cases, we use Mathematica
mostly for the sake of convenience, and
the calculations are so straightforward and precisely
described that I didn't see the need to implement them in Java.
The reader who knows Mathematica or (say) Matlab or Sage
could re-implement these calculations in about an hour.

The (non-Mathematica) programs come with
graphical user interfaces, and each one has built
in documentation and debuggers. The interfaces
let the user watch the programs in action,
and check in many ways that they are operating
correctly.  The computer programs
are distributed in $4$ directories, corresponding
to the $4$ main phases of the proof. 

 The overall computer program is large and
complicated, but mostly because of the elaborate
graphical user interfaces.  These interfaces
are extremely useful for debugging purposes,
and also they give a good feel for the
calculations.  However, they make the code
harder to inspect.  To remedy this situation,
I have made separate directories which just
contain the code which goes into the computational
tests.  In these directories all the superfluous
code associated to the interfaces has been
stripped away.  The reader can launch the tests
either from the interface or from a command line
in these stripped-down directories.  The tests
do the same thing in each case, though with
the graphical interface there is more room
for varying the conditions of the test
(away from what is formally needed) for the
purposes of experimentation and control.

I think that a competent programmer would be able
to recreate all the code that goes into the
tests, either by reading the descriptions in
this monograph or else by inspecting the code
itself.

\newpage

\chapter{Outline of the Proof}

In this chapter we outline the proof of
the main results.  The proof has
$4$ main ingredients,
Under-Approximation, Divide and Conquer,
Symmetrization, and Endgame.
We discuss these ingredients in turn.
To avoid irritating trivialities,
we consider two configurations the same
if they are isometric.  Given the function $f$
and a $5$-point configuration $X$ we
often write $f(X)$ in place
of ${\cal E\/}_f(X)$ for ease of notation.

\section{Under-Approximation}
\label{forcing}

Let $I \subset \R$ denote an interval, which we
think of an interval of power law exponents.
We say that a triple $(\Gamma_2,\Gamma_3,\Gamma_4)$ of
potentials is {\it forcing\/} 
on the interval $I$ if the following
implication holds for any configuration $5$-point
configuration $X$:
\begin{equation}
\Gamma_i(X)>\Gamma_i(\T) \hskip 10 pt
i=2,3,4 \hskip 15 pt 
\Longrightarrow \hskip 15 pt R_s(X)>R_s(\T) \hskip 10 pt
\forall s \in I.
\end{equation}

Recall that $G_k(r)=(4-r^2)^k$.
Tumanov's observation is that $(G_2,G_3,G_5)$ is
forcing on $(-2,0)$ and $(0,2]$.  Tumanov does not
supply a proof for his observation, but we will prove
a more extensive related result.  Consider the following
functions:
\begin{equation}
\matrix{
G_5^{\flat}& = &G_5-25 G_1 \cr
G_{10}^{\#}& =  &G_{10}+13 G_5 + 68 G_2 \cr
G_{10}^{\#\#}& = &G_{10}+28 G_5 + 102 G_2}
\end{equation}
In \S 3-6 we prove the following result.

\begin{lemma}[Forcing]
\label{Tumanov}
The following is true.
\begin{enumerate}
\item $(G_2,G_4,G_6)$ is forcing on $(0,6]$.
\item $(G_2,G_5,G_{10}^{\#\#})$ is forcing on $[6,13]$.
\item  $(G_2,G_5^{\flat},G_{10}^{\#})$ is forcing on 
$[13,15.05]$.
\item $(G_2,G_3,G_5)$ is forcing on $(-2,0)$.
\end{enumerate}
\end{lemma}

\noindent
{\bf Remark:\/} \newline
(i)
We list the negative case last because it
is something of a distraction.  The
triple $(G_2,G_3,G_5)$ is also forcing on
$(0,2]$, as Tumanov observed, but we ignore
this case because it is redundant. \newline
(ii)
These kinds of under-approximation results,
in some form, are used in
many papers on energy minimization.  See
[{\bf BDHSS\/}] for very general ideas like this.

\section{Divide and Conquer}
\label{bigsmall}

Using the divide-and-conquer algorithm
mentioned in the introduction, together with
a local analysis of the Hessian, we give a
rigorous, computer-assisted proof of
the following result.
\begin{theorem}[Big]
The TBP is the unique global minimizer w.r.t.
$$G_3,G_4,G_5^{\flat},G_6,G_{10}^{\#\#}.$$
\end{theorem}

We note that $G_5$ is a positive combination of
$G_1$ and $G_5^{\flat}$, so $G_5$ also satisfies
the conclusion of the Big Theorem.
The Big Theorem fails for $G_7,G_8,...$.
This failure is what led me to search for more complicated
combinations like $G_{10}^{\#}$.

Combining the Big Theorem and the Forcing Lemma, we get
\begin{corollary}
\label{most}
The TBP is the unique global minimizer w.r.t. $R_s$
for all $s \in (-2,0) \cup (0,13]$.
\end{corollary}

Corollary \ref{most} contains
Theorem \ref{aux2}.  For
Theorem \ref{aux1} we just have to deal with the
exponent interval $(13,15+25/512]$.  The TBP is not the
global minimizer for $G_{10}^{\#}$, but we
can still squeeze information out
of this function. We work in stereographic
coordinates, as in Equation \ref{stereo}.
Precisely, we have the correspondence:
\begin{equation}
V_0,V_1,V_2,V_3,(0,0,1) \in S^2 \hskip 6 pt
\Longleftrightarrow \hskip 6 pt
p_0,p_1,p_2,p_3 \in \R^2, \hskip 20 pt
p_k=\Sigma(V_k).
\end{equation}
We write $p_k=(p_{k1},p_{k2})$.  We always normalize
so that $p_0$ lies in the positive $x$-axis and
$\|p_0\| \geq \|p_k\|$ for $k=1,2,3$.
Let ${\bf SMALL\/}$ denote those $5$-point configurations
which are represented by $4$-tuples $\{p_0,p_1,p_2,p_3\}$
such that 
\begin{enumerate}
\item $\|p_0\| \geq \|p_k\|$ for $k=1,2,3$.
\item $512 p_0 \in [433,498] \times [0,0]$.
\item $512 p_1 \in [-16,16] \times [-464,-349]$.
\item $512 p_2 \in [-498,-400] \times [0,24]$.
\item $512 p_3 \in [-16,16] \times [349,364]$.
\end{enumerate}
This domain is pretty tight.
We tried hard to get as far away from the TBP
configuration as possible.

\begin{center}
\resizebox{!}{3.2in}{\includegraphics{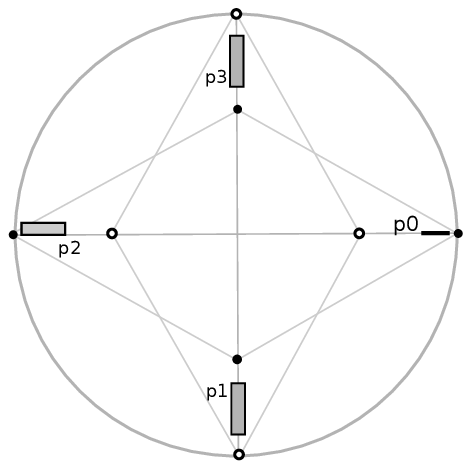}}
\newline
{\bf Figure 2.1:\/} The sets defining {\bf SMALL\/}.
\end{center}

Figure 2.1 shows a picture of the sets
corresponding to the definition of {\bf SMALL\/}.
The grey circle is the unit circle.
Note that $\T \not \in {\bf SMALL\/}$.
The $4$ black dots and the $4$ white dots are the two
nearby normalized TBP configurations. The inner dots
are $\sqrt 3/3$ units from the origin.

\begin{theorem}[Small]
\label{small}
If $G_{10}^{\#}(X) \leq G_{10}^{\#}(\T)$ then
either $X=\T$ or
$X \in {\bf SMALL\/}$.
\end{theorem}

Combining this result with the Forcing Lemma, we get
\begin{corollary}
Let $X$ be any $5$-point configuration.
Let $$s \in \bigg[13,15+\frac{25}{512}\bigg].$$
If $R_s(X) \leq R_s(\T)$ then either $X=\T$ or
$X \in {\bf SMALL\/}$.
\end{corollary}

We cannot plot the $7$-dimensional region {\bf SMALL\/},
but we can plot a canonical slice. 
Let {\bf K4\/} denote the
$2$-dimensional set of configurations
whose stereographic
projections are rhombi with points
in the coordinate axes.
(Here {\bf K4\/} stands for ``Klein $4$ symmetry''.) 
Let
\begin{equation}
{\bf SMALL4\/}={\bf K4\/} \cap {\bf SMALL\/}.
\end{equation}
This is the set of configurations
$(p_0,p_1,p_2,p_3)$ such that
$$
-p_2=p_0=(x,0), \hskip 10 pt
-p_1=p_3=(0,y), \hskip 10 pt
x \geq y, \hskip 10 pt \hskip 10 pt
x,y \in \bigg[\frac{348}{512},\frac{495}{512}\bigg].
$$
The point $(x,y)=(1,1/\sqrt 3)$, outside our slice,
represents the TBP.  The region
${\bf SMALL4\/}$
is the unshaded region below the gray diagonal in
Figure 2.2.  Figure 2.2 shows a plot of
$G_{10}^{\#}$ as a function of $x$ and $y$.
The darkly shaded region is where 
$G_{10}^{\#}(x,y) \leq G_{10}^{\#}(\T)$.
Note that this dark region is contained in the
union of ${\bf SMALL4\/}$ and its
reflection across the diagonal.  This plot doesn't prove anything
but it serves as a sanity check on our calculations.

\begin{center}
\resizebox{!}{3.5in}{\includegraphics{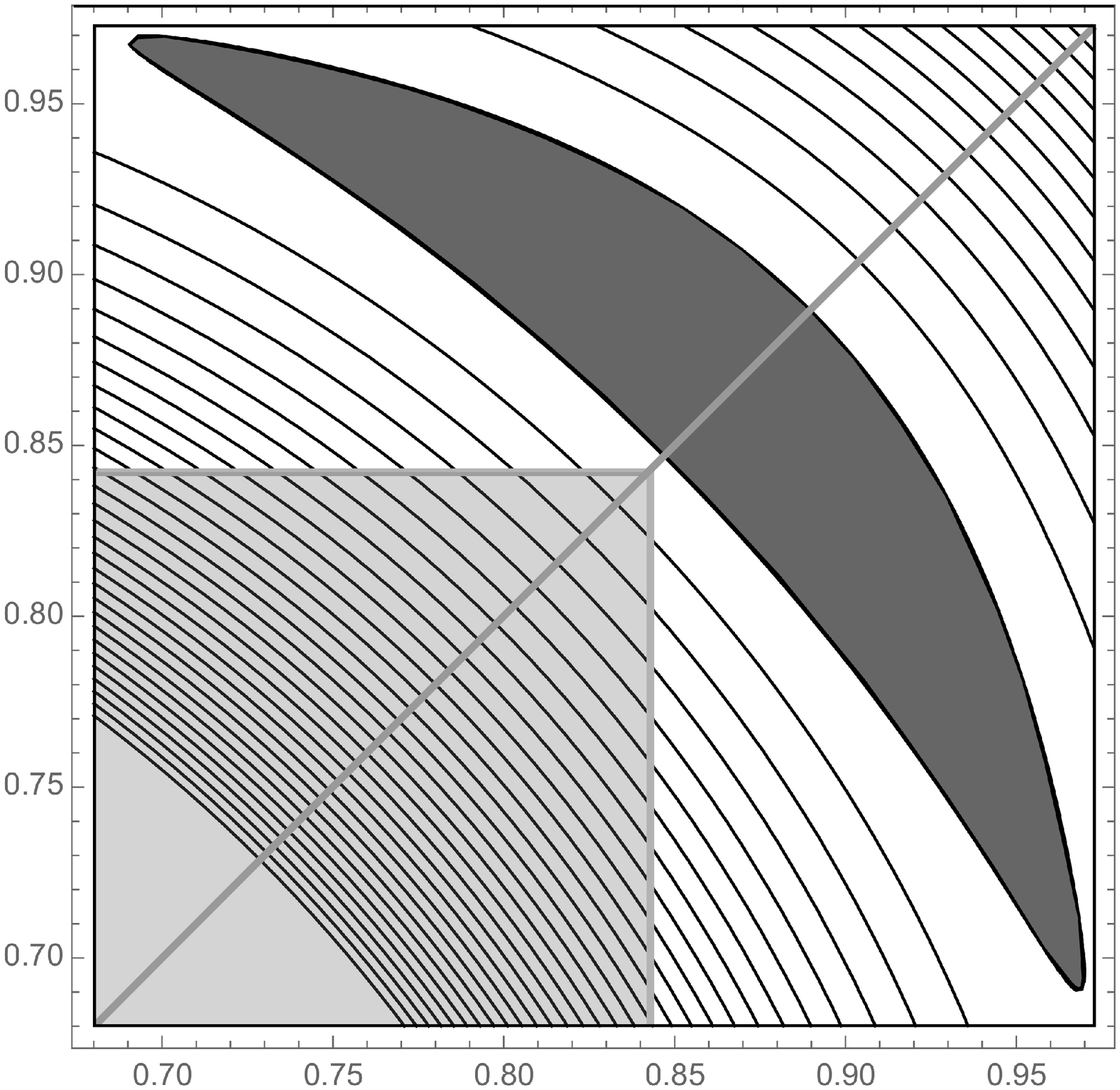}}
\newline
{\bf Figure 2.2:\/} A plot inside {\bf SMALL4\/}.
\end{center}

It remains to deal with configurations in
{\bf SMALL\/}.  We take an approach which
crucially involves the slice {\bf SMALL4\/}.

\section{Symmetrization}
\label{symm0}

Here is the heart of the matter:
We will describe an energy decreasing retraction
from ${\bf SMALL\/}$ into ${\bf K4\/}$.
For some
other choices of the number of points -- e.g. $6$ points -- there
is a configuration -- e.g. the octahedron -- whose global symmetry
dominates over all other configurations.  For $5$ points
there are two kinds of fairly symmetric configurations,
the TBP and FPs, whose
influences compete and weaken each other.
However, if we align these configurations so that
parts of their symmetry groups coincide -- i.e., both
configurations are members of {\bf K4\/} -- then this aligned
symmetry dominates everything nearby.  The existence
of an energy decreasing symmetrization is a precise
formulation of this idea.  Philosophically speaking,
we make progress through the critical region by
aligning the competitors and making them unite their
forces.

We consider the following map from
{\bf SMALL\/} to
{\bf K4\/}.
We start with the configuration $X$ having points
$(p_1,p_2,p_3,p_4) \in {\bf SMALL\/}$.
We let $$(p_1',p_2',p_3',p_4')$$ be the
configuration which is obtained by
rotating $X$ about the origin so that
$p_0'$ and $p_2'$ lie on the same horizontal
line, with $p_0'$ lying on the right.
This is just a very slight rotation.
We then define

\begin{equation}
-p_2^*=p_0^*=\bigg(\frac{p'_{01}-p'_{21}}{2},0\bigg),
\hskip 30 pt
-p_3^*=p_1^*=\bigg(0,\frac{p'_{12}-p'_{32}}{2}\bigg)
\end{equation}

The points $(p_0^*,p_1^*,p_2^*,p_3^*)$ define
the symmetrized configuration $X^*$.
After a brutal struggle, I was finally able
to prove the following result:

\begin{lemma}[Symmetrization]
Let $s \in [12,15+25/512]$ and
suppose that $X \in {\bf SMALL\/}$. Then
$R_s(X^*) \leq R_s(X)$ with
equality iff $X=X^*$.
\end{lemma}

Combining the Symmetrization Lemma and the
Small Theorem, we get the following result:

\begin{lemma}[Dimension Reduction]
Let $s \in [13,15+25/512]$.  Suppose that
$R_s(X) \leq R_s(\T)$ and
$X \not = \T$. Then $X \in {\bf SMALL4\/}$.
\end{lemma}

The Dimension Reduction Lemma practically finishes the
proofs of the Main Theorems.  It leaves us with the
exploration of a $2$-dimensional rectangle in
the configuration space.
Before we get to the endgame, I want to discuss
the Symmetrization Lemma in some detail.

One strange thing is that the map $X \to X^*$ is not
clearly related to spherical geometry.
Rather, it is a linear projection with respect to
the stereographic
coordinates we impose on the moduli space.
On the other hand, the beautiful 
formula in Equation \ref{beauty}
suggests (to me) that there is something
deeper going on.
I found the map $X \to X^*$ experimentally,
after trying many
other symmetrization schemes, some related to
the geometry of $S^2$ and some not.  You might
say that the map defined above is the winner
of a Darwinian competition.

Experimentally, the conclusion of the
Symmetrization Lemma seems to be more
robust than our result suggests.  It seems
to hold in a domain somewhat larger than
{\bf SMALL\/} and for all $s \geq 2$.
However, our proof is quite delicate and uses
the precise domain {\bf SMALL\/} crucially.

Many ideas go into the proof of the
Symmetrization Lemma, some geometric, some
analytic, and some algebraic.  I'd say that
the key idea is to use algebra to implement
some geometric ideas.  For instance, one step
in proving a geometric inequality about triangles
involves showing that some
polynomials in $5$ variables, having about $3000$ terms each,
are positive on on unit cube in $R^5$.
These functions turn out to be positive dominant
in the sense of \S \ref{posdom0}, and so we have
a short positivity certificate available.
Naturally, I would prefer a clean and
simple proof of the Symmetrization Lemma, but
the crazy argument I give is the best I can do
right now.

\section{Endgame}
\label{endgame}

The Dimension Reduction Lemma
makes the endgame fairly routine.
What remains is the exploration of
the $2$-dimensional rectangle
{\bf SMALL4\/}.  It is useful to
introduce a much smaller set.
Define 
\begin{equation}
I=(a,b), \hskip 30 pt
a=\frac{55}{64}, \hskip 40 pt
b=\frac{56}{64}.
\end{equation}
Let {\bf TINY4\/} be denote the square consisting
of configurations such that
\begin{equation}
-p_2,p_0 \in I \times \{0\}, \hskip 30 pt
-p_1,p_3 \in \{0\} \times I.
\end{equation}

We capture $\shin$ in a tiny
parameter interval $[\alpha,\beta]$ where
\begin{equation}
\alpha=15+\frac{24}{512}, \hskip 30 pt
\beta=15+\frac{25}{512}.
\end{equation}
We prove the following result by calculations
which are similar to those which prove
the Big Theorem and the Small Theorem:

\begin{lemma}
\label{most0}
Suppose that one of the following is true:
\begin{enumerate}
\item $s \in [13,\alpha]$ and  
$X \in {\bf SMALL4\/}$.
\item $s \in [\alpha,\beta]$ and
$X \in {\bf SMALL4\/}-{\bf TINY4\/}$.
\end{enumerate}
Then $R_s(X)>R_s(\T)$.
\end{lemma}

Combining Lemma \ref{most0} with the Dimension
Reduction Lemma, we see that 
the TBP must be the unique global minimizer w.r.t.
$R_s$ when $s \in [13,\alpha]$ and that the global
minimizer when $s \in [\alpha,\beta]$ is either
the TBP or some member of {\bf TINY4\/}.

To finish the proof, we use another symmetrization
trick.  The plane is foliated by parabolas having
equations of the form
\begin{equation}
\gamma(t)=(t_0+t- \lambda t^2,t_0-t-\lambda t^2),
\hskip 30 pt \lambda=2.
\end{equation}
This foliation induces a natural retraction from
$\R^2$ to the main diagonal.  We map each
point on $\gamma(t)$ to $\gamma(0)$.  When we
identify {\bf TINY4\/} with a subset of $\R^2$
by taking the coordinates $(p_{01},p_{32})$, we
get a retraction from {\bf TINY4\/} onto the
subset consisting of pyramids with square base.
We call this map the {\it parabolic retraction\/}.

\begin{lemma}
\label{tiny}
Let $s \in [\alpha,\beta]$ and
$X \in {\bf TINY4\/}$.
Let $X^*$ denote the image of $X$ under the
parabolic retraction.  Then $R_s(X^*) \leq R_s(X)$,
with equality iff $X=X^*$.
\end{lemma}

\noindent
{\bf Remark:\/}
Lemma \ref{tiny} is quite delicate. The analogous result
fails if we replace $\lambda=2$ by $\lambda=1$ or
$\lambda=4$.  I found the parabolic retraction after
a lot of trial and error.
\newline

Lemma \ref{tiny}
finishes the proof of Theorem \ref{aux1}.
We actually know more at this point:
\begin{enumerate}
\item The TBP is the unique minimizer w.r.t
$R_s$ when $s \in (0,\alpha]$. 
\item When $s \in [\alpha,\beta]$ any minimizer
w.r.t. $R_s$ is either the TBP or an FP.
\item When $s=\beta$ some FP is the unique
global minimizer.  
\end{enumerate}
The third statement follows from
Lemma \ref{small}, 
Lemma \ref{tiny}, and a check that a
certain FP has lower energy than the TBP
at the parameter $\beta$.
In \S \ref{finalXX} we deduce the Main
Theorem from these two three statements and
some straightforward tricks with calculus.

\section{Computing the Constant}
\label{cc}

We close this chapter with a sketch of how
an ideal computer would compute $\shin$.
The $R_s$ energy of the TBP is
\begin{equation}
R(s)=6 \times 2^{-s/2} + 3 \times 3^{-s/2} + 4^{-s/2}.
\end{equation}
Holding $s$ fixed, let $\phi(t)$ be the
$R_s$ energy of the FP with vertices
\begin{equation}
(\pm \sqrt{1-t^2},0,t), \hskip 30 pt (0,\pm \sqrt{1-t^2},t),
\hskip 30 pt (0,0,1).
\end{equation}
We have
\begin{equation}
\phi(t)=2\bigg(4-4t^2\bigg)^{-s/2} +
4\bigg(2-2t\bigg)^{-s/2} + 
4\bigg(2-2t^2\bigg)^{-s/2}.
\end{equation}
Let $I=(-3/16,-1/16)$.
One check that for each
$s \in [15,15.1]$, that:
\begin{itemize}
\item $\phi$ has it minimum in $I$.
\item $|\phi'(t)|<1/10$ in $I$. That is, $\phi$ is $(1/10)$-Lipschitz.
\item $\phi''(t)>0$ on $I$. That is, $\phi$ is convex on $I$.
\end{itemize}
Since $\phi$ is convex in $I$, one can compute its minimum
value $Q(s)$ to arbitrary precision using a simple iterative
algorithm. Namely, \begin{itemize}
\item Let $a_{0,0}=-3/16$ and $a_{0,1}=-1/16$.
\item Let $b_k=(a_{k,0}+a_{k,1})/2$.
\item If $\phi'(b_k)>0$ let $a_{k+1,0}=a_{k,0}$
and $a_{k+1,1}=b_k$.
\item If $\phi'(b_k)<0$ let $a_{k+1,0}=b_k$ and
$a_{k+1,1}=a_{k,1}$.
\end{itemize}
If $\phi'(b_k)=0$ for some $k$ then
$Q(s)=\phi(b_k)$.  Otherwise 
$\{b_k\}$ converges and $Q(s)=\lim \phi(b_k)$.

The above iteration is a step in a similar
algorithm which computes $\shin$.  First I
describe the {\it comparison step\/}.
Assuming that $s \not = \shin$ then one of two
things is true:
\begin{enumerate}
\item There is some $k$ for which $\phi(b_k)<R(s)$.
In this case $s>\shin$.
\item There is some $k$ for which
$$
\phi(b_k)-R(s)>\frac{a_{k,1}-a_{k,0}}{10}.
$$
In this case $Q(s)>R(s)$ and $s<\shin$.
\end{enumerate}

Here is the algorithm which computes $\shin$.
\begin{itemize}
\item Set $a_{0,0}=15+(1/32)$ and $a_{0,1}=15+(3/32)$.
\item Let $b_k=(a_{k,0}+a_{k,1})/2$.
\item If $b_k<\shin$ let $a_{k+1,0}=b_k$ and $a_{k+1,1}=a_{k,1}$.
\item If $b_k>\shin$ let $a_{k+1,0}=a_{k,0}$ and $a_{k+1,1}=b_k$.
\end{itemize}
Then $a_{k,0} \to \shin$ from below and
$a_{k,1} \to \shin$ from above.  There is
one {\it caveat\/}:
If $\shin=b_k$ for some $k$, the algorithm will 
simply run forever at the comparison step and we 
will have detected an unspeakably improbable
coincidence.

A rigorous implementation of this algorithm would be tedious,
on account of the non-integer power
operations, but a
lazy floating-point implementation  yields
$$\shin=15.0480773927797...$$
I dropped off the last few digits to be sure about
floating point error.

\newpage

\part{Under-Approximation}
\chapter{Outline of the Forcing Lemma}

\section{A Four Chapter Plan}

The next $4$ chapters are devoted to the
proof of the Forcing Lemma.  In this
chapter we will explain the ideas
in the proof and leave certain
details unproved.  In place of proofs,
we will show computer
plots.  In \S 4 and \S 5 we develop
the machinery 
needed to take care of the details.
Once the machinery is assembled we
will finish the proof of the Forcing
Lemma in \S 6.
We split things up this way so that
the reader can first get a feel for the
Forcing Lemma through numerical evidence
and pictures.  Perhaps the reader will
appreciate the need for the technical
machinery we use to convert the
numerical evidence into rigorous proofs.
(If not, then I would welcome a shorter
proof!)

\section{The Regular Tetrahedron}
\label{regtet}

We first explain a toy version of the
Forcing Lemma, which works for the regular
tetrahedron.   I learned this from the paper 
[{\bf CK\/}], and I think that the idea
goes back to Yudin [{\bf Y\/}].  
It is the source for many results
about energy minimization, including the Forcing Lemma.

We set $\Gamma_0=G_0$, the function that is identically one,
and $\Gamma_1=G_1$.  Again, we have
\begin{equation}
G_k(r)=(4-r^2)^k.
\end{equation}

\begin{lemma}
\label{triv}
A configuration of $N$ points on $S^2$ has
minimal $\Gamma_1$ energy if and only if
the center of mass of the configuration
is the origin.
\end{lemma}

\startproof
For vectors $V,W \in S^2$ we have
\begin{equation}
\Gamma_1(V,W)=4-(V-W) \cdot (V-W) = 2+2 V \cdot W.
\end{equation}
The total energy of a $N$-point configuration
$\{V_i\}$ is
$$\sum_{i<j} 2+2 V_i \cdot V_j =N(N-2) + \sum_{i,j} V_i \cdot V_j = N(N-2) + \bigg\|\sum_i V_i\bigg \|^2.$$
Obviously, this is minimized if and only if $\sum V_i=0$.
\endproof

\begin{corollary} 
Suppose that $\phi:(0,4] \to \R$ is a convex
decreasing function and $R(r)=\phi(r^2)$.
Then the regular tetrahedron is the unique
minimizer w.r.t. $R$ amongst all $4$-point
configurations.
\end{corollary}

\startproof
The distance between any pair of points in the
regular tetrahedron is $\sqrt{8/3}$.  Since
$\Gamma_1$ is linear in $r^2$ and the graph 
has negative slope,
we can find a combination
$\Gamma=a_0\Gamma_0+a_1\Gamma_1$, with
$a_1>0$, such that 
$\Gamma(r) \leq R(r)$ with
equality if and only if $r=\sqrt{8/3}$.
Supposing that $X$ is some configuration and $X_0$ is the
regular tetrahedron, we have
\begin{equation}
R(X) \geq \Gamma(X)=a_0+a_1\Gamma_1(X) \geq a_0+a_1\Gamma_1(X_0)=
\Gamma(X_0)=R(X_0).
\end{equation}
We are using the fact that $X_0$ is a minimizer for $G_1$,
though not uniquely so.
The first inequality is strict unless $X=X_0$.
Here we set $R(X)={\cal E\/}_R(X)$, etc., for ease of notation.
\endproof

\section{The General Approach}

Now we consider $5$-point configurations.
We know that the TBP is a minimizer for
$\Gamma_0$ and $\Gamma_1$, though not
uniquely so.  Suppose that
$\Gamma_2,\Gamma_3,\Gamma_4$ are $3$ more
functions.  We have in mind the various
triples which arise in the Forcing Lemma.

The distances involved in $\T$ are 
$\sqrt 2,\sqrt 3,\sqrt 4$.  (Writing
$\sqrt 4$ for $2$ makes the equations look nicer.)
Let $R$ be some function defined on $(0,2]$.
Suppose we can find a combination
\begin{equation}
\Gamma=a_0\Gamma_0+...+a_4 \Gamma_4, \hskip 30 pt
a_1,a_2,a_3,a_4 > 0.
\end{equation}
such that
\begin{equation}
\Gamma(x) \leq R(x), \hskip 30 pt
\Gamma(\sqrt m)=f(\sqrt m), \hskip 20 pt m=2,3,4.
\end{equation}
Suppose also that
$\Gamma_j(X)>\Gamma_j(\T)$ for $j=2,3,4$. Then
$R(X)>R(\T)$. The proof is just like what we
did for the tetrahedron:
\begin{equation}
R(X) \geq \Gamma(X)=\sum a_i \Gamma_i(X)>
\sum a_i \Gamma_i(\T)=\Gamma(\T)=R(\T).
\end{equation}
Note that we have not given the minimal hypotheses
needed to get this inequality to work.  

Here is how we find the coefficients $\{a_i\}$.
We impose the $5$ conditions
\begin{itemize}
\item $\Gamma(x)=R(x)$ for $x=\sqrt 2,\sqrt 3,\sqrt 4$.
\item $\Gamma'(x)=R'(x)$ for $x=\sqrt 2,\sqrt 3$.
\end{itemize}
Here $R'=dR/dx$ and $\Gamma'=d\Gamma/dx$.
These $5$ conditions give us $5$ linear
equations in $5$ unknowns.
In the cases described below, the
associated matrix is invertible and there is a unique
solution.   The proof of the Forcing Lemma,
in each case, involves solving the matrix
equation, checking positivity of the coefficients,
and checking the under-approximation property.
In this chapter we will present the solutions
to the matrix equations, and show plots of the
coefficients and sample plots of the under
approximations. 

\section{Case 1}

The triple is $G_2,G_4,G_6$ and
$s \in (0,6]$ and $R=R_s$, where
$R_s(r)=r^{-s}$. We will suppress the
constant $s$ whenever it seems possible to
do it without causing confusion.
When we want to denote the dependence of $\Gamma$ on the
parameter, we write $\Gamma_{(s)}$.  That is,
$\Gamma_{(s)}=\sum a_i(s) \Gamma_i$.
We need to find the coefficients $a_0,a_1,a_2,a_3,a_4$,
and we will also keep track of the quantity

\begin{equation}
\delta=2\Gamma'(2)-2R'(2).
\end{equation}

The solution is given by

{\tiny
$$
\left[\matrix{
a_0 \cr a_1 \cr a_2 \cr a_3 \cr a_4\cr \delta}\right]=
\frac{1}{792}\left[\matrix{
0 & 0 & 792 & 0 & 0 & 0 \cr
792 & 1152 & -1944 & -54 & -288 & 0 \cr 
-1254 & -96 & 1350 & 87 & 376 & 0 \cr
528 & -312 & -216 & -39 & -98 & 0 \cr
-66 & 48 & 18 & 6 & 10 & 0 \cr
-6336 & -9216 & 15552 & 432 & 2304 & 792}\right]
 \left[\matrix{
2^{-s/2}\cr 3^{-s/2}\cr 4^{-s/2} \cr
s2^{-s/2}\cr s3^{-s/2}, \cr s4^{-s/2}}\right].
$$
\/}

The left side of Figure 3.1 shows a graph of
$$80a_1, \hskip 20 pt 200 a_2,
\hskip 20 pt 2000 a_3, \hskip 20 pt
10000 a_4,$$
considered as functions of the exponent $s$.
Here $a_1,a_2,a_3,a_4$ are colored darkest to lightest.
The completely unimportant positive multipliers
are present so that we get a nice picture.
 It turns out that $a_3$ goes
negative between $6$ and $6.1$, so the interval
$(0,6]$ is fairly near to the maximal
interval of positivity. 

\begin{center}
\resizebox{!}{2.5in}{\includegraphics{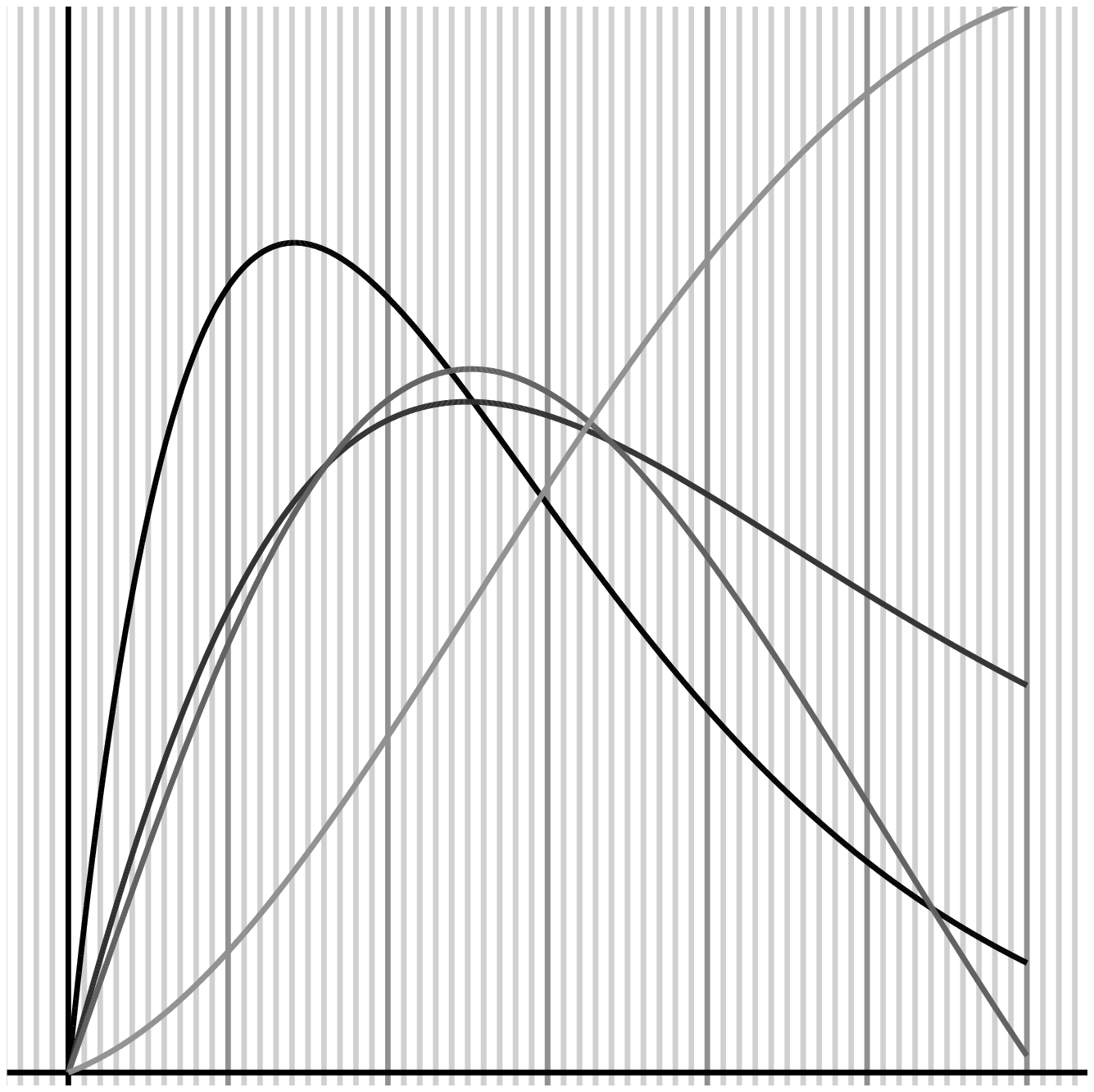}}
\resizebox{!}{2.5in}{\includegraphics{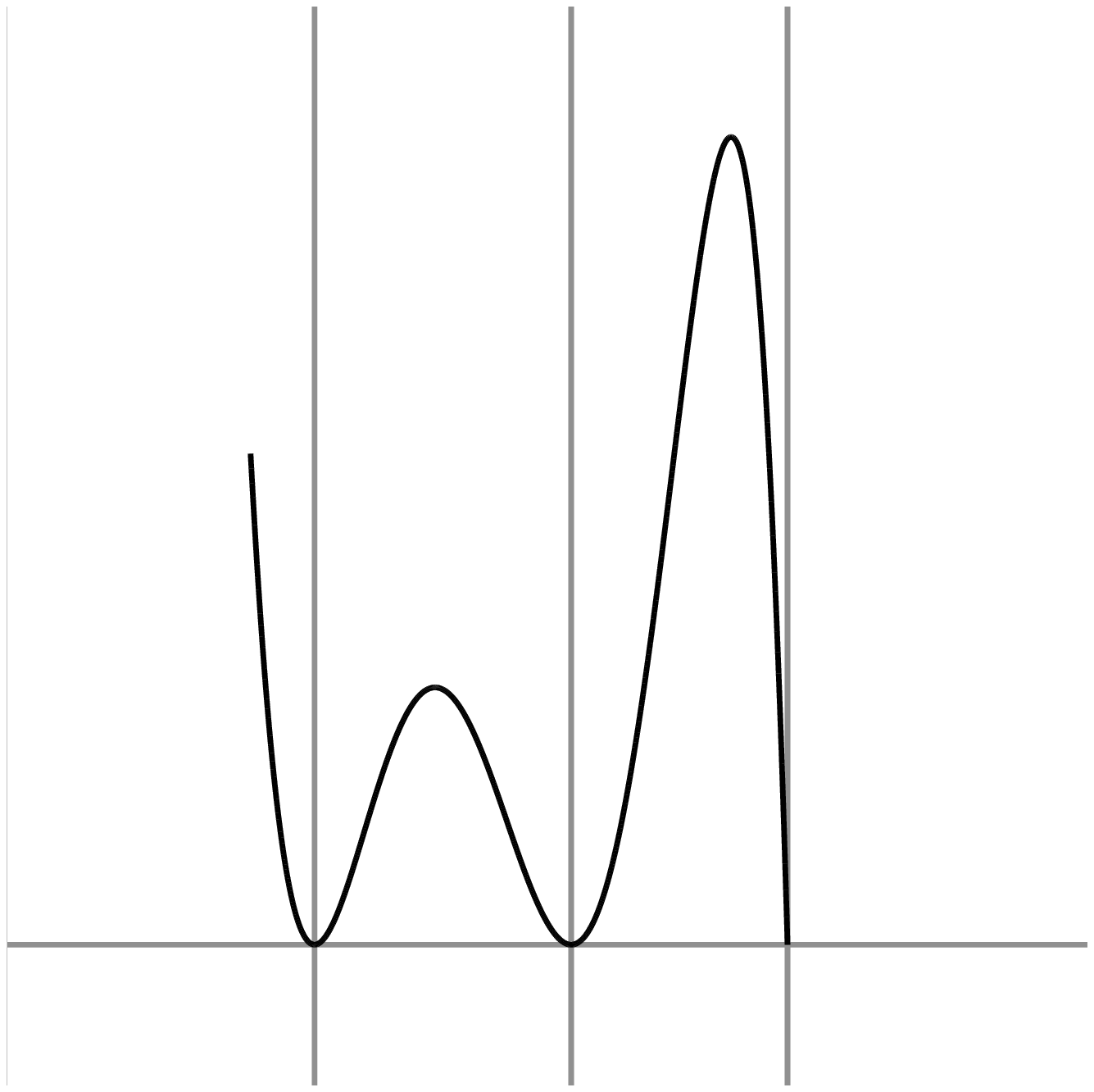}}
\newline
{\bf Figure 3.1:\/} Plots for Case 1.
\end{center}

The right hand of Figure 3.1
shows some positive multiple of
$$H_3(r)=\bigg(1-\frac{\Gamma_{(3)}(r)}{R_3(r)}\bigg)$$
plotted from $r=4/3$ to $r=2$. 
For $r<4/3$ the graph is even more positive,
but it is hard to fit the whole thing into a
nice picture. The horizontal line is the $x$-axis and
the vertical lines are $x=\sqrt 2,\sqrt 3,\sqrt 4$.
The fact that $H_{3} \geq 0$
implies that $\Gamma_{(3)}(r) \leq R_{3}(r)$ for all
$r \in (0,2]$, as desired.
 The reader can use our computer program
and get a much more interactive picture.
In particular, you can adjust the endpoint
of the plot, and see what happens for other values
of $s$, 

 Since
$R_3>0$, the fact that $H_3 \geq 0$ implies that
$\Gamma_{(3)}(r) \leq R_3(r)$ for all $r \in (0,2]$.
The right hand plot looks similar for other
parameters $s \in (0,6]$, as the reader can
see from my program. Note that $H_3$ has negative
slope at $r=2$ and this corresponds to $\delta(3)>0$.

\section{Case 2}

Here the triple is $G_2,G_5,G_{10}^{\#\#}$.
We have

{\tiny
$$
\displaystyle
\left[\matrix{
a_0 \cr a_1 \cr a_2 \cr a_3 \cr \widehat a_4\cr \delta}\right]=
\frac{1}{268536}\left[\matrix{
0& 0& 268536& 0& 0& 0 \cr 
88440& 503040& -591480& -4254& -65728& 0 \cr 
-77586& -249648& 327234& 2361& 65896& 0\cr 
41808& -19440& -22368& -2430& -9076& 0\cr 
-402& 264& 138& 33& 68& 0\cr 
-707520& -4024320& 4731840& 34032& 525824& 268536}\right] 
\left[\matrix{
2^{-s/2}\cr 3^{-s/2}\cr 4^{-s/2} \cr
s2^{-s/2}\cr s3^{-s/2}, \cr s4^{-s/2}}\right].
$$
\/}
Figure 3.2 does for Case 2 what Figure 3.1 does for Case 1.
This time the left hand side plots
$$500 a_1 \hskip 20 pt
500 a_2, \hskip 20 pt
5000 a_3, \hskip 20 pt
500000 a_4.$$
for $s \in [6,13]$.
The coefficients $a_1,a_2,a_3$ go negative for $s$ just 
a tiny bit larger than $13$.  I worked hard to find
the function $\Gamma_4=G_{10}+28G_5+102 G_2$ so that
we could get all the way up to $s=13$.
The right hand side shows a plot of $H_{10}$
from $r=5/4$ to $r=2$.  Since $H_{10}$ has negative
slope at $r=2$ we have $\delta(10)>0$.

\begin{center}
\resizebox{!}{2.5in}{\includegraphics{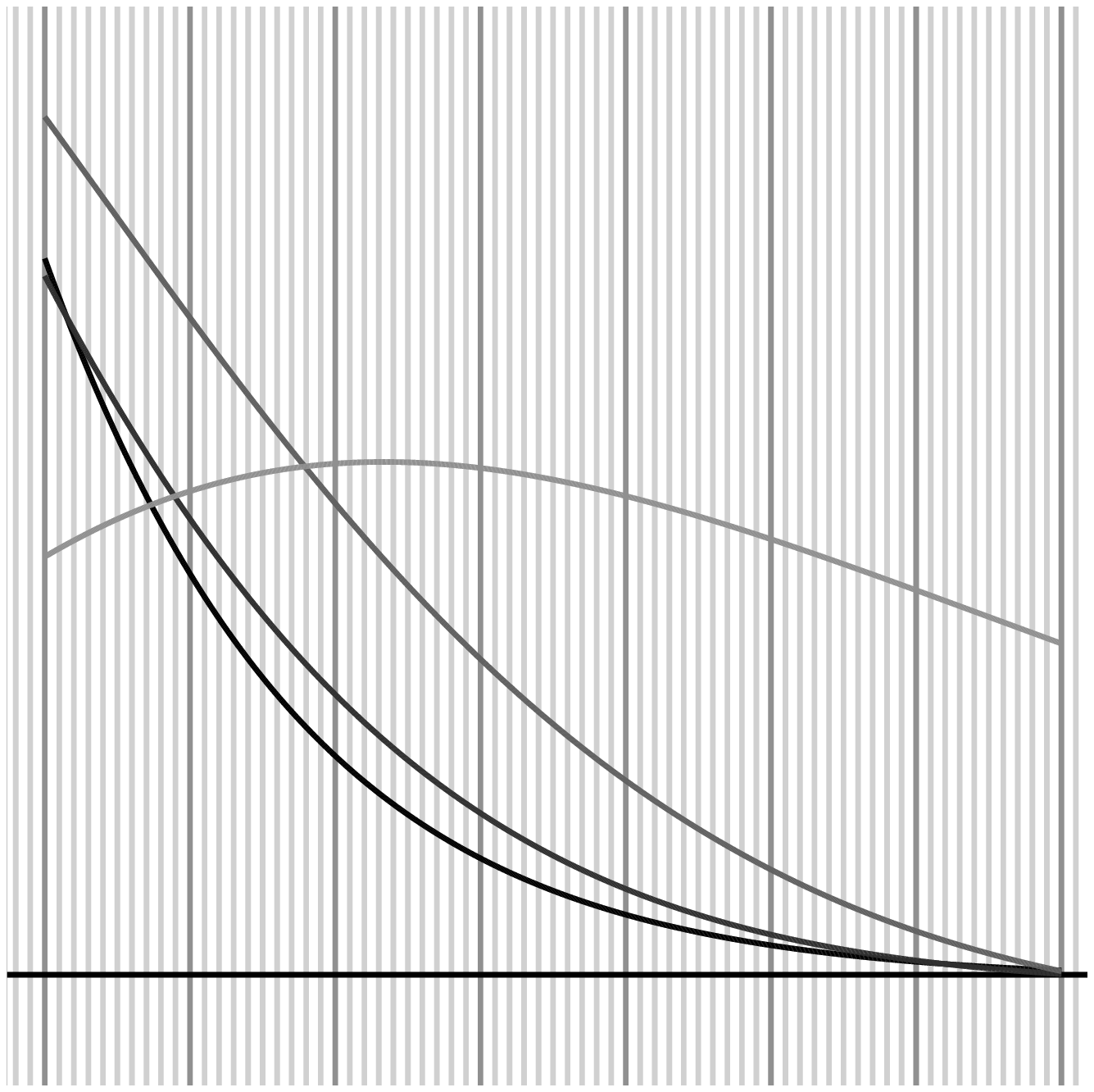}}
\resizebox{!}{2.5in}{\includegraphics{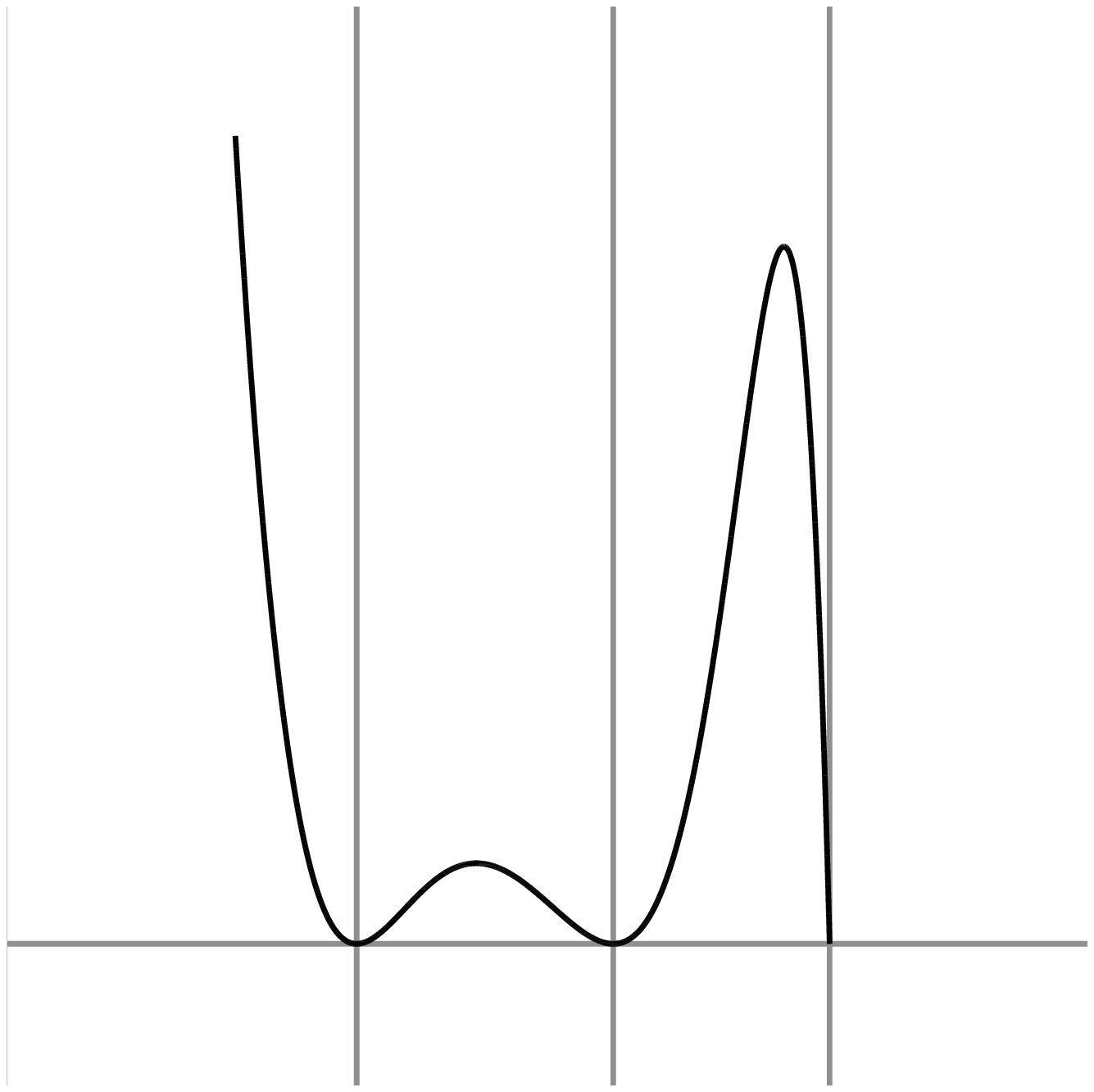}}
\newline
{\bf Figure 3.2:\/} Plot for Case 2.
\end{center}

The interested reader can use our program to
see what these coefficients look like for
parameter value $s<6$. 

\newpage

\section{Case 3}

Here the triple is $G_2,G_5^{\flat},G_{10}^{\#}$.
We have

{\tiny
$$
\displaystyle
\left[\matrix{
a_0 \cr a_1 \cr a_2 \cr a_3 \cr \widehat a_4\cr \delta}\right]=
\frac{1}{268536}
\left[\matrix{
0&0&268536&0&0&0 \cr
947112&131520&-1078632&-50694&-259072&0  \cr
-91254&-240672&331926&3483&68208&0  \cr
35778&-15480&-20298&-1935&-8056&0  \cr
-402&264&138&33&68&0  \cr
174268608&24199680&-198468288&-9327696&-47669248&268536}
\right]
$$ 
\/}
This matrix is quite similar to the one in the previous
case, because we are essentially still taking combinations
of $G_0,G_1,G_2,G_5,G_{10}$.  We are just grouping the
functions differently.
Figure 3.3 does for Case 3 what Figure 3.2 does for Case 2.
This time we plot
$$500 a_1 \hskip 20 pt
15000 a_2, \hskip 20 pt
20000 a_3, \hskip 20 pt
500000 a_4,$$
for $s \in [13,16]$. 
The coefficients $a_1,a_2,a_3$ go negative for $s$ just 
a tiny bit larger than $15.05$.   In particular,
everything up to and including our cutoff of
$5+25/512$ is covered.
The right hand side shows a plot of $H_{14}$ from
$r=5/4$ to $r=2$.  Since $H_{14}$ has negative
slope at $r=2$ we have $\delta(14)>0$.

\begin{center}
\resizebox{!}{2.5in}{\includegraphics{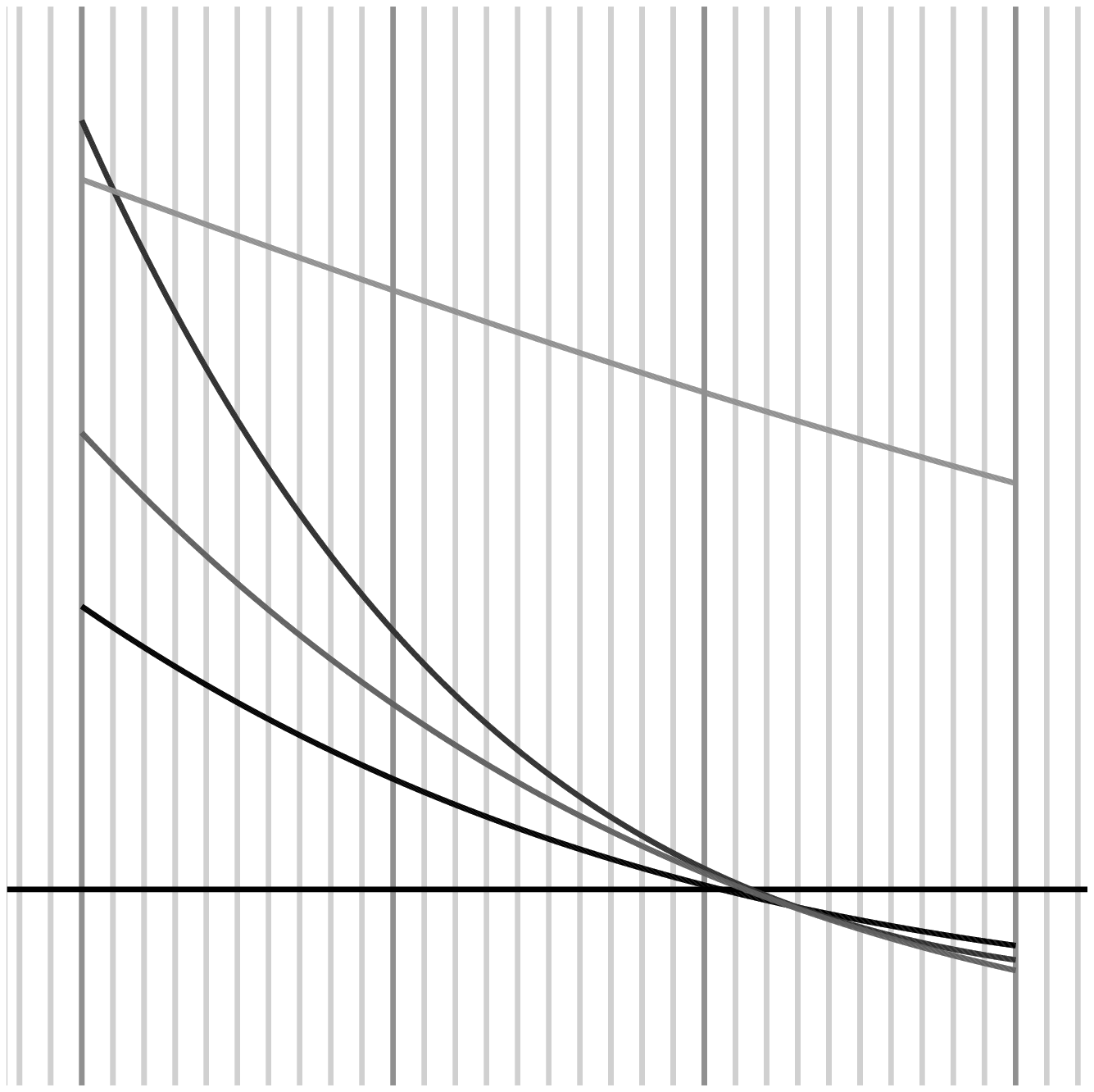}}
\resizebox{!}{2.5in}{\includegraphics{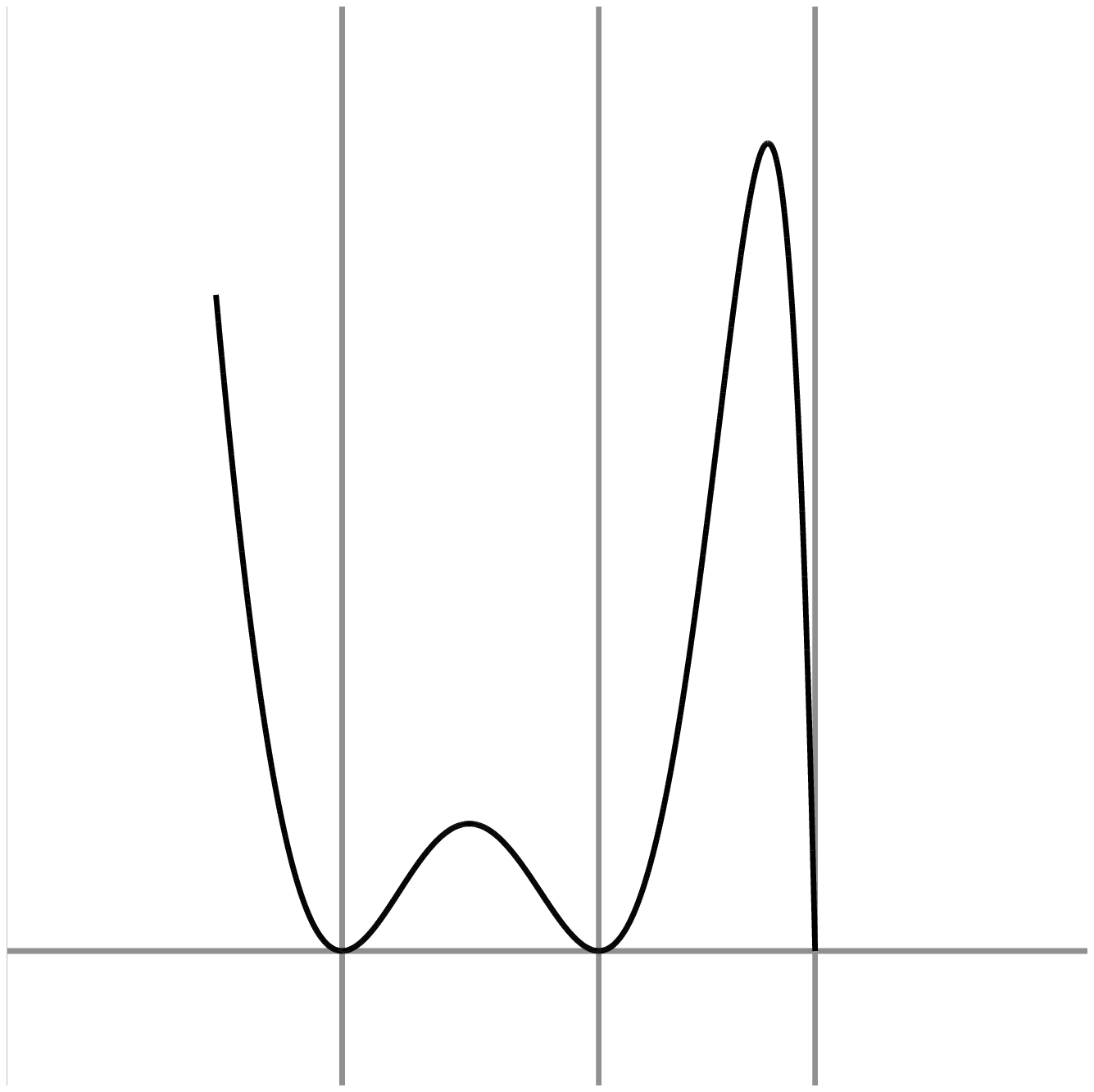}}
\newline
{\bf Figure 3.3:\/} Plot for Case 3.
\end{center}

The interested reader can use our program to
see what these coefficients look like for
parameter value $s<13$. 

\section{Case 4}

Now we consider the negative case.
The triple is $(G_2,G_3,G_5)$.
The solution here is
{\tiny
$$
\left[\matrix{
a_0 \cr a_1 \cr a_2 \cr a_3 \cr a_4\cr \delta}\right]=
\frac{1}{144} \left[\matrix{
0 & 0 & -144 & 0 & 0 & 0 \cr 
-312 & -96 & 408 & 24 & 80 & 0 \cr 
684 & -288 & -396 & -54 & -144 & 0 \cr 
-402 & 264 & 138 & 33 & 68 & 0 \cr 
30 & -24 & -6 & -3 & -4 & 0 \cr 
2496 & 768 & -3264 & -192 & -640 & -144}\right]
 \left[\matrix{
2^{-s/2}\cr 3^{-s/2}\cr 4^{-s/2} \cr
s2^{-s/2}\cr s3^{-s/2}, \cr s4^{-s/2}}\right].
$$
\/}

The left side of Figure 3.4 shows a graph of
$$2a_1, \hskip 20 pt 100 a_2,
\hskip 20 pt 1000 a_3, \hskip 20 pt
10000 a_4,$$

\begin{center}
\resizebox{!}{2.5in}{\includegraphics{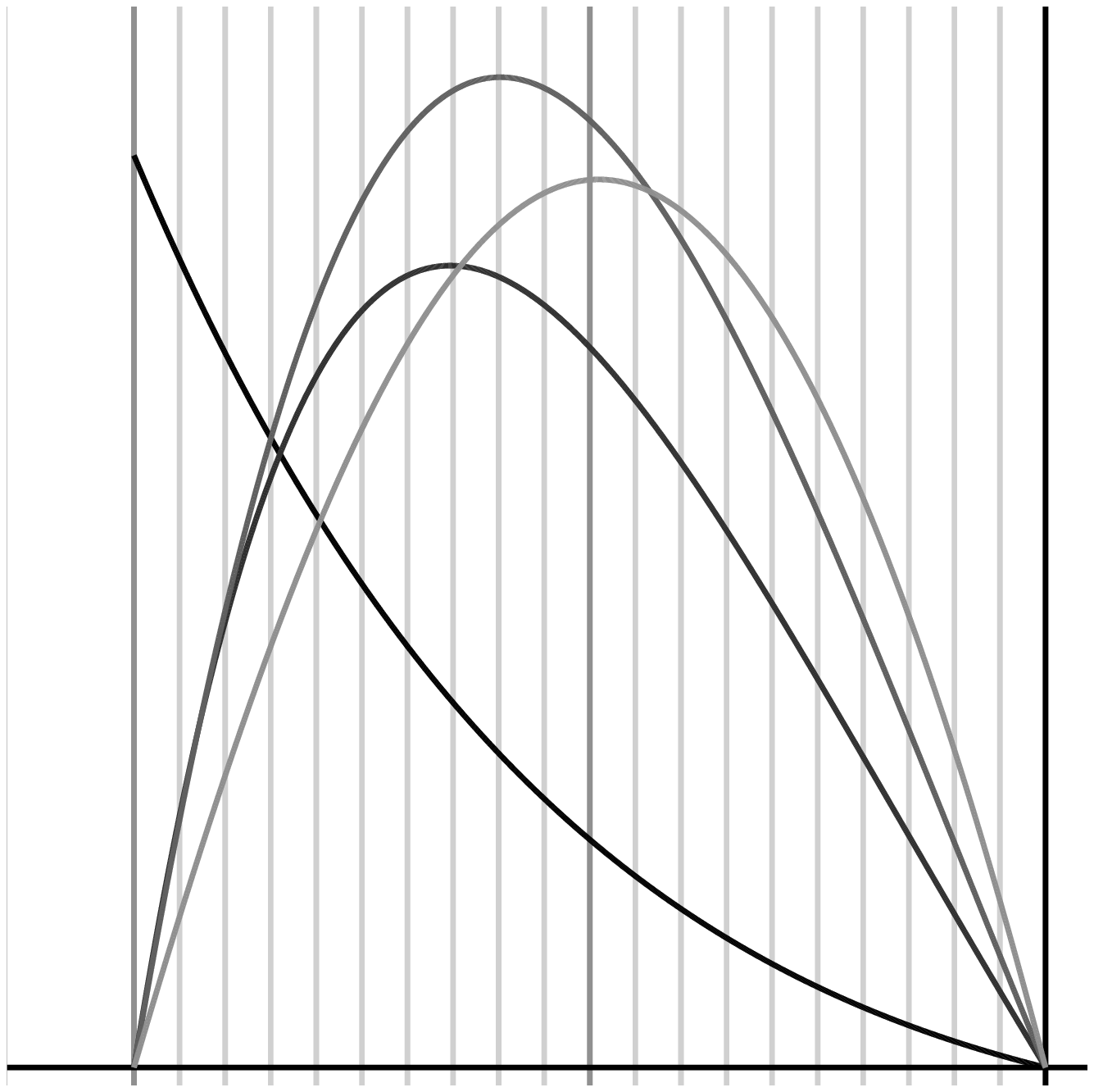}}
\resizebox{!}{2.5in}{\includegraphics{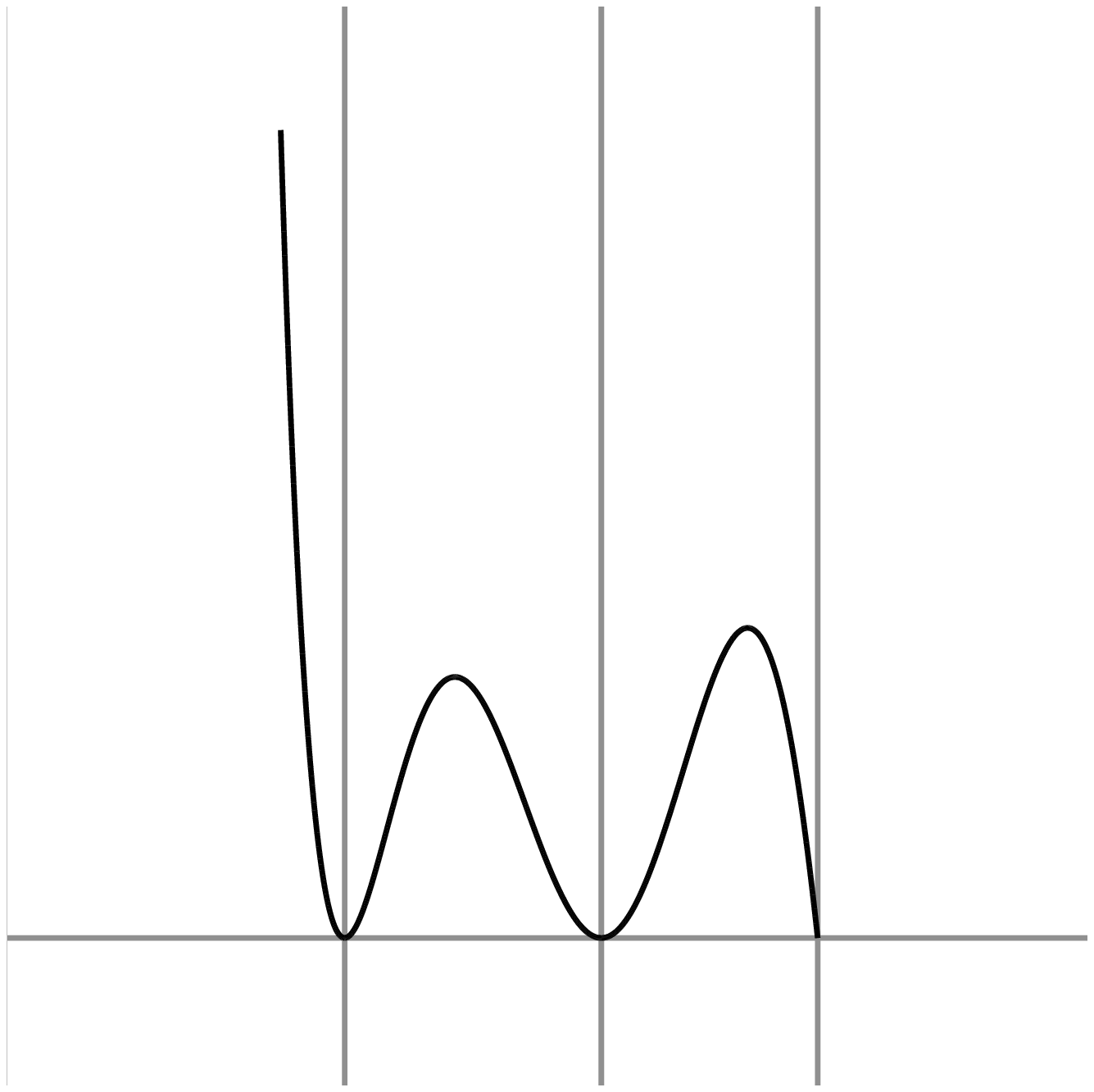}}
\newline
{\bf Figure 3.4:\/} Plots for Case 4.
\end{center}

The right side of 
Figure 3.1 shows a plot of some constant times
\begin{equation}
H_{-1}=\bigg(\frac{\Gamma_{(-1)}(r)}{R_{-1}(r)}-1\bigg)
\end{equation}
Note that $H_{-1}$ has negative slope
at $r=2$, and this corresponds to $\delta(-1)>0$.

\section{Have Fun Playing Around}

I encourage you to have fun playing around
with my computer program.
The program allows you to specify 
$5$ functions of the form
$$A G_a + B G_b + C G_c,$$
and then it will show the above kinds of
plots.  Here
$A,B,C $ are integers and $a,b,c$ are positive integers.
The program also does a rough and ready calculation
to decide where $\T$ has lower energy than any FP,
so you can adjust the combinations so as to get
things which would be useful for proving that
the the TBP is a minimizer for power laws within
a certain exponent range.  

My program will also carry out rigorous positivity proofs
along the lines of what we explain
in the next three chapters.
In this way you can cook up your own versions of
the Forcing Lemma.
For instance, the triple
$$\Gamma_2,\Gamma_4,\Gamma_8+4\Gamma_{4}+32 \Gamma_2$$
is forcing on $[0,11]$.  I found all these crazy
combinations using my program.

\newpage
\chapter{Polynomial Approximations}
\label{polyapx}

This is the first of two chapters
which develop machinery to prove
the positivity claims needed for
the Forcing Lemma -- i.e. that the
functions which appear to be positive
(or non-negative)
from the plots in the previous
chapter really are positive (or non-negative).
The functions involved in the Forcing Lemma
are not polynomials.  Thus, in this chapter
we discuss a method for approximating these
functions with polynomials.  Basically, the
method combines Taylor's Theorem with
the notion of an {\it interval polynomial\/}.

\section{Rational Intervals}
\label{rationalinterval}

We define a {\it rational interval\/} to be
an interval of the form $I=[L,R]$ where 
$L,R \in \Q$ and $L \leq R$.
For each operation $* \in \{+,-,\times\}$ we
define
\begin{equation}
I_1 * I_2=[\min(S),\max(S)],
\hskip 25 pt S=\{L_1*L_2,L_1*R_2,R_1*L_2,R_2*L_2\}.
\end{equation}
These operations are commutative.
The definition is such that
$r_j \in I_j$ for $j=1,2$ then $r_1*r_2 \in I_1*I_2$.
Moreover, $I_1*I_2$ is the minimal
interval with this property.  The minimality
property implies that our laws are
both associative and distributive:
\begin{itemize}
\item $(I_1+I_2) \pm I_3=I_1+(I_2 \pm I_3)$.
\item $I_1 \times (I_2 \pm I_3)=(I_1 \times I_2) \pm (I_1 \times I_3)$.
\end{itemize}
We also can raise a rational interval to a 
nonnegative integer power:
\begin{equation}
I^k=I \times ... \times I \hskip 30 pt {\rm k\ times\/}.
\end{equation}

With the operations above, the set of rational
intervals forms a {\it commutative semi-ring\/}.
The $0$ element is $[0,0]$ and the $1$-element
is $[1,1]$.  Additive inverses do not make sense,
and so we don't get a ring. 
Even though the semi-ring structure is
nice and even beautiful, we technically do not need it.
The main thing the semi-ring structure does is
allow us to write complicated algebraic
expressions without specifying the exact order
of operations needed to evaluate the expressions.
Were some of these laws to fail we would just
pick some order of operations and this would suffice
for our proofs.  Indeed, when we perform interval
arithmetic based on floating point operations we
lose the semi-ring structure but still retain what
we need for our proofs.

\section{Interval Polynomials}

An {\it inteval polynomial\/} is an
expression of the form 
\begin{equation}
I_0+I_1t+ ... + I_n t^n.
\end{equation}
in which each coefficient is an interval
and $t$ is a variable meant to be taken in $[0,1]$.
Given the rules above, interval polynomials
may be added, subtracted or multiplied, in the
obvious way.  The set of interval polynomials
again forms a commutative semi-ring.

We think of an ordinary polynomial as an interval
polynomial, {\it via\/} the identification
\begin{equation}
a_0+...+a_nt^n \hskip 10 pt
\Longrightarrow \hskip 10 pt
[a_0,a_0]+...+[a_n,a_n]t^n.
\end{equation}
(Formally we have an injective semi-ring
morphism from the polynomial ring into the
interval polynomial semi-ring.)
We think of a constant as an interval
polynomial of degree $0$.  Thus, if we have some
expression which appears to involve constants,
ordinary polynomials, and interval polynomials,
we interpret everything in sight as an 
interval polynomial and then perform the
arithmetic operations needed to simplify
the expression.

Let $\cal P$ be the above interval polynomial.
We say that $\cal P$
{\it traps\/} the ordinary
polynomial 
\begin{equation}
C_0+C_1t+...+C_nt^n
\end{equation}
of the same degree
degree if $C_j \in I_j$ for all $j$.
We define the {\it min\/} of an interval polynomial
to be the polynomial whose coefficients are
the left endpoints of the intervals.  We
define the {\it max\/} similarly.   If $\cal P$
is an interval polynomial which traps an
ordinary polynomial, then 
\begin{itemize}
\item $\cal P$ traps ${\cal P\/}_{\rm min\/}$.
\item $\cal P$ traps ${\cal P\/}_{\rm max\/}$.
\item For all $t \in [0,1]$ we have
 ${\cal P\/}_{\rm min\/}(t) \leq P(t) \leq  {\cal P\/}_{\rm mat\/}(t)$.
\end{itemize}

Our arithmetic operations are such that if
the polynomial
${\cal P\/}_j$ traps the polynomial $P_j$ for
$j=1,2$, then
${\cal P\/}_1 * {\cal P\/}_2$ traps
$P_1*P_2$.  Here $* \in \{+,-,\times\}$.

\section{A Bound on the Power Functions}

We are going to use Taylor's Theorem with 
Remainder to approximate the power functions.
As a preliminary step, we need some
{\it priori\/} bound on the remainder term
over a suitable interval.  We work in
the interval $[-2,16]$, which covers
all the parameters of interest to us and
a bit more.

\begin{lemma}[Remainder Bound]
$$
\bigg|\frac{d^{12}}{ds^{12}} m^{-s/2} \bigg|<1 \hskip 30 pt
\forall s \in [-2,16], \hskip 30pt
\forall m=2,3,4.
$$
\end{lemma}

\startproof
Note that $m^{-s/2} \leq m \leq 4$ when
$s \in [-2,16]$.  The max occurs when $s=-2$.
Also, one can check that $\log(4)<\sqrt 2$. Hence
$\log(4)^{12}<64.$  Hence
$$
\frac{d^{12}}{ds^{12}} m^{-s/2}=
\frac{m^{-s/2}}{4096 \log(m)}<\frac{m^{-s/2}}{64} \leq
\frac{4}{64}<1.$$
\endproof

\noindent
{\bf Remark:\/}
Note that we get the much better bound of $1/16$.
Whenever possible, we will be generous with our
bounds because we want our results to be robustly true.

\section{Approximation of Power Combos}
\label{eig}

Suppose that $Y=(a_2,a_3,a_4,b_2,b_3,b_4)$ is a $6$-tuple
of rational numbers.   The coefficients produced
by the Forcing Lemma construction have the form
\begin{equation}
C_Y(x)=
a_2  2^{-s/2} +
a_3  3^{-s/2} +
a_4  4^{-s/2} +
b_2 s 2^{-s/2} +
b_3 s 3^{-s/2} +
b_4 s 4^{-s/2}
\end{equation}
evaluated on the interval $[-2,16]$.
We call such expressions {\it power combos\/}.

For each even integer $2k=-2,...,16$.
we will construct rational polynomials
$A_{Y,2k,-}$ and $A_{Y,2k,+}$ and
$B_{Y,2k,-}$ and $B_{Y,2k,+}$ such that
\begin{equation}
\label{trap1}
A_{Y,2k,-}(t)  \leq C_Y(2k-t) \leq A_{Y,2k,+}(t), \hskip 30 pt
t \in [0,1].
\end{equation}
\begin{equation}
\label{trap2}
B_{Y,2k,-}(t)  \leq C_Y(2k+t) \leq B_{Y,2k,+}(t), \hskip 30 pt
t \in [0,1].
\end{equation}
We ignore the cases $(2,-)$ and $(16,+)$.

The basic idea is to use Taylor's Theorem with
Remainder:
\begin{equation}
\label{tayser}
m^{-s/2}=\sum_{j=0}^{11} \frac{(-1)^j \log(m)^j}{m^k 2^j j!} (s-2k)^j+
\frac{E_s}{12!} (s-2k)^{12}.
\end{equation}
Here $E_s$ is the
$12$th derivative of $m^{-s/2}$ evaluated at some point
in the interval.  Note that the only dependence on $k$ is
the term $m^k$ in the denominator, and this is a
rational number.

The difficulty with this approach is that
the coefficients of the above 
Taylor series are not rational. We get around
this trick by using interval polynomials.
We first pick specific intervals which
trap $\log(m)$ for $m=2,3,4$.  We choose
the intervals
$$
L_2=\bigg[\frac{25469}{36744},\frac{7050}{10171}\bigg], \hskip 10 pt
L_3=\bigg[\frac{5225}{4756},\frac{708784}{645163}\bigg], \hskip 10 pt
L_4=\bigg[\frac{25469}{18372},\frac{345197}{249007}\bigg].
$$
Each of these intervals has width about $10^{-10}$.
I found them using Mathematica's Rationalize function.
It is an easy exercise to check that
$\log(m) \in L_m$ for $m=2,3,4$.  
According to our Remainder Lemma, we always have
$|E_s|<1$ in the series expansion from
Equation \ref{tayser}.  Fixing $k$ we introduce
the interval Taylor series

\begin{equation}
\label{tayser2}
A_{m}(t)=\sum_{j=0}^{11} \frac{(+1)^j (L_m)^j}{m^k 2^j j!}t^j+
\left[-\frac{1}{12!},\frac{1}{12!}\right] t^{12}.
\end{equation}

\begin{equation}
\label{tayser3}
B_m(t)=\sum_{j=0}^{11} \frac{(-1)^j (L_m)^j}{m^k 2^j j!}t^j+
\left[-\frac{1}{12!},\frac{1}{12!}\right] t^{12}.
\end{equation}

By construction $A_m$ traps the Taylor series
expansion from Equation \ref{tayser} when it
is evaluated at $t=s-2k$ and $t \in [0,1]$.
Likewise 
$B_m$ traps the Taylor series expansion 
from Equation \ref{tayser} when it is
evaluated at $t=2k-s$ and $t \in [0,1]$.
Define
$$
A_Y(t)=
a_2  A_2(t)+
a_3  A_3(t)+
a_4  A_4(t)+$$
\begin{equation}
b_2 (2k-t) A_2(t)+
b_3 (2k-t) A_3(t)+
b_4 (2k-t) A_4(t)
\end{equation}

\noindent
{\bf Remark:\/} It took an embarrassing
amount of trial and error before I realized
that the coefficients $b_m(2k-t)$ are
correct.
\newline

Likewise, define

$$
B_Y(t)=
a_2  B_2(t)+
a_3  B_3(t)+
a_4  B_4(t)+$$
\begin{equation}
b_2 (2k+t) B_2(t)+
b_3 (2k+t) B_3(t)+
b_4 (2k+t) B_4(t)
\end{equation}
By construction, $A_Y(t)$ traps
$C_Y(2k-t)$ when $t \in [0,1]$ and
$B_Y(t)$ traps $C_Y(2k+t)$ when
$t \in [0,1]$.
\newline
\newline

Finally, we define
$$
A_{Y,2k,-}=(A_Y)_{\rm min\/}, \hskip 15 pt
A_{Y,2k,-}=(A_Y)_{\rm max\/},$$
\begin{equation}
B_{Y,2k,-}=(B_Y)_{\rm min\/}, \hskip 15 pt
B_{Y,2k,-}=(B_Y)_{\rm max\/}.
\end{equation}
By construction these polynomials
satisfy Equations \ref{trap1} and
\ref{trap2} respectively for each $k=-1,...,8$.
These are our under and over approximations.

The reader can print out these approximations,
either numerically or exactly, using my program.
The reader can also run the debugger and see that
these functions really do have the desired
approximation properties.

\newpage
\chapter{Positive Dominance}
\label{posdom0}

In the previous chapter we explained how to
approximate the expressions that arise
in the Forcing Lemma by polynomials.
When we make these approximations, we
reduce the Forcing Lemma to statements that
various polynomials are positive (or
non-negative) on various domains.
The next chapter explains how this goes.
In this chapter we explain how we 
prove such positivity results about
polynomials.   I will explain my Positive
Dominance Algorithm.  It is entirely
possible to me that this algorithm is
known to some community of mathematicians,
but I thought of it myself and have not
seen it in the literature.

\section{The Weak Positive Dominance Criterion}
\label{WPD}

Let 
\begin{equation}
\label{poly}
P(x)=a_0+a_1x+...+a_nx^n
\end{equation}
be a polynomial with real coefficients.  Here we describe
a method for showing that $P \geq 0$ on $[0,1]$,

Define
\begin{equation}
\label{sum1}
A_k=a_0+ \cdots + a_k.
\end{equation}
We call $P$ {\it weak positive dominant\/} (or
{\it WPD\/} for short) if
$A_k \geq 0$ for all $k$ and $A_n>0$.

\begin{lemma}
\label{PD1}
If $P$ is weak positive dominant, then
$P>0$ on $(0,1]$.
\end{lemma}

\startproof
The proof goes by induction on the degree of $P$.
The case $\deg(P)=0$ follows from the fact that $a_0=A_0>0$.
Let $x \in (0,1]$.
We have
$$P(x)=a_0+a_1x+x_2x^2+ \cdots + a_nx^n \geq $$
$$a_0x+a_1x+a_2x^2+ \cdots + a_nx^n=$$
$$x(A_1+a_2x+a_3x^2+ \cdots a_nx^{n-1})=xQ(x)> 0$$
Here $Q(x)$ is weak positive dominant and has degree $n-1$.
\endproof

We can augment this criterion by subdivision.
Given $I=[a,b] \subset \R$, let $A_I$ be one of
the two affine maps which carries $[0,1]$ to $I$. We
call the pair $(P,I)$ {\it weak positive dominant\/} if
$P \circ A_I$ is WPD. If $(P,I)$ is WPD
then $P \geq 0$ on $(a,b]$, by Lemma \ref{PD1}. 
For instance, if $P$ is WPD on $[0,1/2]$ and
$[1/2,1]$ then $P>0$ on $(0,1)$.  

\section{The Positive Dominance Algorithm}

As we suggested briefly above, the
WPD {\it criterion\/} is really a step in
a recursive subdivision
algorithm.  In explaining the algorithm, I will
work with positive dominance
rather than weak positive dominance, because
the main application 
involves positive dominance.   (For
most of our proofs, the WPD criterion is enough.)
I will treat the case when the domain is the unit cube.
I have used this method extensively in other contexts.
See e.g. [{\bf S2\/}] and [{\bf S3\/}].
In particular, in [{\bf S3\/}] I explain
the method for arbitrary polytopes.

Given a multi-index 
$I=(i_1,...,i_k) \in (\N \cup \{0\})^k$ we
let 
\begin{equation}
x^I=x_1^{i_1}...x_k^{i_k}.
\end{equation}
Any polynomial $F \in \R[x_1,...,x_k]$
can be written succinctly as
\begin{equation}
F=\sum a_I X^I, \hskip 30 pt a_I \in \R.
\end{equation}
If $I'=(i_1',...,i_k')$ we write
$I' \leq I$ if $i'_j \leq i_j$ for
all $j=1,...,k$.
We call $F$ {\it positive dominant\/}
(PD)
if
\begin{equation}
\label{summa}
A_I:=\sum_{I' \leq I} a_{I'}> 0
\hskip 30 pt \forall I,
\end{equation}

\begin{lemma}
\label{PD2}
If $P$ is PD, then $P> 0$ on $[0,1]^k$.
\end{lemma}

\startproof
When $k=1$ the proof is the
same as in Lemma \ref{PD1},
once we observe that also $P(0)>0$.
Now we prove the general case.
Suppose the the coefficients of $P$ are
$\{a_I\}$.
We write
\begin{equation}
P=f_0+f_1x_k+...+f_mx_k^m,
\hskip 20 pt f_j \in \R[x_1,...,x_{k-1}].
\end{equation} 
Let $P_j=f_0+...+f_j$.
A typical coefficient in $P_j$ has
the form
\begin{equation}
\label{expandout}
b_J=\sum_{i=1}^j a_{Ji},
\end{equation}
where $J$ is a multi-index of length $k-1$
and $Ji$ is the multi-index of length $k$
obtained by appending $i$ to $J$.
From equation \ref{expandout} and
the definition of PD, the fact that
$P$ is PD implies that $P_j$ is
PD for all $j$.
\endproof

The positive dominance criterion is not that useful in
itself, but it feeds into a powerful divide-and-conquer
algorithm. We define the maps
$$
A_{j,1}(x_1,...,x_k)=(x_1,...,x_{j-1}\frac{x_j+0}{2},x_{i+1},...,x_k),
$$
\begin{equation}
\label{polysub}
A_{j,2}(x_1,...,x_k)=(x_1,...,x_{j-1}\frac{x_j+1}{2},x_{j+1},...,x_k),
\end{equation}

We define the $j$th {\it subdivision\/} of $P$ to be the
set
\begin{equation}
\label{SUB}
\{P_{j1},P_{j2}\}=
\{P \circ A_{j,1},P \circ A_{j,2}\}.
\end{equation}

\begin{lemma}
\label{half}
$P > 0$ on $[0,1]^k$ if and only if
$P_{j1}> 0$ and $P_{j2}> 0$ on
$[0,1]^k$.
\end{lemma}

\startproof
By symmetry, it suffices to take $j=1$.
Define
\begin{equation}
[0,1]^k_1=[0,1/2] \times [0,1]^{k-1}, \hskip 30 pt
[0,1]^k_2=[1/2,1] \times [0,1]^{k-1}.
\end{equation}
Note that
\begin{equation}
A_1([0,1]^k)=[0,1]^k_1,\hskip 30 pt
B_1 \circ A_1([0,1]^k)=[0,1]^k_2.
\end{equation}
Therefore,
$P > 0$ on $[0,1]^k_1$ 
if and only if $P_{j1} > 0$ on $[0,1]^k$. Likewise
$P > 0$ on $[0,1]^k_2$ if and only if
if $P_{j2} > 0$ on $[0,1]^k$.
\endproof

Say that a {\it marker\/} is a non-negative
integer vector in $\R^k$. 
Say that the {\it youngest entry\/} in the the marker
is the first minimum entry going from left to right. 
The {\it successor\/} of a marker is the marker obtained
by adding one to the youngest entry. For instance,
the successor of $(2,2,1,1,1)$ is $(2,2,2,1,1)$.
Let $\mu_+$ denote the successor of $\mu$.

We say that a {\it marked polynomial\/} is a pair
$(P,\mu)$, where $P$ is a polynomial and
$\mu$ is a marker.  Let $j$ be the position of the
youngest entry of $\mu$.  We define the
{\it subdivision\/} of $(P,\mu)$ to be the
pair
\begin{equation}
\{(P_{j1},\mu_+),(P_{j2},\mu_-)\}.
\end{equation}
Geometrically, we are cutting the domain in half
along the longest side, and using a particular rule
to break ties when they occur.
Now we have assembled the ingredients needed
to explain the algorithm.
\newpage
\noindent
{\bf Divide-and-Conquer Algorithm:\/}
\begin{enumerate}
\item Start with a list LIST of marked polynomials.
Initially, LIST consists only of the marked polynomial
$(P,(0,...,0))$.
\item Let $(Q,\mu)$ be the last element of LIST.
We delete $(Q,\mu)$ from LIST and test whether
$Q$ is positive dominant.
\item Suppose $Q$ is positive dominant.
We go back to Step 2 if LIST is not empty.
Otherwise, we halt.
\item Suppose $Q$ is not positive dominant.
we append to LIST the two marked polynomials
in the subdivision of $(Q,\mu)$ and then go
to Step 2.
\end{enumerate}

If the algorithm halts, it constitutes a proof
that $P > 0$ on $[0,1]^k$.  Indeed, the
algorithm halts if and only if $P>0$ 
on $[0,1]^k$.
\newline
\newline
{\bf Parallel Version:\/}
Here is a variant of the algorithm.
Suppose we have a list $\{P_1,...,P_m\}$ of
polynomials and we want to show that at least
one of them is positive at each point of
$[0,1]^k$.  We do the following
\begin{enumerate}
\item Start with $m$ lists LIST$(j)$ for $j=1,...,m$
 of marked polynomials.
Initially, LIST$(j)$ consists only of the marked polynomial
$(P_j,(0,...,0))$.
\item Let $(Q_j,\mu)$ be the last element of LIST$(j)$.
We delete $(Q_j,\mu)$ from LIST$(j)$ and test whether
$Q_j$ is positive dominant.  We do this for $j=1,2,...$
until we get a success or else reach the last index.
\item Suppose {\it at least one\/}
$Q_j$ is positive dominant.
We go back to Step 2 if LIST$(j)$ is not empty.
(All lists have the same length.)
Otherwise, we halt.
\item Suppose none of $Q_1,...,Q_m$ is positive dominant.
For each $j$ we append to LIST$(j)$ the two marked polynomials
in the subdivision of $(Q_j,\mu)$ and then go
to Step 2.
\end{enumerate}
If this algorithm halts it constitutes a proof that
at least one $P_j$ is positive at each point
of $[0,1]^k$.
\newline
\newline
{\bf The Weak Version:\/}
We can run the positive dominance algorithm
using the weak positive dominance criterion
in place of the positive dominance criterion.
If the algorithm runs to completion, it
establishes that the given polynomial is
non-negative on the unit cube, and
almost everywhere positive.

\section{Discussion}

For polynomials in $1$ variable, the
method of Sturm sequences counts
the roots of a polynomial in any
given interval. An early version of
this paper used Sturm sequences (to
get more limited results) but
I prefer the positive dominance criterion.
The calculations for the positive
dominance criterion are much simpler
and easier to implement than Sturm
sequences.  

Another $1$ dimensional
technique that sometimes works is
Descartes' Rule of Signs or (relatedly) 
Budan's Theorem.  Both these techniques
involve examining the coefficients and
looking for certain patterns.  In this
way they are similar to positive dominance.
However, they do not work as well.
I could imagine placing Budan's Theorem
inside an iterative algorithm, but I
don't think it would work as well as the PDA.

There are generalizations of Sturm
sequences to higher dimensions, and
also other positivity criteria (such as
the Handelman decomposition) but
I bet they don't work as well as the
positive dominance algorithm.
Also, I don't see how to do 
the parallel positive dominance
algorithm with these other methods.

The positive dominance algorithm works
so well that one can ask why I didn't
use it to deal with the calculations
related to the Big and Small Theorems
discussed in \S \ref{bigsmall}.
For such calculations, the polynomials
involved, when expanded out, have a
vast number of terms.  The calculation
is not feasible.  I briefly tried to
set up the problem for the energy $G_3$
and already the number of terms was
astronomical.

\newpage
\chapter{Proof of The Forcing Lemma}

\section{Positivity of the Coefficients}

\noindent
{\bf Case 1:\/}
We consider $a_1$ on $(0,6]$ in detail.
\begin{itemize}
\item We set $$Y=(792, 1152, -1944, -54, -288, 0),$$ the row
of the relevant
  matrix corresponding to $a_1$.  By construction, $a_1(x)=C_Y(x)$,
the power combo defined in \S \ref{eig}.
\item We verify that the $6$
under-approximations
$$A_{Y,0,+},A_{Y,2,-},A_{Y,2,+},A_{Y,4,-},A_{4,+},A_{6,-}$$ are 
either WPD on $[0,1]$ or WPD on each of
$[0,1/2]$ and $[1/2,1]$ in all cases.
Hence these functions are
positive on $(0,1]$.  See
\S \ref{WPD}.
\item Since
$A_{Y,k,+}(t) \leq a_1(k+t)$ for $t \in [0,1]$ we
see that $a_1>0$ on $(k,k+1]$ for $k=0,2,4$.
\item Since 
$A_{Y,k,-}(t) \leq a_1(k-t)$ for $t \in [0,1]$ we
see that $a_1>0$ on $[k,k+1)$ for $k=1,3,5$.
\item We check by direct inspection that $a_1(x)>0$ at
the values $x=1,2,3,4,5,6$.
\end{itemize}

The same argument works for $a_2,a_3,a_4,\delta$ on
$(0,6]$. 
We conclude that $a_1,a_2,a_3,a_4,\delta>0$ on
$(0,6]$.  By construction $a_0(s)=2^{-s}>0$.
So, all coefficients are positive on $(0,6]$.
The statement 
$\delta>0$ means $\Gamma'(2)>R'(2)$ for all
$s \in (0,6]$.
\newline
\newline
{\bf Case 2:\/}
We do the same thing on the interval $[6,13]$,
using the intervals $[6,7],...,[12,13]$.
Again, every polynomial in sight is WPD on $[0,1]$,
and in fact PD on $[0,1]$.   Since we just check
the WPD condition, we also check that our
functions are positive at the integer values in
$[6,13]$ by hand.
We conclude 
that $a_1,a_2,a_3,a_4,\delta>0$ on $[6,13]$.
\newline
\newline
{\bf Case 3:\/}
This case is different, because we are
checking positivity on an interval whose
right endpoint is not an endpoint.
In all cases, we use the positive dominance
algorithm with the subdivision variant
that causes the intervals to be treated
in order.  The algorithm fails on
the interval $[13,16]$ but for each
coefficient it passes on $[13,14]$
and $[14,15]$ and then we check the
following:
\begin{enumerate}
\item The polynomial under-approximating $a_1$ is
also PD on $[15,481/32]$ and
$[481/32,963/64]$ and $[963/64,3853/256]$.
\item The polynomial under-approximating $a_2$ is
PD on $[15,121/8]$.
\item The polynomial under-approximating $a_3$ is
PD on $[15,121/8]$.
\item The polynomial under-approximating $a_4$ is
PD on $[15,16]$.
\item The polynomial under-approximating $\delta$ is
also PD on $[15,481/32]$ and
$[481/32,963/64]$ and $[964/64,1721/128]$.
\end{enumerate}
The numbers
$$3853/256,
\hskip 20 pt
121/8,
\hskip 20 pt
1971/128$$
all (barely) exceed $15.05$.
Hence $a_1,a_2,a_3,a_4,\delta>0$ on
$[13,15.05]$. 
\newline
\newline
{\bf Case 4:\/}
This case is just like Case 1, except
that here our function is $R_s(r)=-r^{-s}$.
The means that $a_0(s)=-2^{-s}<0$.  The same argument as in Case 1 
shows that $a_1,a_2,a_3,a_4,\delta>0$ on
$(-2,0)$.
We also check, using the same method that
$$\Gamma_{(s)}(0)=a_0+4a_1+16a_2+256a_3+1024a_5<0, \hskip 30 pt
\forall\ s \in (-2,0).$$ 

\section{Under Approximation: Case 1}

We have $s \in (0,6]$.
as we can.  In particular, we set $R=R_s$, etc.
 We know
already that $\Gamma(r)>0$
and $R(r)>0$ for all $r \in (0,6]$.
Define
\begin{equation}
\label{under}
H(r)=1-\frac{\Gamma(r)}{R(r)}=1-r^s \Gamma.
\end{equation}
We just have to show (for each value of $s$)
that $H \geq 0$ on $(0,2)$.  Let
\begin{equation}
H'=\frac{dH}{dr}.
\end{equation}

\begin{lemma}
\label{Qdefined}
$H'$ has $4$ simple roots in $(0,2)$.
\end{lemma}

\startproof
We count roots with multiplicity.
We have
\begin{equation}
\label{deriv}
H'(r)= - r^{s-1}(s\Gamma(r)+r\Gamma'(r)).
\end{equation}
Combining Equation \ref{deriv} with the general 
equation
\begin{equation}
rG_k'(r)=2kG_k(r)-8kG_{k-1}(r),
\end{equation}
we see that 
the positive roots of $H'(r)$ are the same as
the positive roots of
\begin{equation}
s\Gamma(r)+r\Gamma'(r)=(12+s)a_4 G_6(r)-48 a_4 G_5(r)+\sum_{k=1}^4 b_k G_k(r)+b_0.
\end{equation}
Here $b_0,...,b_4$ are coefficients we don't care about.
Making the substitution $t=4-r^2$ and dividing through
to make the polynomial monic, we see that the roots of
$H'$ in $(0,2)$ are in bijection with the roots in $(0,4)$ of
\begin{equation}
\label{monic2}
\psi(t)=t^6-\frac{48}{12+s} t^5 + c_4t^4 + c_3t^3+c_2t^2+c_1t+c_0.
\end{equation}
Here $c_0,...,c_4$ are coefficients we don't care about.
The change of coordinates is a diffeomorphism
from $(0,4)$ to $(0,2)$ and so it carries simple
roots to simple roots.

We just have to show that $\psi$ has $4$ simple
roots in $(0,4)$.  Note that the
sum of the $6$ roots of $\psi$ is
$48/(12+s)<4$.  This works because $s>0$ here.
The $4$ roots of $\psi$ we already know
about are $1$ and $2$ and some number
in $(0,1)$ and some number in $(1,2)$.
The sum of these roots exceeds $4$ and
so the remaining two roots cannot also
be positive.  Hence $\psi$ has at most
$5$ roots in $(0,4)$.

To finish the proof, it suffices to show
that $\psi$ has an even number of roots in
$(0,4)$.  Let's call a function
$f$ {\it essentially positive from the left\/}
at the point $x_0$ if there is an infinite
sequence of values converging to $x_0$ from
the left at which $f$ is negative.
We make the same definitions with
{\it right\/} in place of {\it left\/}
and {\it negative\/} in place of
{\it positive\/}.

To prove our claim about the parity of
the number of roots, it suffices to show
that $\psi$ is essentially negative from
the right at $0$ and essentially negative
from the left at $4$.  For this purpose,
it suffices to show that $H'$ is essentially
negative at $0$ and negative at $2$.

Note that $H(\epsilon) \to 1$ as $\epsilon \to 0$
because $R \to \infty$ and $\Gamma$ is bounded.
Also, both $\Gamma$ and $R$ are both positive.
From these properties, we see that
$H'$ is essentially negative from the right
at $0$. At the same time, we know that
$\Gamma(2)=R(2)$ and $\Gamma'(2)>R'(2)$.
But then
\begin{equation}
H'(2)=\frac{\Gamma(2)R'(2)-\Gamma'(2)R(2)}{R^2(2)}=
\frac{R'(2)-\Gamma'(2)}{R(2)}<0.
\end{equation}

This completes the proof that $\psi$ has an
even number of roots in $(0,4)$.  Since there
are at most $5$ such roots, and at least
$4$ distinct roots, there are exactly $4$
roots. Since they are all distinct, they are
all simple.
\endproof

How we mention the explicit dependence on 
$s$ and remember that we are taking about $H_s$.

\begin{lemma}
$H''_s(\sqrt 2)>0$ and
$H''_s(\sqrt 3)>0$ for all 
$s \in (0,6]$.
\end{lemma}

\startproof
We check directly that $$H_{3}''(\sqrt 2)>0,
\hskip 40 pt H''_3(\sqrt 3)>0.$$
It cannot happen that $H_s''(\sqrt 2)=0$
for other $s \in (0,6]$
because then $H'_s$ has a double
root in $(0,2)$.  Hence
$H_s''(\sqrt 2)>0$ for all $s \in (0,6]$.
The same argument shows that
$H_s''(\sqrt 3)>0$ for all $s \in (0,6]$.
\endproof

Now we set $H=H_s$ again.

\begin{lemma}
For all sufficiently small $\epsilon>0$ the
quantities
$$H(0+\epsilon), \hskip 20 pt
H(\sqrt 2 \pm \epsilon), \hskip 20 pt
H(\sqrt 3 \pm \epsilon), \hskip 20 pt
H(2-\epsilon)$$
are positive.
\end{lemma}

\startproof
We have seen already that $\lim_{\epsilon \to 0} H(\epsilon)=1$.
Likewise, we have seen that
$H'(2)<0$ and $H(2)=0$.  So $H(2-\epsilon)>0$ for
all sufficiently small $\epsilon$.
Finally, the case of $\sqrt 2$ and $\sqrt 3$
follows from the previous lemma and the second
derivative test.
\endproof

We have already proved that $H'$ has exactly
$4$ simple roots in $(0,2)$.  
In particular, the interval
$(0,\sqrt 2)$ has no roots of $H'$
and the intervals $(\sqrt 2,\sqrt 3)$
and $(\sqrt 3,2)$ have $1$ root each.
Finally, we know that $H>0$ sufficiently
near the endpoints of all these intervals.
If $H(x)<0$ for some $x \in (0,2)$, then $x$ must
be in one of the $3$ intervals just mentioned,
and this interval contains at least $2$
roots of $H'$.  This is a contradiction.
This shows that $H \geq 0$ on $(0,2)$
and thereby completes Case 1.

\section{Under Approximation: Case 2}
\label{case3}

This time we have $s \in [6,13]$ and
everything is as in Case 1.
All we have to do is show that the polynomial
$\psi$, as Equation \ref{monic2}, has
exactly $4$ simple roots in $(0,2)$.
This time we have
\begin{equation}
\label{monic3}
\psi(t)=t^{10}-\frac{80}{20+s} t^9 + b_8t^8+...+b_0.
\end{equation}
Again these coefficients depend on $s$.
\newline
\newline
{\bf Remark:\/}
$\psi$ only has $7$ nonzero terms and hence
can only have $6$ positive real roots.
The number of positive real roots is bounded
above by Descartes' rule of signs.
Unfortunately, $\psi$ turns
out to be alternating and so Descartes' rule
of signs does not eliminate the case of $6$ roots.
This approach seems useless.
The sum of the roots of $\psi$
is less than $4$, so it might seem as if
we could proceed as in Case 2.  Unfortunately, there are
$10$ such roots and this approach also seems useless.
We will take another approach to proving what we want.
\newline

\begin{lemma}
When $s=6$ the polynomial $\psi$ has $4$ simple
roots in $(0,4)$.
\end{lemma}

\startproof
We compute explicitly that
$$
\psi(t)=t^{10}-\frac{40}{13}t^9+ 
\frac{830304}{5785}t^5 -
\frac{415152}{1157} t^4 +\frac{789255}{1157} t^2
-\frac{3264104}{5785} t + \frac{115060}{1157}.
$$
This polynomial only has $4$ real roots -- the
ones we know about.  The remaining roots are
all at least $1/2$ away from the interval
$(0,4)$ and so even a very crude analysis
would show that these roots do not lie
in $(0,4)$.  We omit the details.

\begin{lemma}
Suppose, for all $s \in [6,13]$ that
$\psi$ only has simple roots in
$(0,4)$.  Then in all cases
$\psi$ has exactly $4$ such roots.
\end{lemma}

\startproof
Let $N_s$ denote the number of simple roots
of $\psi$ at the parameter $s$.  The same
argument as in Cases 1 and 2 shows that
$N_s$ is always even.  Suppose $s$ is
not constant.  Consider the infimal
$u \in (6,13]$ such that $N_u>4$.
The roots of $\psi$ vary continuously with
$s$.  How could more roots move into
$(0,4)$ as $s \to u$?

One possibility is that such a
root $r_s$ approaches from the upper half
plane or from the lower half plane.
That is, $r_s$ is not real for $s<u$.
Since $\psi$ is a real polynomial,
the conjugate $\overline r_s$ is also
a root.  The two roots
$r_s$ and $\overline r_s$ are approaching
$(0,4)$ from either side. But then the the limit
$$\lim_{s \to u} r_s$$ is a double root of
$\psi$ in $(0,4)$.  This is a contradiction.

The only other possibility is that the
roots approach along the real line.
Hence, there must be some $s<u$ such that
both $0$ and $4$ are roots of $\psi$.
But the same parity argument as in Case 1
shows $\psi(0)>0$ for all $s \in [6,13]$. 
\endproof

To finish the proof we just have to show that
$\psi$ only has simple roots in $(0,4)$ for
all $s \in [6,13]$.  We bring the
dependence on $s$ back into our notation and
write $\psi_s$.  It suffices to show that
that $\psi_s$ and $\psi'_s=d\psi_s/dr$ do not simultaneously
vanish on the rectangular domain
$(s,r) \in [6,13] \times [0,4]$.
This is a job for our method of positive
dominance.

We will explain in detail what do on the
smaller domain $$(s,r)=[6,7] \times [0,4].$$
The proof works the same for the remining
$1 \times 4$ rectangles.
The coefficients of $\psi_s$ and $\psi'_s$ are
power combos in the sense of
Equation \ref{eig}.  

We have rational vectors 
$Y_0,...,Y_9$ such that
\begin{equation}
\psi_s(r)=
\sum_{j=0}^9 C_j r^j, \hskip 20 pt
C_j=C_{Y_j}.
\end{equation}

We have the under- and over-approximations:
\begin{equation}
A_j=A_{Y_j,6,+}
\hskip 30 pt
B_j=B_{Y_j,6,+}
\end{equation}

We then define $2$-variable under- and over-approximations:
\begin{equation}
\underline \psi(t,u)= \sum_{i=0}^9 A_i(t) (4u)^i,
\hskip 30 pt
\overline \psi(t,u)= \sum_{i=0}^9 B_i(t) (4u)^i.
\end{equation}

We use $4u$ in these sums because we want our
domains to be the unit square.  We have
\begin{equation}
\label{overunder}
\underline \psi(t,u) \leq 
\psi_{6+t}(4u) \leq
\overline \psi(t,u), 
\hskip 30 pt
\forall (t,u) \in [0,1]^2.
\end{equation}

Now we do the same thing for $\psi'_s$.
We have rational vectors
$Y_0',...,Y_8' \in \Q^8$ which
work for $\psi'$ in place of $\psi$, and
this gives under- and over-approximations
$\underline \psi'$ and $\overline \psi'$ which
satisfy the same kind of equation as
Equation \ref{overunder}.

We run the parallel positive dominance algorithm
on the set of functions
$\{\underline \phi,-\overline \phi,
\underline \phi',-\overline \phi'\}$
and the algorithm halts.  This constitutes
a proof that at least one of these
functions is positive at each point.
But then at least one of $\phi_s(r)$
or $\phi'_s(r)$ is nonzero for each 
$s \in [6,13]$ and each $r \in [0,4]$.
Hence $\psi_s$ only has simple roots
in $(0,4)$.  This completes our proof
in Case 2.

\section{Under-Approximation: Case 3}

This time we work in the interval $[13,16]$.
For both Case 3 and Case 4 we are finding
some linear combination of the functions
$G_0,G_1,G_2,G_5,G_{10}$ which matches the
relevant power law at 
$\sqrt 2,\sqrt 3,\sqrt 4$.  Even though
the coefficients and the functions are
different in the two cases, the final
linear combination is the same.  Thus,
we just have to re-do Case 3 on the
intervals $[13,14]$, $[14,15]$, and
$[15,16]$.  We do this, and it works
out the same way.  This completes
our proof in Case 3.

\section{Under Approximation: Case 4}
\label{ua}

This case is just like Case 1 except for
some irritating sign changes.  Mostly we
will explain the differences between this
case and Case 1.
Here we show that
$\Gamma_{(s)}(r) \leq R_s(r)$ for all $r \in (0,2]$
and all $s \in (-2,0)$.
Here we have $R_s<0$ so we want show that
$\Gamma_{(s)}$ is even more negative. This
time we define
\begin{equation}
H=\frac{\Gamma}{R}-1.
\end{equation}
We want to show that $H \geq 0$ on $(0,2)$.

As in Case 1, we first prove that $H'$ has
exactly $4$ simple roots in $(0,2)$.
We take the same steps as in Case 1,
and this leads to the polynomial
\begin{equation}
\label{monic}
\psi(t)=t^5-\frac{40}{10+s} t^4 + b_3t^3+b_2t^2+b_1t+b_0.
\end{equation}
The polynomial $\psi$ has $5$ roots counting 
multiplicities, and also $4$ distinct roots
we know about in $(0,4)$.  We just have to
show that $\psi$ has an even number of roots.
As in Case 1, it suffices to show that
$H'$ is essentially negative from the right at $0$
and $H'(2)<0$.  The calculation that $H'(2)<0$
is the same as in Case 1.
Since $R<0$ on $(0,2)$ and
$R(0)=0$ and $\Gamma(0)<0$ we see that
$H(r) \to \infty$ as $r \to 0$.  
But this shows that $H'$ is essentially
negative from the right at $0$.
This gives us our $4$ simple roots
at every parameter $s \in (-2,0)$.

We make explicit checks that
$$H''_{(-1)}(\sqrt 2)>0, \hskip 30 pt
H''_{(-1)}(\sqrt 3)>0.$$
Given this information, the rest of the
proof is exactly like Case 1.

\newpage
\part{Divide and Conquer}
\chapter{The Configuration Space}

\section{Normalized Configurations}
\label{pn}

The next $5$ chapters are devoted to proving
the Big Theorem and the Small Theorem
from \S \ref{bigsmall}.  In this chapter
we study the moduli space in which
the calculations take place.

We consider $5$-point configurations
$\widehat p_0,...,\widehat p_4 \in S^2$ with
$\widehat p_4=(0,0,1)$.  For $k=1,2,3,4$ let
$p_k=\Sigma(\widehat p_k)$ where $\Sigma$
is the stereographic projection from
Equation \ref{stereo}. 
Before we go further, we mention that inverse
stereographic projection is given by
\begin{equation}
\label{inversestereo}
\Sigma^{-1}(x,y)=\bigg(\frac{2x}{1+x^2+y^2},\frac{2y}{1+x^2+y^2},
1-\frac{2}{1+x^2+y^2}\bigg).
\end{equation}
This is the inverse of the map $\Sigma$ in Equation \ref{stereo}.

Setting
$p_k=(p_{k1},p_{k2})$, we 
normalize so that
\begin{itemize}
\item $p_0$ lies in the positive $x$-axis. That is, $p_{01}>0$ and $p_{02}=0$.
\item $\|p_0\| \geq \max(\|p_1\|,\|p_2\|,\|p_3\|)$.
\item $p_{12} \leq p_{22} \leq p_{32}$ and $0 \leq p_{22}$.
\end{itemize}
If, in addition, we are working with a monotone decreasing
energy function, we add the condition that
\begin{equation}
p_{12} \leq 0.
\end{equation}
If this condition fails then our normalization implies
that all the points lie in the same hemisphere.  But then
one can decrease the energy by reflecting one of the points
across the hemisphere boundary.  This argument can fail for the one
energy function we consider which is not monotone,
namely $G_5^{\flat}$.

We will sometimes augment these conditions in
some circumstances, but the above is the baseline.
We call such configurations {\it normalized\/}.

\section{Compactness Results}

\begin{lemma}
\label{farout}
\label{confine}
Let $\Gamma$ be any of 
$G_3,G_4,G_5^{\flat},G_5,G_6,G_{10}^{\#},G_{10}^{\#\#}$.
Any normalized minimizer w.r.t. $\Gamma$
has $p_0 \in [0,4]$ and $\|p_k\| \leq 3/2$ for $k=1,2,3$.
\end{lemma}

\startproof
The TBP has $6$ bonds of length $\sqrt 2$,
and $3$ bonds of length $\sqrt 3$, and one
bond of length $2$.  Hence, 
\begin{equation}
G_k(\T)=3(2^{k+1}+1), \hskip 30 pt k=1,2,3,...
\end{equation}

Suppose first that $\Gamma$ is one of the
energies above, but not $G_3$ or
$G_5^{\flat}$.  These two require special consideration.
If $\|p_0\| \geq 4$ then the distance
from $\widehat p_0$ to $(0,0,1)$ is at most
$d=\sqrt{4/17}$.  We check this by computing
$\Sigma^{-1}(4,0)$, which would be the farthest
point from $(0,0,1)$ with the given constraints.
We check by direct calculation that
$\Gamma(d)>\Gamma(\T)$ in all cases.
This single bond contributes too much to 
the energy all by itself.

If the second condition fails then 
we have $\|p_k\|>3/2$ for some $k$ and
also $\|p_0\|>3/2$.
The distance from $\widehat p_k$ and $\widehat p_0$ from
$\widehat p_4$ is at most $d'=4/\sqrt 13$.  In all cases
we compute that $2\Gamma(d')>\Gamma(\T)$.  For example
$G_4(\T)=99$ whereas $2\Gamma(d') \approx 117.616$.
This deals with all cases except $G_3$ and $G_5^{\flat}$.
\newline

Now we deal with $G_3$.
We have $G_3(d)>G_3(\T)$.  This gives $\|p_0\| \leq 4$ as
in the preceding case.  However, now we see that
$2G_3(d')=42.4725$ whereas $G_3(\T)=51$.  We will
make up the rest of the energy by looking
at the bonds which do not involve $\widehat p_4$.

The distance between
any two points for the $4$-point minimizer is
at least $1$, because otherwise the energy of the
single bond, $3^3=27$, would exceed the energy 
$6 \times (4/3)^4 \in (14,15)$
of the regular tetrahedron.
Now observe that the function $G_3$ is convex
decreasing on the interval $[1,2]$.  We can
replace $G_3$ by a function $G_3'$ which
is convex decreasing on $(0,2]$ and
agrees with $G_3$ on $[1,2]$.  
For our minimizer, we have
$$\sum_{i<j<4} G_3(\|\widehat p_i-\widehat p_j\|)=
\sum_{i<j<4} G'_3(\|\widehat p_i-\widehat p_j\|) \geq^*
6 G'_3(\sqrt{8/3})>14.$$
The starred equality comes from the fact that the
regular tetrahedron is the global minimizer,
amongst $4$-point configurations, for any
convex decreasing potential.
In conclusion,
the sum of energies in the bonds not involving
$\widehat p_4$ is at least $14$.
Since $14+42>51$,
our configuration could not be an minimizer
if we have $\|p_k\|>3/2$ for some $k=1,2,3$.
\newline

Finally, we deal with $G_5^{\flat}$.
The preceding arguments used the fact that
$\Gamma \geq 0$ on $(0,2]$.  We don't have
this property for 
$G_5^{\flat}$ but it suffices to prove the result for
$\Gamma=G_5^{\flat}+30$.
A bit of calculus shows that $\Gamma>0$ on $(0,2]$.
The preceding arguments seemed to use the fact
that $\Gamma$ is decreasing, and this is not the
case in the present situation.  However, all we
really used is that $\Gamma$ is decreasing on
the intervals $(0,d]$ and $(0,d']$.  A bit
of  calculus shows that our
$\Gamma$ here is decreasing on $(0,1]$ and
we have $d,d'<1$. 

We compute 
$$\Gamma(\T)=G_5(\T)-25 G_1(\T)+300=120.$$
 We also compute that
$\Gamma(d)>692>\Gamma(\T)$ and
$2\Gamma(d')>716>\Gamma(\T).$
The rest of the proof is the same as in the
previous cases.
\endproof

Now we give a lengthy discussion which
contributes nothing to the proof but
sheds some light on what we are doing.
The precise bounds given above are
not really vital to our proof, though
they do turn out to be very convenient.
All we would really need, to get our
calculations off the ground, is some
(reasonable) constant $C$ such that
$\|p_0\|<C$ for any normalized configuration.
Given such a bound, we could reconfigure our
calculations (and estimates) to handle a
wider, but still compact,
domain of configurations.  The point I am trying to
make is that our approach is more robust than
 Lemma \ref{confine} might indicate.

Now let me turn the discussion to the function $G_2$,
which plays a special role in our proof because it
appears in all cases of the Forcing Lemma.
Fortunately, Tumanov has a very nice proof that the TBP is
the minimizer w.r.t. $G_2$, but one might wonder
what we would have done without Tumanov's result.
When we simply confine our points to the
domain given by the lemmas above, our calculations
run to completion for $G_2$, and the local analysis
of the Hessian turns out to work fine as well.
Again, if we had any reasonable bound on $\|p_0\|$
we could get the Big Theorem to work for $G_2$.

The arguments above do not work for $G_2$, because
$G_2(\T)=27$, whereas a
single bond contributes at most $16$ to the energy and
the regular tetrahedron has $G_2$ energy $32/3$.
Since $16+32/3=26+\frac{2}{3}<27$, we cannot use the
argument above to get a bound.
One way to scramble around this problem is to observe
that the argument fails only if two points are very
close together and either one makes something very
close to a regular tetrahedron with the other $3$.
An examination of the equations in \S \ref{regtet}
would eventually yield such a result.
In this case, an examination of more of the bonds would
reveal the energy to be more than $27$,
and we would get our desired bound.
Fortunately, we don't have to wander down this path.
\section{The TBP Configurations}
\label{TBP}

The TBP has two kinds of points, the two at the
poles and the three at the equator.
When $\infty$ is a polar point, the points
$p_0,p_1,p_2,p_3$ are, after suitable permutation,
\begin{equation}
(1,0), \hskip 20 pt
(-1/2,-\sqrt 3/2), \hskip 20 pt
(0,0), \hskip 20 pt
(-1/2,+\sqrt 3/2).
\end{equation}
We call this the {\it polar configuration\/}.
When $\infty$ is an equatorial point, the points
$p_0,p_1,p_2,p_3$ are, after suitable permutation,
\begin{equation}
(1,0), \hskip 20 pt (0,-1/\sqrt 3),
\hskip 20 pt
(-1,0), \hskip 20 pt
(0,1/\sqrt 3).
\end{equation}
We call this the {\it equatorial configuration\/}.
We can visualize the two configurations together 
in relation to the regular $6$-sided star. The black
points are part of the polar configuration and
the white points are part of the equatorial configuration.
The grey point belongs to both configurations.
The points represented by little squares are polar 
and the points represented by little disks are equatorial.
The beautiful pattern made by these two configurations
is not part of our proof, but it is nice to contemplate.

\begin{center}
\resizebox{!}{2.5in}{\includegraphics{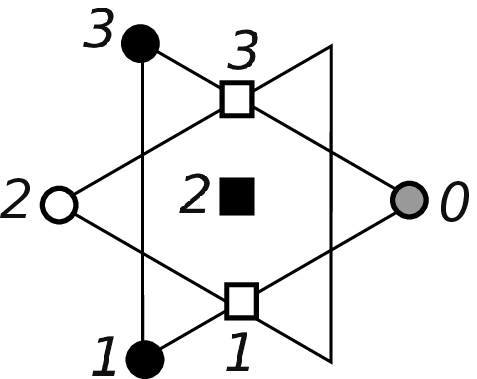}}
\newline
{\bf Figure 7.1:\/} Polar and equatorial versions of the
TBP.

\end{center}We prefer \footnote{In an earlier version of this paper,
which just proved the Big Theorem, I preferred the
polar TBP and used a trick to avoid searching near 
the equatorial TBP.  On further reflection, I found
it better to switch my preferences.  Either approach
works for the Big Theorem, but for the Small Theorem
it is much better to prefer the equatorial TBP.}
 the equatorial configuration, and we will use
a trick to avoid ever having to search very near the
polar configuration.

For each $k=0,...,4$ we introduce the quantity
\begin{equation}
\delta_k=\max_j \|\widehat p_j-\widehat p_k\|^2.
\end{equation}
We square the distance just to get a rational function.
We say that a normalized configuration is
{\it totally normalized\/} if 
\begin{equation}
\delta_4 \leq \delta_j, \hskip 30 pt j=0,1,2,3.
\end{equation}
This condition is saying that the points
in the configuration are bunched up around
$(0,0,1)$ as much as possible.  In particular,
the polar TBP is not totally normalized
and the equatorial TBP is totally normalized.
By relabeling our points we can assume that
our configurations are totally normalized.

\section{Dyadic Blocks}
\label{DYAD}

For the moment we find convenient to
only require that $p_k \in [-2,2]^2$ for $k=1,2,3$.
Later on, we will enforce the stronger condition
given by the lemmas above.  Define
\begin{equation}
\square=[0,4] \times [-2,2]^2 \times [-2,2]^2 \times [-2,2]^2.
\end{equation}
Any minimizer of any of the energies we consider
is isometric to one which is represented by a
point in this cube.  This cube is our universe.

In $1$ dimension, the {\it dyadic subdivision\/} of a
line segment is the union of the two segments obtained
by cutting it in half.  In $2$ dimensions, the
{\it dyadic subdivision\/} of a square is the
union of the $4$ quarters that result in cutting
the the square in half along both directions.
See Figure 7.2.

\begin{center}
\resizebox{!}{.7in}{\includegraphics{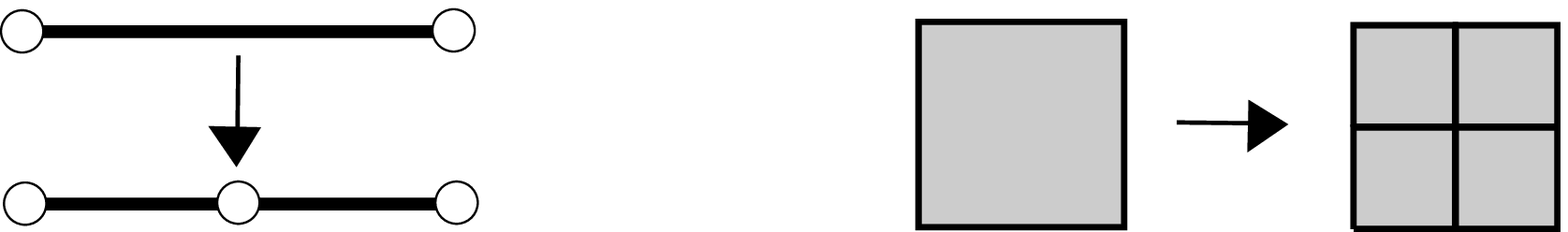}}
\newline
{\bf Figure 7.2:\/} Dyadic subdivision
\end{center}

We say that a {\it dyadic segment\/} is any segment
obtained from $[0,4]$ by applying dyadic subdivision
recursively.  We say that a {\it dyadic square\/} 
is any square obtained from $[-2,2]^2$ by applying
dyadic subdivision recursively.
We count $[0,4]$ as a dyadic segment and
$[-2,2]^2$ as a dyadic square.
\newline
\newline
{\bf Hat and Hull Notation:\/}
Since we are going to be switching back and forth between
the picture on the sphere and the picture in $\R^2$, we
want to be clear about when we are talking about
solid bodies, so to speak, and when we are talking about
finite sets of points.
We let $\langle X \rangle$ denote the convex hull
of any Euclidean subset. Thus, we think of a dyadic
square $Q$ as the set of its $4$ vertices and
we think of $\langle Q \rangle$ as the
solid square having $Q$ as its vertex set.
Combining this with our notation for stereographic
projection, we get the following notation, which
we will use repeatedly.
\begin{itemize}
\item $\widehat Q$ is a set of $4$ co-circular points
on $S^2$.
\item $\langle \widehat Q \rangle$ is a convex quadrilateral
whose vertices are $\widehat Q$.
\item $\widehat{\langle Q \rangle}$ is a ``spherical patch''
on $S^2$, bounded by $4$ circular arcs.
\end{itemize}

\noindent
{\bf Good Squares:\/}
A dyadic square is {\it good\/} if it
is contained in $[-3/2,3/2]^2$ and has
side length at most $1/2$.
Note that a good dyadic square cannot cross the 
coordinate axes.  The only dyadic square which
crosses the coordinate axes is $[-2,2]^2$, and
this square is not good.
Our computer program will only do spherical
geometry calculations on good squares.
\newline
\newline
{\bf Dyadic Blocks:\/}
We define a {\it dyadic block\/} to be a
$4$-tuple $(Q_0,Q_1,Q_2,Q_3)$, where
$Q_0$ is a dyadic segment and $Q_i$ is a
dyadic square for $j=1,2,3$.  We say that
a block is {\it good\/} if each of its $3$ component
squares is good.
By Lemma \ref{farout}, any energy minimizer for
$G_k$ is contained in a good block. Our
algorithm in \S \ref{ALG} 
quickly chops up the blocks in
$\square$ so that only good ones are considered.

The product
\begin{equation}
\langle B \rangle=\langle Q_0 \rangle \times
\langle Q_1 \rangle \times
\langle Q_2 \rangle \times
\langle Q_3 \rangle
\end{equation}
is a rectangular solid in the
configuration space $\square$.
On the other hand, the product
\begin{equation}
B= Q_0  \times
 Q_1  \times
 Q_2  \times
 Q_3 
\end{equation}
is the collection of $128$ vertices of
$\langle B \rangle$. 
We call these the
{\it vertex configurations\/} of the block.
\newline
\newline
{\bf Definition:\/}
We say that a configuration
$p_0,p_1,p_2,p_3$ is {\it in\/} the block $B$
if $p_i \in \langle Q_i \rangle$ for $i=0,1,2,3$.
In other words, the point in $\square$ representing
our configuration is contained in $\langle B \rangle$.
Sometimes we will say that this configuration is
{\it associated to\/} the block.
\newline
\newline
{\bf Sudvidision of Blocks:\/}
There are $4$ obvious subdivision operations we can
perform on a block. 
\begin{itemize}
\item The operation $S_0$ divides
$B$ into the two blocks
$(Q_{00},Q_1,Q_2,Q_3)$ and
$(Q_{01},Q_1,Q_2,Q_3)$.
Here $(Q_{00},Q_{01})$ is the dyadic subdivision of $Q_0$.
\item the operation $S_1$ divides
$B$ into the $4$ blocks
$(Q_0,Q_{1ab},Q_2,Q_3)$, where
$(Q_{100},Q_{101},Q_{110},Q_{111})$ is
the dyadic subdivision of $Q_1$.
\end{itemize}
The operations $S_2$ and $S_3$ are similar
to $S_1$.

The set of dyadic blocks has a natural poset structure
to it, and basically our algorithm does a depth-first-search
through this poset, eliminating solid blocks either
due to symmetry considerations or due to energy 
considerations. 

\section{A Technical Result about Dyadic Boxes}

The following result will be useful for
the estimates in \S \ref{sep}.  We state it
in more generality than we need, to show
what hypotheses are required, but we note
that good dyadic squares
satisfy the hypotheses.  We only
care about the result for
good dyadic squares and for
good dyadic segments. In the
case of good dyadic segments, the
lemma is obviously true.

\begin{lemma}
\label{mono}
Let $\widehat Q$ be a rectangle whose
sides are parallel to the coordinate
axes and do not cross the coordinate
axes.  Then the points of
$\widehat{\langle Q \rangle}$ closest
to $(0,0,1)$ and farthest from $(0,0,1)$ are
both vertices.
\end{lemma}

\startproof
Put the metric on $\R^2 \cup \infty$ which makes
stereographic projection an isometry.  By symmetry,
the metric balls about $\infty$ are the complements
of disks centered at $0$.  The smallest disk centered
at $0$ and containing $\langle Q \rangle$ must have a
vertex on its boundary.  Likewise, the largest disk
centered at $0$ and disjoint from the interior
of $\langle Q \rangle$ must have a vertex in its
boundary.  This latter result uses the fact that
$\langle Q \rangle$ does not cross the
coordinate axes.  These statements are equivalent
to the statement of the lemma.
\endproof

\section{Very Near the TBP}
\label{checkTBP}

Our calculations will depend on 
a pair $(S,\epsilon_0)$, both powers of two.
We always take $S=2^{30}$ and
$\epsilon=2^{-18}$.

Define the {\it in-radius\/} of a cube
to be half its side length.  Let
$P_0$ denote the configuration representing the
totally normalized polar TBP and let
$B_0$ denote the cube centered at $P_0$
ahd having in-radius $\epsilon_0$.
Note that $B_0$ is not a dyadic
block.  This does not bother us. 
Here we give a sufficient condition for
$B \subset B_0$, where $B$ is some dyadic block.

Let $a=1/\sqrt 3$.  For each choice of $S$ we
compute a value $a^*$ such that
$Sa^* \in \Z$ and $|a-a^*|<1/S$.  There are
two such choices, namely
\begin{equation}
\frac{{\rm floor\/}(Sa)}{S}, \hskip 30 pt
\frac{{\rm floor\/}(Sa+1)}{S}.
\end{equation}  
In practice, our program sets $a^*$ to be the
first of these two numbers, but we want to
state things in a symmetric
way that works for either choice.

We define
$B'_0=Q'_0 \times Q'_1 \times Q'_2 \times Q'_3$ where
$$Q_0'=[1-\epsilon_0,1+\epsilon_0],$$
$$Q_1'=[-\epsilon_0,\epsilon_0]  
\times [-S^{-1}-a^*-\epsilon_0,S^{-1}-a^*+\epsilon_0]$$
$$Q_2'=[-1-\epsilon_0,-1+\epsilon_0] \times [-\epsilon_0,\epsilon_0]$$
$$Q_3'=[-\epsilon_0,\epsilon_0]  
\times [-S^{-1}+a^*-\epsilon_0,S^{-1}+a^*+\epsilon_0]$$
By construction, we have
$B'_0 \subset B_0$.
Given a block $B=Q_0 \times Q_1 \times Q_2 \times Q_3$,
the condition
$Q_i \subset Q_i'$ for all $i$ implies
that $B \subset B_0$.   
We will see in \S \ref{computer} that this is an
exact integer calculation for us.

\newpage

\chapter{Spherical Geometry Estimates}

\section{Overview}

In this chapter we define the basic quantities
that go into the Energy Theorem, Theorem
\ref{ENERGY}.
We will persistently use the hat and hull
notation defined in \S \ref{DYAD}.  Thus:
\begin{itemize}
\item When $Q$ is a dyadic square,
$\langle \widehat Q \rangle$ is a
convex quadrilateral in space whose vertices
are $4$ co-circular points on $S^2$, and
$\widehat{\langle Q \rangle}$ is a
subset of $S^2$ bounded by $4$ circular arcs.
\item When $Q$ is a dyadic segment, 
$\langle \widehat Q \rangle$ is a
segment whose endpoints are on $S^2$, and
$\widehat{\langle Q \rangle}$ is an arc of a great circle.
\end{itemize}

Here is a summary of the quantities we will
define in this chapter.  The first three
quantities are not rational functions of the
inputs, but our estimates only use the
squares of these quantities.
\newline
\newline
{\bf Hull Diameter:\/} $d(Q)$ will be the diameter
$\langle \widehat Q \rangle$.  
\newline
\newline
{\bf Edge Length:\/} $d_1(Q)$ will be the length
of the longest edge of 
$\langle \widehat Q \rangle$. 
\newline
\newline
{\bf Circular Measure:\/} Let $D_Q \subset \R^2$
denote the disk containing 
$Q$ in its boundary.  $d_2(Q)$ is the
diameter of $\widehat D_Q$. 
\newline
\newline
{\bf Hull Separation Constant:\/}
$\delta(Q)$ will be a constant such that every
point in
$\widehat{\langle Q \rangle}$ is within
$\delta(Q)$ of a point of
$\langle \widehat Q \rangle$.
This quantity is a rational function
of the coordinates of $Q$.
\newline
\newline
{\bf Dot Product Bounds:\/}
We will introduce a (finitely computable, rational)
 quantity
$(Q \cdot Q')_{\rm max\/}$ which
has the property that
$$V \cdot V' \leq
(Q \cdot Q')_{\rm max\/}$$ for all 
$V \in \widehat{\langle Q \rangle} \cup \langle \widehat Q \rangle$
and
$V' \in \widehat{\langle Q' \rangle} \cup \langle \widehat Q' \rangle$.

\section{Some Results about Circles}

Here we prove a few geometric facts about
circles and stereographic projection.

\begin{lemma}
\label{diam}
Let $Q$ be a good dyadic square or a dyadic
segment.  The circular arcs bounding
$\widehat{\langle Q \rangle}$ lie in
circles having diameter at least $1$.
\end{lemma}

\startproof
Let $\Sigma$ denote stereographic projection.
In the dyadic segment case, $\widehat{\langle Q \rangle}$
lies in a great circle.
In the good dyadic square case, each edge of $\langle Q \rangle$
lies on a line $L$ which contains a point $p$ at most $3/2$ from the origin.
But $\Sigma^{-1}(p)$ is at least $1$ unit from $(0,0,1)$.
Hence $\Sigma^{-1}(L)$, which limits on $(0,0,1)$ and
contains $\Sigma^{-1}(p)$,
has diameter at least $1$.  The set
$\Sigma^{-1}(L \cup \infty)$ is precisely
the circle extending the relevant edge of 
$\widehat{\langle Q \rangle}$
\endproof

Let $D \subset \R^2$ be a disk of radius $r \leq R$
centered at a point which is $R$ units
from the origin. 
Let $\widehat D$ denote the corresponding disk
on $S^2$.  We consider $\widehat D$ as a subset
of $\R^3$ and compute its diameter in with
respect to the Euclidean metric on $\R^3$.

\begin{lemma}
\begin{equation}
{\rm diam\/}^2(\widehat D)=
\frac{16r^2}{1+2r^2+2R^2+(R^2-r^2)^2}.
\end{equation}
\end{lemma}

\startproof
By symmetry it suffices to consider the
case when the center of $D$ is $(R,0)$.  The
diameter is then achieved by the two
points $V=\Sigma^{-1}(R-r,0)$ and
$W=\Sigma^{-1}(R+r,0)$. 
 The formula comes
from computing $\|V-W\|^2$ and simplifying.
\endproof

We introduce the functions
\begin{equation}
\chi(D,d)=\frac{d^2}{4D}+\frac{d^4}{2D^3}
\hskip 30 pt
\chi^*(D,d)=\frac{1}{2}(D-\sqrt{D^2-d^2}).
\end{equation}
The second of these is a function closely
related to the geometry of circles.
This is the function we would use if we had
an ideal computing machine.  However, since
we want our estimates to all be rational
functions of the inputs, we will use the
first function.  We first prove an approximation
lemma and then we get to the main point.

\begin{lemma}
If $0 \leq d \leq D$ then $\chi^*(D,d) \leq \chi(D,d)$.
\end{lemma}

\startproof
If we replace $(d,D)$ by $(rd,rD)$ then both sides
scale up by $r$.  Thanks to this homogeneity, it suffices
to prove the result when $D=1$.  We have
$\chi(1,1)=3/4>1/2=\chi^*(1,1)$.  So, if suffices
to prove that the equation
$$\chi(1,d)-\chi^*(1,d)=\frac{d^2}{4}+\frac{d^4}{2} - 
\frac{1}{2}(1-\sqrt{1-d^2})$$
has no real solutions in $[0,1]$ besides $d=0$.
Consider a solution to this equation.
Rearranging the equation, we get $A=B$ where
$$A=-\frac{1}{2}\sqrt{1-d^2},
\hskip 30 pt B=\frac{d^2}{4}+\frac{d^4}{2} - 1/2.$$
An exercise in calculus shows that the
only roots of $A^2-B^2$ in $[0,1]$ are $0$ and
$\sqrt{1/2(\sqrt 8-1)}>.95$.  On the other hand
$A<0$ and $B>0$ on $[.95,1]$.
\endproof

Now we get to the key result.  This result holds in
any dimension, but we will apply it once to the
$2$-sphere, and once to circles contained in suitable
planes in $\R^3$.

\begin{lemma}
\label{circle}
Let $\Gamma$ be a round sphere of diameter $D$, contained
in some Euclidean space. Let $B$ be the
ball bounded by $\Gamma$.
Let $\Pi$ be a hyperplane which intersects
$B$ but does
not contain the center of $B$.
Let $\gamma=\Pi \cap B$ and let
$\gamma^*$ be the smaller of the
two spherical caps on $\Gamma$
bounded by $\Pi \cap \Gamma$.
Let $p^* \in \gamma^*$ be a point.
Let $p \in \gamma$ be the point so
the line $\overline{pp^*}$ contains
the center of $B$. Then
$\|p-p^*\| \leq  \chi(D,d)$.
\end{lemma}

\startproof
The given distance is maximized when
$p^*$ is the center of $\gamma^*$ and
$p$ is the center of $\gamma$. 
In this case it suffices by symmetry
to consider the situation in $\R^2$,
where $\overline{pp^*}$ is the perpendicular
bisector of $\gamma$.
Setting $x=\|p-p^*\|$, we have
\begin{equation}
\label{cross}
x(D-x)=(d/2)^2.
\end{equation}
This equation comes from a well-known theorem
from high school geometry concerning the lengths
of crossing chords inside a circle.  When we solve
Equation \ref{cross} for $x$, we see that
$x=\chi^*(D,d)$.  The previous lemma finishes the proof.
\endproof

\section{The Hull Approximation Lemma}
\label{hullsep}

\noindent
{\bf Circular Measure:\/}
When $Q$ is a dyadic square or segment, we
define
\begin{equation}
d_2(Q)={\rm diam\/}(\widehat D_Q),
\end{equation}
Where $D_Q \subset \R^2$ is such that
$Q \subset D_Q$.
So $d_2(Q)$ is the diameter of the small spherical cap
which contains $\widehat Q$ in its boundary.
Note that $\widehat{\langle Q \rangle} \subset \widehat D_Q$
by construction.  We call $d_2(Q)$ the {\it circular measure\/} of $Q$.
\newline
\newline
{\bf Hull Separation Constant:\/}
Recall that $d_1(Q)$ is the maximum
side length of $\langle Q \rangle$.  When
When $Q$ is a dyadic segment, we define
$\delta(Q)=\chi(2,d_2)$. 
When $Q$ is a good dyadic
square,
We define
\begin{equation}
\delta(Q)=\max\Big(\chi(1,d_1), \chi(2,d_2)\Big).
\end{equation} 
This definition makes sense, because
$d_1(Q) \leq 1$ and
$d_2(Q) \leq \sqrt 2<2$.  The point
here is that $\Sigma^{-1}$ is $2$-Lipschitz
and $Q$ has side length at most $1/2$.
We call $\delta(Q)$ the 
{\it Hull approximation constant\/}
of $Q$.

\begin{lemma}[Hull Approximation]
Let $Q$ be a dyadic segment or a
good dyadic square.
Every point of the spherical patch
$\widehat{\langle Q \rangle}$ is within
$\delta(Q)$ of a point of the
convex quadrilateral
$\langle \widehat Q \rangle$.
\end{lemma}

\startproof
Suppose first that $Q$ is a dyadic segment.
$\widehat{\langle Q \rangle}$ is
the short arc of a great circle
sharing endpoints with
$\langle \widehat Q \rangle$, a chord
of length $d_2$.  
By Lemma \ref{circle} each
point of on the circular arc is
within $\chi(2,d_2)$ of a
point on the chord.

Now suppose that $Q$ is a good dyadic square.
Let $O$ be the origin in $\R^3$.
Let $H \subset S^2$ denote
the set of points such that the segment
$Op^* \in H$ intesects
$\langle \widehat Q \rangle$ in a point $p$.  
Here $H$ is the itersection with the cone
over $\langle \widehat Q \rangle$ with $S^2$.
\newline
\newline
{\bf Case 1:\/}
Let $p^* \in \widehat{\langle Q \rangle} \cap H$.
Let $p \in \langle Q \rangle $ be such that
the segment $Op^*$ contains $p$. Let $B$
be the unit ball. Let
$\Pi$ be the plane containing $\widehat Q$.
Note that $\Pi \cap S^2$ bounds the spherical 
$\widehat D_Q$ which contains
$\widehat Q$ in its boundary. Let
$$
\Gamma=S^2, \hskip 20 pt
\gamma^*=\widehat D_Q, \hskip 20 pt
\gamma=\Pi \cap B.$$
The diameter of $\Gamma$ is $D=2$.  Lemma
\ref{circle} now tells us that
$\|p-p^*\| \leq \chi(2,d_2)$.
\newline
\newline
{\bf Case 2:\/}
Let $p^*\in \widehat{\langle Q\rangle}-H$.
The sets $\widehat{\langle Q\rangle}$ and
$H$ are both bounded by $4$ circular
arcs which have the same vertices. $H$
is bounded by arcs 
of great circles and $\widehat{\langle Q \rangle}$ is
bounded by arcs 
of circles having diameter at least $1$.
The point $p^*$ lies between an edge-arc $\alpha_1$ of
$H$ and an edge-arc $\alpha_2$ of $\widehat{\langle Q \rangle}$ 
which share both
endpoints.  Let $\gamma$ be the line segment
joining these endpoints.  The diameter of
$\gamma$ is at most $d_1$.

Call an arc of a circle {\it nice\/} if it is
contained in a semicircle, and if the circle
containing it has diameter at least $1$.
The arcs $\alpha_1$ and $\alpha_2$ are both nice.
We can foliate the region between $\alpha_1$ and
$\alpha_2$ by arcs of circles. These circles are
all contained in the intersection of $S^2$ with
planes which contain $\gamma$.  Call this
foliation $\cal F$.   We get this foliation by
rotating the planes around their common
axis, which is the line through $\gamma$.

Say that an ${\cal F\/}$-{\it circle\/} is
a circle containing an arc of $\cal F$.
Let $e$ be the edge of $\langle Q \rangle$
corresponding to $\gamma$.  
Call $(e,Q)$ {\it normal\/} if
$\gamma$ is never the diameter of an
$\cal F$-circle.  If $(e,Q)$  is normal, 
then the diameters of
the $\cal F$-circles interpolate
monotonically between the diameter of
$\alpha_1$ and the diameter of $\alpha_2$.
Hence, all $\cal F$-circles have diameter
at least $1$.  At the same time, if
$(e,Q)$   is normal, then all arcs of
$\cal F$ are contained in semicircles,
by continuity.  In short, if $(e,Q)$ is normal,
then all arcs of $\cal F$ are nice.
Assuming that $(e,Q)$ is normal,
let $\gamma^*$ be the arc in $\cal F$
which contains $p^*$. Let $p \in \Gamma$
be such that the line $pp^*$ contains
the center of the circle $\Gamma$ containing
$\gamma^*$.  Since
$\gamma^*$ is nice, Lemma \ref{circle} 
says that
$\|p-p^*\| \leq \chi(D,d_1) 
\leq \chi(1,d_1).$
\newline

To finish the proof,
we just have to show that $(e,Q)$ is normal.
We enlarge the set
of possible pairs we consider, by allowing
rectangles in $[-3/2,3/2]^2$ having sides parallel
to the coordinate axes and maximum side
length $1/2$.  The same arguments as above,
Lemma \ref{diam} and the $2$-Lipschitz nature
of $\Sigma^{-1}$, show that
$\alpha_1'$ and $\alpha_2'$ are still nice
for any such pair $(e',Q')$.

If $e'$ is the long side of
a $1/2 \times 10^{-100}$ rectangle $\langle Q' \rangle$ contained
in the $10^{-100}$-neighborhood of the coordinate
axes, then $(e',Q')$ is normal:
The arc $\alpha'_1$ is very nearly the
arc of a great circle and the angle between
$\alpha_1'$ and $\alpha_2'$ is very small,
so all arcs ${\cal F\/}'$ 
are all nearly arcs of great circles.
If some choice $(e,Q)$ is not normal, then
by continuity, there is a choice $(e'',Q'')$
in which $\gamma''$ is the diameter of
one of the two boundary arcs of ${\cal F\/}''$.
There is no other way to switch from normal
to not normal.
But this is absurd because the
boundary arcs, $\alpha_1''$ and $\alpha_2''$, are
nice.
\endproof

\section{Dot Product Estimates}
\label{sep}

Let $Q$ be a dyadic segment or a good dyadic square.
Let $\delta$ be the hull separation constant of $Q$.
Let $\{q_i\}$ be the points of $Q$.
We make all the same definitions for a second
dyadic square $Q'$.
We define

\begin{equation}
(Q \cdot Q')_{\rm max\/}=\max_{i,j}(\widehat q_i \cdot \widehat q_j')
+\delta+\delta'+\delta\delta'.
\end{equation}
\begin{equation}
(Q \cdot \{\infty\})_{\rm max\/}=\max_{i}\ \widehat q_i \cdot (0,0,1)
\end{equation}

\noindent
{\bf Connectors:\/}
We say that a {\it connector\/} is a line segment connecting
a point on $\widehat{\langle Q \rangle}$ to any of its closest
points in $\langle \widehat Q \rangle$.
We let $\Omega(Q)$ denote the set of connectors
defined relative to $Q$.  By the Hull Approximation Lemma,
each $V \in \Omega(Q)$ has the 
form $W+\delta U$ where $W \in \langle \widehat Q \rangle$
and $\|U\| \leq 1$.  

\begin{lemma}
\label{dotmin}
$V \cdot V' \leq
(Q \cdot Q')_{\rm max\/}$
for all $(V,V') \in \Omega(Q) \times \Omega(Q')$
\end{lemma}

\startproof
Suppose $V \in  \langle \widehat Q \rangle$ and
$V' \in  \langle \widehat Q' \rangle$.
Since the dot product is bilinear, the restriction of the dot product
to the convex polyhedral set
$\langle \widehat Q \rangle \times \langle \widehat Q' \rangle$
takes on its extrema at vertices. Hence
$V \cdot V' \leq
\max_{i,j} q_i \cdot q_i'$.
In this case, we get the desired inequality
whether or not $Q'=\{\infty\}$.

Suppose $Q' \not = \{\infty\}$ and 
$V, V'$ are arbitrary.
We use the decomposition mentioned above:
\begin{equation}
V=W+\delta U, \hskip 30 pt
V'=W'+\delta' U, \hskip 30 pt
W \in \langle \widehat Q \rangle, \hskip 30 pt
W' \in \langle \widehat Q' \rangle.
\end{equation}
But then, by the Cauchy-Schwarz inequality,
$$
|(V \cdot V')-(W \cdot W')|=
|V \cdot \delta' U' + V' \cdot \delta U + \delta U \cdot \delta' U'|
\leq \delta + \delta' + \delta \delta'.$$
The lemma now follows immediately from this
equation, the previous case applied to $W,W'$, and
the triangle inequality.

Suppose that $Q'=\{\infty\}$.
We already know the result when
$V \in \langle \widehat Q \rangle$.
When
$V \in \widehat{\langle Q \rangle}$
we get the better bound above from
Lemma \ref{mono}
and from the fact that
the dot product $V \cdot (0,0,1)$ varies
monotonically with the distance from $V$ to $(0,0,1)$
and {\it vice versa\/}.  Now we know the
result whenever $V$ is the endpoint of a
connector.  By the linearity of the
dot product, the result holds also
when $V$ is an interior point of
a connector.
\endproof

\section{Disordered Blocks}
\label{disorder}

$B=(Q_0,Q_1,Q_2,Q_3)$ be a good block.
As usual $Q_4=\{\infty\}$.  For
each index $j$ we define
\begin{equation}
(B,k)_{\rm min}=\min_{j \not = k} (Q_j \cdot Q_k)_{\rm min\/},
\hskip 30 pt
(B,k)_{\rm max}=\min_{j \not = k} (Q_j \cdot Q_k)_{\rm max\/}.
\end{equation}
We say that $B$ is {\it disordered\/} if there is
some $k \in \{0,1,2,3\}$ such that
\begin{equation}
(B,4)_{\rm max\/}<(B,k)_{\rm min\/}
\end{equation}

\begin{lemma}
If $B$ is disordered, then no configuration in $B$ is
totally normalized.
\end{lemma}

\startproof
In this situation, every configuration $\{p_j\}$ in $B$
has the following property. 
\begin{equation}
\min_{i \not = 4} p_4 \cdot p_i < \min_{i \not = k} p_k \cdot p_j.
\end{equation}
But square-distance is (linearly) decreasing as a
function of the dot product.
Hence
\begin{equation}
\delta(4)=\max_{i \not = 4} \|p_4-p_i\|^2 > 
\max_{i \not = k} \|p_k-p_i\|^2=\delta(k).
\end{equation}
In short $\delta(4)>\delta(k)$ for every configuration
associated to $B$.  This is to say that none of these
configurations is totally normalized.
\endproof

\section{Irrelevant Blocks}
\label{irrelevantX}

Call a block {\it irrelevant\/} if no
configuration in the interior of the
block is totally normalized. 
Call a block
{\it relevant\/} if it is not irrelevant.
Every relevant configuration in
the boundary of an irrelevant block is
also in the boundary of a relevant block.
So, to prove the Big and Small Theorems, we can
ignore the irrelevant blocks.

Disordered blocks are irrelevant,
but there are other irrelevant blocks.
Here we give a criterion for a good block $B$
to be irrelevant.
Given a box $Q_j$, let $\overline Q_{jk}$ and
$\underline Q_{jk}$
denote the maximum and
minimum $k$th coordinate of a point in $Q_j$.
If at least one of the following holds,
The good block is irrelevant
provided that at least one of the following
holds.
\begin{enumerate}
\item $\min(|\underline Q_{k1}|,|\overline Q_{k1}|) \geq \overline Q_{01}$
for some $k=1,2,3$.
\item $\min(|\underline Q_{k2}|,|\overline Q_{k2}|) \geq \overline Q_{01}$
for some $k=1,2,3$.
\item $\underline Q_{12} \geq 0$, provided we
 have a monotone decreasing energy.
\item $\overline Q_{22} \leq 0$.
\item $\overline Q_{32} \leq \underline Q_{22}$.
\item $\overline Q_{22} \leq \underline Q_{12}$.
\item $B$ is disordered.
\end{enumerate}

Conditions $1$ and $2$ each
imply that there is some index
$k \in \{1,2,3\}$ such that all points in the interior of
$B_k$ are farther from the origin than all
points in $B_0$. Condition $3$ implies that all points in the interior
of $B_1$ lie above the $x$-axis.  This violates our
normalization when the energy function is monotone decreasing.
Condition $4$ implies that all points in the interior
of $B_2$ lie below the $x$-axis.
Condition $5$ implies that all points in the interior
of $B_3$ lie below all points in the interior of $B_2$.
Condition $6$ implies that all points in the interior
of $B_2$ lie below all points in the interior of $B_1$.
We have already discussed Condition 7.
Thus, if any of these conditions holds, the block is
irrelevant.

\newpage

\chapter{The Energy Theorem}

\section{Main Result}
\label{ee}
\label{subdivision}

We think of the energy potential
$G=G_k$ as being a function on $(\R^2 \times \infty)^2$,
via the identification $p \leftrightarrow \widehat p$.
We take $k \geq 1$ to be an integer.

Let $\cal Q$ denote the set of dyadic squares
$[-2,2]^2$ together with the dyadic segments in $[0,4]$,
together with $\{\infty\}$.
When $Q=\{\infty\}$ the corresponding
hull separation consstant and the hull diameter are $0$.

Now we are going to define a function
$\epsilon: {\cal Q\/} \times {\cal Q\/} \to [0,\infty)$.
First of all, for notational convenience we set
$\epsilon(Q,Q)=0$ for all $Q$,
When $Q,Q' \in \cal Q$ are unequal, we define
\begin{equation}
\label{EPSILON}
\epsilon(Q,Q')=
\frac{1}{2} k(k-1) T^{k-2}d^2+
2kT^{k-1} \delta
\end{equation}
Here
\begin{itemize}
\item $d$ is the diameter of $\widehat Q$. 
\item $\delta=\delta(Q)$ is the hull approximation constant for $Q$.
See \S \ref{hullsep}.
\item
$T=T(Q,Q')=2+2(Q \cdot Q')_{\rm max\/}$.  See \S \ref{sep}.
\end{itemize}
This is a rational function in the coordinates of $Q$ and $Q'$.
The quantities $d^2$ and $\delta$ are essentially quadratic in
the side-lengths of $Q$ and $Q'$.  Note that we have
$\epsilon(\{\infty\},Q')=0$ but
$\epsilon(Q,\{\infty\})$ is nonzero when $Q \not = \{\infty\}$.

Let $B=(Q_0,Q_1,Q_2,Q_3)$. For notational
convenience we set $Q_4=\{\infty\}$.  
We define
\begin{equation}
\label{errorsum}
{\bf ERR\/}(B)=\sum_{i=0}^3
\sum_{j=0}^4 \epsilon(Q_i,Q_j).
\end{equation}

\begin{theorem}
\label{ENERGYX}
$$\min_{v \in \langle B \rangle} {\cal E\/}_k(v)>
\min_{v \in B} {\cal E\/}_k(v)-{\bf ERR\/}(B), \hskip 10 pt
\max_{v \in \langle B \rangle} {\cal E\/}_k(v) 
<\max_{v \in B} {\cal E\/}_k(v)+
{\bf ERR\/}(B).$$
\end{theorem}

\noindent
{\bf Remark:\/} 
A careful examination of our proof will reveal that,
for the case of max, one can get away with just
using the second term in the definition of $\epsilon$
in Equation \ref{EPSILON}.  However, we are so rarely
interested in this case that we will let things
stand as they are.

\section{Generalizations and Consequences}
\label{gen}

Theorem \ref{ENERGYX} suffices to deal with
$G_3,G_4,G_6$, but we need a more
general result to deal with $G_5^{\flat}$,
$G_{10}^{\#}$ and $G_{10}^{\#\#}$.
Suppose we have some energy of the form
\begin{equation}
F=\sum_{k=1}^N a_k G_k
\end{equation}
where $a_1,...,a_N$ is some sequence of
numbers, not necessarily positive.

Suppressing the dependence on $F$, we define
\begin{equation}
\epsilon(Q_i,Q_j)=\sum |a_k|\  \epsilon_k(Q_i,Q_j),
\end{equation}
where $\epsilon_k(Q_i,Q_j)$ is the above
quantity computed with respect to $G_k$.
We then define {\bf ERR\/} exactly as in
Equation \ref{errorsum}.  Here is an
immediate consequence of Theorem
\ref{ENERGYX}.  (The other statement
of Theorem \ref{ENERGYX} holds as well,
but we don't care about it.) 

\begin{theorem}[Energy]
\label{ENERGY}
$$\min_{v \in \langle B \rangle} {\cal E\/}(v) >
\min_{v \in B } {\cal E\/}(v)-
{\bf ERR\/}(B).$$
\end{theorem}

\startproof
Choose some $v \in B$ and $z \in \langle B \rangle$.
By the triangle inequality and Theorem \ref{ENERGYX},
\begin{equation}
|{\cal E\/}(v)-{\cal E\/}(z)| \leq
\sum |a_k| |{\cal E\/}_k(v)-{\cal E\/}_k(z)|
\leq \sum |a_k|\ {\bf ERR\/}_k={\bf ERR\/}.
\end{equation}
This holds for all pairs $v,z$.  The conclusion
of the Energy Theorem follows immediately.
\endproof

\begin{corollary}
\label{SUPER}
Suppose that $B$ is a block such that
\begin{equation}
\min_{v \in B} {\cal E\/}(v)-{\bf ERR\/}(B)>{\cal E\/}({\rm TBP\/}).
\end{equation}
Then all configurations in $B$ have higher
energy than the TBP.
\end{corollary}

We can write
\begin{equation}
{\bf ERR\/}(B)=\sum_{i=0}^3 {\bf ERR\/}_{i}(B), \hskip 30 pt
{\bf ERR\/}_{i}(B)=
\sum_{j=0}^4 \epsilon(Q_i,Q_j).
\end{equation}
We define the
{\it subdivision recommendation\/} 
to be the index
$i \in \{0,1,2,3\}$ for which ${\bf ERR\/}_{i}(B)$
is maximal. In the extremely unlike event that
two of these terms coincide, we pick the smaller
of the two indices to break the tie.  The subdivision
recommendation feeds into the algorithm
described in \S \ref{ALG}.

\section{A Polynomial Inequality}

The rest of the chapter is devoted to proving
Theorem \ref{ENERGYX}. 
One building block of 
Theorem \ref{ENERGYX} is
the case $M=4$ of the following inequality.

\begin{lemma}
\label{ineq}
Let $M \geq 2$ and $k=1,2,3...$
Suppose \begin{itemize}
\item $0 \leq x_1 \leq ... \leq x_M$.
\item $\sum_{i=1}^M \lambda_i=1$ and $\lambda_i \geq 0$ for all $i$.
\end{itemize}
Then
\begin{equation}
\label{INEQ}
0 \leq \sum_{i=1}^M \lambda_i x_i^k-
\bigg( \sum_{I=1}^M \lambda_i x_i\bigg)^k \leq 
\frac{1}{8} k(k-1) x_M^{k-2}\ (x_M-x_1)^2.
\end{equation}
\end{lemma} 

The lower bound is a trivial consequence of
convexity, and both bounds are trivial when
$k=1$.  So, we take $k=2,3,4,...$ and prove
the upper bound.
I discovered Lemma \ref{ineq} experimentally.  

\begin{lemma}
The case $M=2$ of Lemma \ref{ineq} implies the
rest.
\end{lemma}

\startproof
Suppose that $M \geq 3$.
We have one degree of freedom when we keep 
$\sum \lambda_i x_i$ constant and try to vary
$\{\lambda_j\}$ so as to maximize the left
hand side of the inequality.  The right hand
side does not change when we do this, and 
the left hand side varies linearly. Hence,
the left hand size is maximized when
$\lambda_i=0$ for some $i$. But then
any counterexample
to the lemma for $M \geq 3$ gives rise to
a counter example for $M-1$.
\endproof

In the case $M=2$, we set
$a=\lambda_1$.
Both sides of the inequality in
Lemma \ref{ineq} are homogeneous of
degree $k$, so it suffices to consider the
case when $x_2=1$.  We set $x=x_1$.
The inequality of is then
$f(x) \leq g(x)$, where
\begin{equation}
\label{McM}
f(x)=
(a x^k + 1-a)-(a x + 1-a)^k;
\hskip 30 pt
g(x)=\frac{1}{8}k(k-1)(1-x)^2.
\end{equation}
This is supposed to hold for
all $a, x \in [0,1]$.

The following argument is due to
C. McMullen, who figured it out
after I told him about the inequality.

\begin{lemma}
Equation \ref{McM} holds for
all $a,x \in [0,1]$ and all $k \geq 2$.
\end{lemma}

\startproof
Equation \ref{McM}. We think of $f$ as
a function of $x$, with $a$ held fixed.
Since $f(1)=g(1)=1$, it suffices
to prove that $f'(x) \geq g'(x)$ on $[0,1]$.
Define
\begin{equation}
\phi(x)=akx^{k-1}, \hskip 30 pt
b=(1-a)(1-x).
\end{equation}
We have
\begin{equation}
-f'(x)=\phi(x+b) - \phi(x).
\end{equation}
Both $x$ and $x+b$ lie in $[0,1]$. So,
by the mean value theorem there is some
$y \in [0,1]$ so that
\begin{equation}
\frac{\phi(x+b)-\phi(x)}{b}=\phi'(y)=
ak(k-1)y^{k-2}.
\end{equation}
Hence
\begin{equation}
-f'(x)=b\phi'(y)=a(1-a)k(k-1)(1-x)y^{k-2}
\end{equation}
But
$a(1-a) \in [0,1/4]$ and
$y^{k-2} \in [0,1].$
Hence
\begin{equation}
-f'(x) \leq \frac{1}{4}k(k-1)(1-x)=-g'(x).
\end{equation}
Hence $f'(x) \geq g'(x)$ for all $x \in [0,1]$.
\endproof

\noindent
{\bf Remark:\/} Lemma \ref{ineq} has the following
motivation.  The idea behind the Energy Theorem is that
we want to measure the deviation of the energy function
from being linear, and for this we would like a
quadratic estimate.  Since our energy $G_k$ involves
high powers, we want to estimate these high powers by
quadratic terms.

\section{The Local Energy Lemma}

Let $Q=\{q_1,q_2,q_3,q_4\}$ be the vertex set 
of $Q \in \cal Q$. We allow
for the degenerate case that $Q$ is a
line segment or $\{\infty\}$.  In this case we
just list the vertices multiple times,
for notational convenience.

Note that every point in the convex
quadrilateral $\langle \widehat Q \rangle$
is a convex average of the vertices.
For each $z \in \langle Q \rangle$, there is a
some point $z^* \in \langle \widehat Q\rangle$
which is as close as possible to
$\widehat z \in \widehat{\langle Q \rangle}$.  There are
constants $\lambda_i(z)$ such that
\begin{equation}
\label{zstar}
\label{starz}
z^*=\sum_{i=1}^4 \lambda_i(z)\ \widehat q_j,
\hskip 30 pt
\sum_{i=1}^4 \lambda_i(z)=1.
\end{equation}
We think of the $4$ functions
$\{\lambda_i\}$ as a partition of unity on $\langle Q \rangle$.
The choices above might not be unique,
but we make such choices once and for all for each $Q$.
We call the assignment $Q \to \{\lambda_i\}$
the {\it stereographic weighting system\/}.

\begin{lemma}[Local Energy]
Let $\epsilon$ be the function defined in
the Theorem \ref{ENERGYX}.
Let $Q,Q'$ be distinct members of $\cal Q$. Given
any $z \in Q$ and $z' \in Q'$, 
\begin{equation}
\bigg|G(z,z') - \sum_{i=1}^4 \lambda_i(z)G(q_i,z')\bigg| \leq \epsilon(Q,Q').
\end{equation}
\end{lemma}

\startproof
For notational convenience, we set $w=z'$.
Let
\begin{equation}
X=\big(2+2z^*\cdot \widehat w\big)^k.
\end{equation}
The Local Energy Lemma follows from the triangle
inequality and the following two inequalities
\begin{equation}
\label{zoop}
\bigg|\sum_{i=1}^4 \lambda_i G(q_i,w) - X\bigg| \leq
\frac{1}{2} k(k-1)T^{k-2}d^2
\end{equation}
\begin{equation}
\label{zoop2}
|X-G(z,w)| \leq
2kT^{k-1}\delta.
\end{equation}
We will establish these inequalities in turn.
\newline
\newline
Let $q_1,q_2,q_3,q_4$ be the vertices of $Q$.
Let $\lambda_i=\lambda_i(z)$.  
We set
\begin{equation}
x_i=4-\|\widehat q_i-\widehat w\|^2=2+2\widehat q_i \cdot \widehat w, \hskip 30 pt i=1,2,3,4.
\end{equation}
Note that $x_i \geq 0$ for all $i$. 
We order so that
$x_1 \leq x_2 \leq x_3 \leq x_4$. We have
\begin{equation}
\label{term1}
\sum_{i=1}^4 \lambda_i(z)G(q_i,w)=\sum_{i=1}^4 \lambda_i x_i^k,
\end{equation}
\begin{equation}
\label{term2}
X=(2+2\widehat z^* \cdot \widehat w)^k=
\bigg(\sum_{i=1}^4 \lambda_i (2+\widehat q_i \cdot \widehat w)\bigg)^k=
\bigg(\sum_{i=1} \lambda_i x_i\bigg)^k.
\end{equation}

By Equation \ref{term1}, Equation \ref{term2}, and
the case $M=4$ of Lemma \ref{ineq},
\begin{equation}
\label{boundX}
\bigg|\sum_{i=1}^4 \lambda_i G(q_i,w) - X\bigg|=
\bigg|\sum_{i=1}^4 \lambda_i x_i^k - \bigg(\sum_{i=1}^4 \lambda_i x_i\bigg)^k\bigg|
\leq
\frac{1}{8} k(k-1)x_4^{k-2}(x_4-x_1)^2.
\end{equation}

By Lemma \ref{dotmin}, we have
\begin{equation}
\label{subs1}
x_4=2+2(\widehat q_4 \cdot \widehat w) \leq
2+2(Q \cdot Q')_{\rm max\/}=T.
\end{equation}
Since $d$ is the diameter of  $\langle \widehat Q \rangle$
and $\widehat w$ is a unit vector, 
\begin{equation}
\label{subs2}
x_4-x_1=2 \widehat w \cdot (\widehat q_4-\widehat q_1)
\leq 2\|\widehat w\| \|\widehat q_4-\widehat q_1\|=
2 \|\widehat q_4-\widehat q_1\| \leq 2d.
\end{equation}
Plugging Equations \ref{subs1} and
\ref{subs2} into Equation
\ref{boundX}, we get Equation \ref{zoop}.

Now we establish Equation \ref{zoop2}.
Let $\gamma$ denote the unit speed line segment connecting
$\widehat z$ to $z^*$.  Note that the length $L$ of
$\gamma$ is at most $\delta$, by the Hull Approximation Lemma.
Define
\begin{equation}
f(t)=
\bigg(2+2 \widehat w \cdot \gamma(t)\bigg)^k.
\end{equation}
We have $f(0)=X$. Since $\widehat w$ and $\gamma(1)=\widehat z$
are unit vectors, $f(L)=G(z,w)$. Hence
\begin{equation}
\label{integ}
X-G(z,w)=f(0)-f(L), \hskip 30 pt L \leq \delta.
\end{equation}
By the Chain Rule,
\begin{equation}
f'(t)=\big(2 \widehat w \cdot \gamma'(t)\big) \times
k\bigg(2+2\widehat w \cdot \gamma(t)\bigg)^{k-1}.
\end{equation}
Note that $|2\widehat w \cdot \gamma'(t)| \leq 2$
because both of these vectors are unit vectors.
$\gamma$ parametrizes one of the connectors from
Lemma \ref{dotmin}, so 
\begin{equation}
\label{diff}
|f'(t)| \leq 2k\bigg(2+2\widehat w \cdot \gamma(t)\bigg)^{k-1}  \leq 
2k\bigg(2+2(Q \cdot Q')_{\rm max\/}\bigg)^{k-1}=2kT^{k-1}.
\end{equation}
Equation \ref{zoop2} now follows from Equation \ref{integ},
Equation \ref{diff}, and integration.
\endproof

\section{From Local to Global}

Let $\epsilon$ be the function from
the Energy Theorem.  
Let $B=(Q_0,...,Q_N)$ be a list of
$N+1$ elements of $\cal Q$.  We care about
the case $N=4$ and $Q_4=\{\infty\}$, but
the added generality makes things clearer.
Let 
$q_{i,1},q_{i,2},q_{i,3},q_{i,4}$ be the vertices of $Q_i$.
The vertices of $\langle B \rangle$ are indexed by
a multi-index $$I=(i_0,...,i_n) \in \{1,2,3,4\}^{N+1}.$$
Given such a multi-index, which amounts to a choice
of vertex of $\langle B \rangle$, 
we define the energy of the corresponding vertex configuration:
\begin{equation}
{\cal E\/}(I)=
{\cal E\/}(q_{0,i_0},...,q_{N,i_N})
\end{equation}
We will prove the following sharper result.

\begin{theorem}
\label{ENERGY2}
Let $z_0,...,z_N \in \langle B \rangle$.
Then
\begin{equation}
\label{AVE}
\bigg|{\cal E\/}(z_0,...,z_N) -
\sum_{I} \lambda_{i_0}(z_0)...\lambda_{i_N}(z_N){\cal E\/}(I) \bigg|
\leq \sum_{i=0}^N \sum_{j=0}^N \epsilon(Q_i,Q_j).
\end{equation}
The sum is taken over all multi-indices.
\end{theorem}

\begin{lemma}
Theorem \ref{ENERGY2} implies 
Theorem \ref{ENERGYX}.
\end{lemma}

\startproof
Notice that
\begin{equation}
\sum_{I} \lambda_{i_0}(z_0)...\lambda_{i_N}(z_N)=
\prod_{j=0}^N \bigg(\sum_{a=1}^4 \lambda_a(z_j)\bigg)=1.
\end{equation}
Therefore
\begin{equation} 
\label{minmax}
\min_{v \in B} {\cal E\/}(v) \leq
\sum_{I} \lambda_{i_0}(z_0)...\lambda_{i_N}(z_N){\cal E\/}(I)
\leq \max_{v \in B} {\cal E\/}(v),
\end{equation}
because the sum in the middle is the convex
average of vertex energies.

We will deal with the min case of 
Theorem \ref{ENERGYX}.  The max case has
the same treatment.
Choose some $(z_1,...,z_N) \in B$
which minimizes ${\cal E\/}$.  We have
$$
0 \leq \min_{v \in B} {\cal E\/}(v)-\min_{v \in \langle B \rangle}
{\cal E\/}(v)=\min_{v \in B} {\cal E\/}(v)-{\cal E\/}(z_0,...,z_N) \leq $$
\begin{equation}
\sum_{I} \lambda_{i_0}(z_0)...\lambda_{i_N}(z_N){\cal E\/}(I)-
{\cal E\/}(z_0,...,z_N) \leq  \sum_{i=0}^N \sum_{j=0}^N \epsilon(Q_i,Q_j).
\end{equation}
The last expression is {\bf ERR\/} when $N=4$ and $Q_4={\infty\/}$.
\endproof

We now prove Theorem \ref{ENERGY2}.

\subsubsection{A Warmup Case}

Consider the case when $N=1$. 
Setting $\epsilon_{ij}=\epsilon(Q_i,Q_j)$, the
Local Energy Lemma gives us
\begin{equation}
G(z_0,z_1) \geq 
\sum_{\alpha=1}^4 \lambda_{\alpha}(z_0)G(q_{0\alpha},z_1)-\epsilon_{01}.
\end{equation}
\begin{equation}
G(q_{0\alpha},z_1) \geq
\sum_{\beta=1}^4 \lambda_{\beta}(z_1)G(q_{1\beta}(z_1),q_{0\alpha}) - \epsilon_{10}.
\end{equation}
Plugging the second equation into the first and using 
$\sum \lambda_{\alpha}(z_0)=1$, we have
\begin{equation}
\label{AVE1}
G(z_0,z_1) \geq \bigg(\sum_{\alpha,\beta} \lambda_{\alpha}(z_0)\lambda_{\beta}(z_1)
G(q_{0\alpha},q_{1\beta})\bigg) - (\epsilon_{01}+\epsilon_{10}).
\end{equation}
Similarly,
\begin{equation}
\label{AVE1B}
G(z_0,z_1) \leq \bigg(\sum_{\alpha,\beta} \lambda_{\alpha}(z_0)\lambda_{\beta}(z_1)
G(q_{0\alpha},q_{1\beta})\bigg) + (\epsilon_{01}+\epsilon_{10}).
\end{equation}
Equations \ref{AVE1} and \ref{AVE1B} are equivalent
to Equation \ref{AVE} when $N=1$.

\subsubsection{The General Case}

Now assume that $N \geq 2$.
We rewrite Equation \ref{AVE1} as follows:
\begin{equation}
\label{AVE2}
G(z_0,z_1) \geq \sum_{A} \lambda_{A_0}(z_0)\lambda_{A_1}(z_1)\ G(q_{0A_0},q_{1A_1})-
(\epsilon_{01}+\epsilon_{10}).
\end{equation}
The sum is taken over multi-indices $A$ of length $2$.

We also observe that
\begin{equation}
\label{AVE3}
\sum_{I'} \lambda_{i_2}(z_2)...\lambda_{i_N}(z_N)=1.
\end{equation}
The sum is taken over all multi-indices
$I'=(i_2,...,i_N)$.
Therefore, if $A$ is held fixed, we have
\begin{equation}
\lambda_{A_0}(z_0)\lambda_{A_1}(z_1) =
\sum'_{I} \lambda_{I_0}(z_0)...\lambda_{I_N}(z_N).
\end{equation}
The sum is taken over all multi-indices of length $N+1$ which 
have $I_0=A_0$ and $I_1=A_1$.
Combining these equations, we have
\begin{equation}
G(z_0,z_1) \geq 
\sum_I \lambda_{I_0}(z_0)... \lambda_{I_N}(z_N) G(q_{0I_0},q_{1I_1})-
(\epsilon_{01}+\epsilon_{10}).
\end{equation}
The same argument works for other pairs of indices, giving
\begin{equation}
\label{AVE4}
G(z_i,z_j) \geq 
\sum_I \lambda_{I_0}(z_0)... \lambda_{I_N}(z_N) G(q_{iI_i},q_{jI_j})-
(\epsilon_{ij}+\epsilon_{ji}).
\end{equation}

Now we interchange the order of summation and observe that
$$
\sum_{i<j}
\bigg( \sum_I \lambda_{I_0}(z_0)...\lambda_{I_N}(z_N)\ G(q_{iI_i},q_{jI_j})\bigg)=$$
$$
\sum_I \sum_{i<j} \lambda_{I_0}(z_0)... \lambda_{I_N}(z_N)\ G(q_{iI_i},q_{jI_j})=
$$
$$
\sum_I \lambda_{I_0}(z_0)... \lambda_{I_N}(z_N) 
\Bigg(\sum_{i<j} G(q_{iI_i},q_{jI_j})\Bigg)=
$$
\begin{equation}
\sum_I \lambda_{I_0}(z_0)... \lambda_{I_N}(z_N)\ {\cal E\/}(I).
\end{equation}

When we sum Equation \ref{AVE4} over all $i<j$, we
get 
\begin{equation}
{\cal E\/}(z_0,...,z_N) \geq
\sum_{I} \lambda_{i_0}(z_0)...\lambda_{i_N}(z_N){\cal E\/}(I)-
 \sum_{i=0}^N \sum_{j=0}^N \epsilon(Q_i,Q_j).
\end{equation}

Similary,

\begin{equation}
{\cal E\/}(z_0,...,z_N) \leq
\sum_{I} \lambda_{i_0}(z_0)...\lambda_{i_N}(z_N){\cal E\/}(I) +
 \sum_{i=0}^N \sum_{j=0}^N \epsilon(Q_i,Q_j).
\end{equation}

These two equations together are equivalent to
Theorem \ref{ENERGY2}.
This completes the proof.

\newpage

\chapter{The Main Algorithm}
\label{ALG}

\section{Grading a Block}
\label{GRADE}

In this section we describe what we mean
by {\it grading\/} a block.  
The algorithm depends
on the constants $(S,\epsilon_0)$ from
\S \ref{checkTBP}.  We take
$S=2^{30}$ and $\epsilon_0=2^{-18}$.

Let $B_0$ denote the cube of in-radius
$\epsilon_0$ about
configuration of $\square$ representing
the normalized TBP.
We fix some energy $G_k$.  We perform the
following tests on a block 
$B=(Q_0,Q_1,Q_2,Q_3)$, in the order listed.
\begin{enumerate}

\item If some component square $Q_i$ of $B$
has side length more than $1/2$ we fail $B$ and
recommend that $B$ be subdivided along the
first such index.  This step guarantees that
we only pass good blocks.

\item If $B$ satisfies the criteria in
\S\ref{irrelevantX}, and hence is irrelevant, we
pass $B$.

\item If we compute that $Q_i \not \subset
[-3/2,3/2]^2$ for some $i=1,2,3$, we pass $B$.
Given the previous step, $Q_i$ is disjoint
from $(-3/2,3/2)^2$.

\item If the calculations in \S \ref{checkTBP} show
that $B \subset B_0$, we pass $B$. 

\item If the calculations in
\S \ref{ee} show that
$B$ satisfies Corollary \ref{SUPER}, we
pass $B$.
Otherwise, we fail $B$ and pass along the
recommended subdivision.
\end{enumerate}

\section{The Divide and Conquer Algorithm}
\label{divideandconquer}

Now we describe the divide-and-conquer algorithm.

\begin{enumerate}
\item Begin with a list LIST of blocks in $\square$.
Initially LIST consists of a single element,
namely $\square$.

\item Let $B$ be the last member of LIST. We delete
$B$ from LIST and then we grade $B$.

\item Suppose $B$ passes.  If
LIST is empty, we halt and declare success.
Otherwise, we return to Step 2.

\item Suppose $B$ fails.  
In this case, we subdivide $B$ along the
subdivision recommendation and we
append to LIST the
subdivision of $B$.  Then we return
to Step 2.
\end{enumerate}

If the algorithm halts with success, it implies
that every relevant block $B$ either lies
in $B_0$ or does not contain a minimizer.

\section{Results of the Calculations}

For $G_3,G_4,G_5^{\flat},G_6,G_{10}^{\sharp\sharp}$, 
I ran the programs (most recently) in late September 2016,
on my 2014 IMAC.

\begin{itemize}

\item For $G_3$ the program finished in about
$1$ hour.

\item For $G_4$ the program finished in about
$1$ hour and $9$ minutes.

\item For $G_5^{\flat}$ the program finished in about
$5$ hours and $40$ minutes.

\item For $G_6$ the program finished in about
$3$ hours and $6$ minutes.

\item For $G_{10}^{\sharp\sharp}$ the program finished in about
$25$ hours and $32$ minutes.
\end{itemize}

These programs are about $5$ times as fast without
the interval arithmetic.

In each case, the program
produces a partition of $\square$
into $N_k$ smaller blocks, each of
which is either irrelevant, contains no
minimizer, or lies in $B_0$. Here
$$
(N_3,N_4,N_5^{\flat},N_6,N_{10}^{\sharp\sharp})=
(4800136,5302730,13247122,13212929,41556654).
$$

These calculations rigorously establish the
following result.

\begin{lemma}
\label{energy1}
Let $B_0 \subset \square$ denote the cube of 
in-radius $2^{-18}$ about the equatorial
version of the TBP.  
If $P \in \square$ is a minimizer with
respect to any of
$G_3,G_4,G_5^{\flat},G_6,G_{10}^{\sharp\sharp}$ then $P$ is
represented by a configuration in $B_0$.
\end{lemma}

Now we turn to the calculation with
$G_{10}^{\sharp}$.  For reference,
here are the inequalities defining {\bf SMALL\/}.
\begin{enumerate}
\item $\|p_0\| \geq \|p_k\|$ for $k=1,2,3$.
\item $512 p_0 \in [433,498] \times [0,0]$.
\item $512 p_1 \in [-16,16] \times [-464,-349]$.
\item $512 p_2 \in [-498,-400] \times [0,24]$.
\item $512 p_3 \in [-16,16] \times [349,364]$.
\end{enumerate}
Conditions 2-5 are obviously linear inequalities
involving dyadic rationals.
To test when a block lies in {\bf SMALL\/} we
are just testing linear inequalities for
$128$ vertices.
Given the dyadic nature of the coefficients,
this turns out to be an exact integer calculation for us.

We run the algorithm above for $G_{10}^{\sharp}$ except that
we replace $B_0$ in Step 4 of the grading with the
set $B_0 \cup {\bf SMALL\/}$.   That is, we first
check if $B \subset B_0$, as before, and then we
check if $B \subset {\bf SMALL\/}$. 
The program finished in about
$35$ hours and $42$ minutes, and 
produced a paritition of size $56256273$.  This
rigorously establishes the following result.

\begin{lemma}
\label{energy2}
Let $B_0 \subset \square$ denote the cube of 
in-radius $2^{-18}$ about the equatorial
version of the TBP.  
If $P \in \square$ is a minimizer with
respect to $G_{10}^{\sharp}$ then either $P$ is
represented by a configuration in {\bf SMALL\/}
or by a configuration in $B_0$.
\end{lemma}

The rest of the chapter discusses the implementation
of the algorithm.

\section{Subdivision Recommendation}

We first discuss a subtle point about
the way the subdivision algorithm runs.
The Energy Theorem is designed to run
with the subdivision algorithm.  As
discussed at the end of \S \ref{gen} and
mentioned again above in the description
of the algorithm, the Energy Theorem comes
with a recommendation for how
to subdivide the block in the way that is
most likely to reduce the error term.
Call this the {\it sophisticated approach\/}.
The sophisticated approach might not always
subdivide the longest side of the block.

As an alternative, we 
could perform the subdivisions
just we se did for the positive dominance
algorithm, always dividing so as to cut
the longest side in half.  Call this the
{\it naive method\/}.  

The sophisticated method is more
efficient.  As an experiment, I ran
the program for $G_3$, using floating
point calculations, using the two
approaches.  The sophisticated approach
completed in $16$ minutes and produced
a partition of size $4799998$.
The naive approach completed in $32$
minutes and produced a partition of size
$11372440$.

\section{Interval Arithmetic}

As we mentioned in the introduction, everything we
need to compute is a rational function of the
coordinates of the block vertices. Thus, the
calculations only involve the operations
plus, minus, times, divide.  In theory, they
could be done with integer arithmetic. 
This seems to slow, so we do the calculations
using interval arithmetic.  Our implementation is
like the in [{\bf S1\/}] but here
we don't need to worry about the square-root function.

Java represents real numbers by {\bf doubles\/}, essentially
according to the scheme
discussed in [{\bf I\/}, \S 3.2.2].  A double is a
$64$ bit string where $11$ of the bits control the
exponent, $52$ of the bits control the binary expansion,
and one bit controls the sign.
The non-negative doubles have a lexicographic ordering, and this ordering
coincides with the usual ordering of the real numbers they
represent.   The lexicographic ordering for the non-positive doubles
is the reverse of the usual ordering of the real numbers they
represent.  To {\it increment\/} $x_+$ of
 a positive double $x$ is the very next double
in the ordering.  This amounts to treating the last $63$ bits of the
string as an integer (written in binary) and adding $1$ to it.
With this interpretation, we have $x_+=x+1$.
We also have the decrement $x_-=x-1$.
Similar operations are defined on the non-positive doubles.
These operations are not defined on the largest and smallest
doubles, but our program never encounters (or comes anywhere near)
these.

Let $\D$ be the set of all doubles.
Let 
\begin{equation}
\R_0=\{x \in \R|\ |x| \leq 2^{50}\}
\end{equation}
Our choice of $2^{50}$ is an arbitrary but convenient cutoff.
Let $\D_0$ denote the set of doubles representing reals in $\R_0$.

According to [{\bf I\/}, 3.2.2, 4.1, 5.6], there is a map
$\R_0 \to \D_0$ which maps each $x \in \R_0$ to some
$[x] \in \D_0$ which is closest to $x$.   In case there
are several equally close choices, the computer chooses one
according to the method
in [{\bf I\/}, \S 4.1].    This ``nearest point projection''
exists on a subset of $\R$ that is much larger
than $\R_0$, but we only need to consider
$\R_0$.  We also have the inclusion $r: \D_0 \to \R_0$, which
maps a double to the real that it represents.   

The basic operations plus, minus, times, divide act on $\R_0$
in the usual way, and
operations with the same name act on $\D_0$.
Regarding these operations, [{\bf I\/}, \S 5] states that
{\it each of the operations shall be performed as if it first produced an
intermediate result correct to infinite precision and with unbounded range, and then 
coerced this intermediate result to fit into the destination's format\/}.  Thus,
for doubles $x,y \in \D_0$ such that $x*y \in \D_0$ as well, we have
\begin{equation}
\label{rule}
x*y=[r(x)*r(y)]; \hskip 20 pt
* \in \{+,-,\times,\div\}.
\end{equation}
The operations on the left hand side represent operations on doubles
and the operations on the right hand side represent operations on
reals.

It might happen that $x,y \in \D_0$ but $x*y$ is not.
To avoid this problem, we make the following
checks before performing any arithmetic operation.
\begin{itemize}
\item For addition and subtraction,
$\max(|x|,|y|)\leq2^{40}$.
\item For multiplication, $\max(|x|,|y|)<2^{40}$ and
$\min(|x|,|y|)<2^{10}$.
\item For division,
$|x|\leq2^{40}$ and $|y|\leq2^{10}$ and
$|y| \geq 2^{-10}$.
\end{itemize}
We set the calculation to abort if any of these
conditions fails.
\newline

For us, an
{\it interval\/} is a pair $I=(x,y)$ of doubles with $x \leq y$ and
$x,y \in \D_0$.   Say
that $I$ {\it bounds\/} $z \in \R_0$ if $x \leq [z] \leq y$.  This is true
if and only if $x \leq z \leq y$.  Define
\begin{equation}
[x,y]_o=[x_-,y_+].
\end{equation}
This operation is well defined for doubles in $\D_0$.
We are essentially {\it rounding out\/} the endpoints of the interval.
Let $I_0$ and $I_1$ denote the left and right endpoints of $I$.
Letting $I$ and $J$ be intervals, we define
\begin{equation}
\label{operate}
I*J = (\min_{ij} I_i *I_j,\max_{ij} I_i*I_j)_o.
\end{equation}
That is, we perform the operations on all the endpoints, order
the results, and then round outward. 
Given Equation \ref{rule}, we the interval
$I*J$ bounded $x*y$ provided that $I$ bounds $x$ and $J$ bounds $y$.
This operation is very similar to what we did
in \S \ref{rationalinterval}, except that here
we are rounding out to take care of the round-off error.
This rounding out process potentially destroys the
semi-ring structure discussed in \S \ref{rationalinterval}
but, as we remarked there, this doesn't bother us.

We also define an interval version of a vector in $\R^3$.
Such a vector consists of $3$ intervals.  The only
operations we perform on such objects are addition, subtraction,
scaling, and taking the dot product.  These operations are
all built out of the arithmetic operations.

All of our calculations come down to proving inequalities 
of the form $x<y$.  We imagine that $x$ and $y$ are the
outputs of some finite sequence of arithmetic operations
and along the way we have intervals
$I_x$ and $I_y$ which respectively bound $x$ and $y$.
If we know that the right endpoint of $I_x$ is less
than the left endpoint of $I_y$, then this constitutes
a proof that $x<y$.  The point is that the whole
interval $I_x$ lies to the left of $I_y$ on the number line.

\section{Integer Calculations}
\label{computer}

Now let me discuss the implementation of the
divide and conquer algorithm.
We manipulate blocks and dyadic squares
using {\bf longs\/}.  These are $64$ bit
integers.  Given a dyadic square $Q$
with center $(x,y)$ and side length $2^{-k}$,
we store the triple 
\begin{equation}
(Sx,Sy,k).
\end{equation}
Here $S=2^{25}$ when we do the calculations
for $G_3,G_4,G_5,G_6$ and
$S=2^{30}$ when we do the calculation for
$\widehat G_{10}$.  The reader can modify
the program so that it uses any power of 
$2$ up to $2^{40}$.  Similarly,
we store a dyadic segment with
center $x$ and side length
$2^{-k}$ as $(Sx,k)$.

The subdivision is then obtained by manipulating
these triples.  For instance, the top right square
in the subdivision of $(Sx,Sy,k)$ is
$$\big(Sx-2^{-k+1}S,Sy-2^{-k+1}S,k+1\big).$$
The scale $2^N$ allows for $N$ such subdivisions
before we lose the property that the squares
are represented by integer triples.
The biggest dyadic square is stored as
$(0,0,-2)$, and each subdivision increases
the value of $k$ by $1$.  We terminate the
algorithm if we ever arrive at a dyadic
square whose center is not an even pair of integers.
We never reach this eventuality when we
run the program on the functions from the
Main Theorem, but it does occur if we try
functions like $G_7$.

Steps 1,3, and 4 of the block grading just involve
integer calculations. For instance, the point
of scaling our square centers by $S$ is that
the inequalities which go into the calculations
in \S \ref{checkTBP} are all integer inequalities.
We are simply clearing denominators.

Most of Step 2 just requires exact integer calculations.
The only part of Step 2 that requires interval
arithmetic calculations is the check that a block
could be disordered, as in \S \ref{disorder}.
In this case, we convert each long $L$ representing
a coordinate of our block into the dyadic
interval
\begin{equation}
[L/S,L/S], \hskip 30 pt L=2^{30}.
\end{equation}
We then perform the calculations using interval
arithmetic.  The first step in the calculation
is to apply inverse stereographic projection,
as in Equation \ref{inversestereo}.  The
remaining steps are as in \S \ref{disorder}.

Step 5 is the main step that requires interval
arithmetic.  We do exactly the same procedure,
and then perform the calculations which go into
the Energy Theorem using interval arithmetic.

\section{A Gilded Approach}

In the interest of speed, we take what might be called
a gilded approach to interval arithmetic.  We run the
calculations using integer arithmetic and
floating point arithmetic until the
point where we have determined that a block passes for
a reason that requires a floating point calculation.
Then we re-do the calculation using the interval
arithmetic.  This approach runs a bit faster, and
only integer arithmetic and interval arithmetic
are used to actually pass blocks.

\section{A Potential Speed Up}

One thing that probably slows our algorithm
down needlessly is the computation of the
constants $d^2$ and $\delta$ from Equation
\ref{EPSILON}. These functions depend
on the associated dyadic square or
segment in a fairly predictable way, and
it ought to be possible to replace the
calculations of $d^2$ and $\delta$ with
an {\it a priori\/} estimate.  Such estimates
might also make it more feasible to run
exact integer versions of the programs,
because the sizes of the rational numbers
in $d^2$ and $\delta$ get very large.  I haven't
tried to do this yet, but perhaps a later
version of the program will do this.

\section{Debugging}

One serious concern about any computer-assisted
proof is that some of the main steps of
the proof reside in computer programs which
are not printed, so to speak, along with
the paper.  It is difficult for one to
directly inspect the code without a serious
time investment, and indeed the interested
reader would do much better simply to reproduce
the code and see that it yields the same results.

The worst thing that could happen is
if the code had a serious bug which caused it
to suggest results which are not actually true.
Let me explain the extent to which I have debugged
the code.  Each of the java programs has a
debugging mode, in which the user can test that
various aspects of the program are running correctly. 
While the debugger does not check every method, it
does check that the main ones behave exactly
as expected.

Here are some debugging features of the program.

\begin{itemize}
\item You can check on random inputs that
the interval arithmetic operations are
working properly.
\item You can check on random inputs that
the vector operations - dot product, addition, etc. - are
working properly.
\item You can check for random dyadic squares that
the floating point and interval arithmetic
measurements match in the appropriate sense.
\item You can select a block of your choice
and compare the estimate from the Energy
Theorem with the minimum energy taken over
a million random configurations in the block.
\item You can open up an auxiliary window
and see the grading step of the algorithm
performed and displayed for a block of
your choosing.
\end{itemize}

\newpage

\chapter{Local Analysis of the Hessian}
\label{local}

The purpose of this chapter is to prove, in all relevant
cases that the Hessian of the energy function
$\cal E$ is positive definite throughout the
neighborhood $B_0$ mentioned in Lemmas
\ref{energy1} and \ref{energy2}.  This
result finishes the proof of the Big Theorem
and the Small Theorem.

\section{Outline}
\label{TBP2}

We begin with a well-known lemma about the Hessian
derivative of a function.  The Hessian is the
matrix of second partial derivatives.

\begin{lemma}
Let $f: \R^n \to \R$ be a smooth function.
Suppose that the Hessian of $f$ is positive
definite in a convex neighborhood $U \subset \R^n$.
Then $f$ has at most one local minimum in $U$.
\end{lemma}

\startproof
The fact that the Hessian is positive definite
in $U$ means that the restriction of $f$ to
each line segment in $U$ is a convex function
in one variable.  
If $F$ has two local minima in $U$ then
let $L$ be the line connecting these minima.
The restriction of $f$ to $L \cap U$ has
two local minima and hence cannot be a convex
function.  This is a contradiction.
\endproof

Now we turn our attention to the Big and Small Theorems.
Let $\Gamma$ be any of the energy functions we have
been considering.
We have the energy map
${\cal E\/}_{\Gamma}: B_0 \to \R_+$
given by
\begin{equation}
{\cal E\/}_{\Gamma}(x_1,...,x_7)=
\sum_{i<j} \Gamma(\Sigma^{-1}(p_i)-\Sigma^{-1}(p_j)).
\end{equation}
Here we have set $p_4=\infty$, and
$p_0=(x_1,0)$ and $p_i=(x_{2i},x_{2i+1})$ for $i=1,2,3$.
As usual $\Sigma^{-1}$ is 
inverse stereographic projection.
See Equation \ref{inversestereo}.

It follows
from symmetry, and also from a direct
calculation, that $\T$ is a critical
point for ${\cal E\/}_{\Gamma}$.  Assuming
that $\Gamma$ is one of
$G_3,G_4,G_5^{\flat},G_6,G_{10}^{\#\#}$, the
Big Theorem now follows from
Lemma \ref{energy1} and
from the statement that the Hessian of
${\cal E\/}_{\Gamma}$ is positive definite
at every point of $B_0$.
Assuming that $\Gamma=G_{10}^{\#\#}$, the
Small Theorem follows from the Big Theorem,
from Lemma \ref{energy2}, and from the
fact that ${\cal E\/}_{\Gamma}$ is positive
definite throughout $B_0$.

To finish the proof of the Big Theorem
and the Small Theorem, we will show that in each
case  ${\cal E\/}_{\Gamma}$ is positive
definite throughout $B_0$.  This is the same
as saying that all the eigenvalues of the
Hessian $H_{\Gamma}$ are positive at each
point of $B_0$.  We will get this result
by a two step process:
\begin{enumerate}
\item We give an explicit lower bound to the
eigenvalues at $\T \in B_0$.
\item We use Taylor's Theorem with Remainder
to show that the positivity persists throughout
$B_0$.
\end{enumerate}
In Step 1, the point $\T$ corresponds to the point
\begin{equation}
(1,\hskip 10 pt  0,-1/\sqrt 3,\hskip 10 pt 
 -1,0,\hskip 10 pt  0,1/\sqrt 3) \in \R^7.
\end{equation}
The spacing is supposed to indicate how the
coordinates correspond to individual points.
In Step 2 there are two parts.  We have to
compute some higher derivatives of ${\cal E\/}_{\Gamma}$
at the TBP, and we have to bound the remainder term.

For ease of notation, we set $H_k=H_{G_k}$ and
$H_5^{\flat}=H_{G_5^{\flat}}$, etc.

\section{Lower Bounds on the Eigenvalues}
\label{eig00}

Let $H$ be a symmetric $n \times n$ real matrix.
$H$ always has an orthonormal basis of eigenvectors,
and real eigenvalues.  $H$ is {\it positive definite\/}
if all these eigenvalues are positive.  This is equivalent
to the condition that $Hv \cdot v>0$ for all nonzero $v$.
More generally, $Hv \cdot v \geq \lambda \|v\|$, where
$\lambda$ is the lowest eigenvalue of $H$.  

\begin{lemma}[Alternating Criterion]
Suppose $\chi(t)$ is the characteristic polynomial of $H$.
Suppose that the coefficients of $P(t)=\chi(t+\lambda)$ are
alternating and nontrivial.  Then the lowest eigenvalue
of $H$ exceeds $\lambda$.
\end{lemma}

\startproof
An alternating polynomial has no negative roots.
So, if $\chi(t+\lambda)=0$ then $t>0$ and
$t+\lambda>\lambda$.
\endproof

Let us work out $H_3$ in detail.  We first compute
numerically that the lowest eigenvalue at $\T$ is
about $14.0091$. This gives something to aim for.
We then compute, using Mathematica,
that
$$ 64 \det(H_3-(t+14)I)=
198935338432 - 21803700961600 t + 2456945274144 t^2 - $$
 $$  98631799232 t^3 + 1740726332 t^4 - 13792044 t^5 + 49024 t^6 - 
   64 t^7
$$
Since this polynomial is alternating we see (rigorously) that
the lowest eigenvalue of $H_3(\T)$ exceeds $14$.
Similar calculations deal with the remaining energies.
Here are the results (with the first one repeated):
\begin{itemize}
\item the lowest eigenvalue of $H_3(\T)$ exceeds $14$.
\item the lowest eigenvalue of $H_4(\T)$ exceeds $40$..
\item the lowest eigenvalue of $H_5^{\flat}(\T)$ exceeds $53$.
\item the lowest eigenvalue of $H_6(\T)$ exceeds $91$.
\item the lowest eigenvalue of $H_{10}^{\#}(\T)$ exceeds $1249$.
\item the lowest eigenvalue of $H_{10}^{\#\#}(\T)$ exceeds $3145$.
\end{itemize}
In all cases the polynomials we get have dyadic rational
coefficients.

\section{Taylor's Theorem with Remainder}
\label{taylor1}

Before we get to the discussion of the variation
of the eigenvalues, we repackage a special case of
Taylor's Theorem with Remainder.  Here are 
some preliminary definitions.

\begin{itemize}
\item Let $P_0 \in \R^7$ be some point.
\item Let $B$ denote some cube of in-radius
$\epsilon$ centered at $P_0$.
\item $\phi: \R^7 \to \R$ be some function.
\item Let
$\partial_I \phi$ be the
partial derivative of $\phi$ w.r.t. a
multi-index $I=(i_1,...,i_7)$. 
\item  Let
$|I|=i_1+\cdots i_7$.  This is the {\it weight\/} of $I$.
\item Let $I!=i_1! \cdots i_7!$.
\item Let $\Delta^I=x^{i_1}...x^{i_7}$.
Here $\Delta=(x_1,...,x_7)$ is some vector.
\item For each positive integer $N$ let
\begin{equation}
\label{biggie}
M_N(\phi)=\sup_{|I|=N}\sup_{P \in B} |\partial_I \phi(P)|,
\hskip 30 pt
\mu_N(\phi)=\sup_{|I|=N}|\partial_I \phi(P_0)|.
\end{equation}
\end{itemize}

Let $U$ be some open neighborhood of $B$.
Given $P \in B$, let $\Delta=P-P_0$.
Taylor's Theorem with Remainder says that there
is some $c \in (0,1)$ such that
$$
\phi(P)= \sum_{a=0}^N 
\sum_{|I|=a} \frac{|\partial_I \phi(P_0)|}{I!}\Delta^I+
\sum_{|I|=N+1} \frac{\partial_I f(P_0+c\Delta)}{I!}\Delta^I
$$

Using the fact that
$$
|\Delta^I| \leq \epsilon^{|I|}, \hskip 30 pt
\sum_{|I|=m} \frac{1}{I!}=\frac{7^{m}}{m!},
$$
and setting $N=4$ we get
\begin{equation}
\label{taylor}
\sup_{P \in B} |\phi(P)| \leq |\phi(P_0)|+
\sum_{j=1}^4 \frac{(7\epsilon)^j}{j!} \mu_j(\phi) +
\frac{(7\epsilon)^5}{(5)!}M_5(\phi).
\end{equation}

\section{Variation of the Eigenvalues}

Let $H_0$ be some positive definite symmetric
matrix and let $\Delta$ be
some other symmetric matrix of the same size.
Recall various definitions of the $L_2$ matrix norm:
\begin{equation}
\|\Delta\|_2=\sqrt{\sum_{ij}\Delta_{ij}^2}=
\sqrt{{\rm Trace\/}(\Delta\Delta^t)}=
\sup_{\|v\|=1} \|\Delta v\|.
\end{equation}

\begin{lemma}[Variation Criterion]
Suppose that $\|\Delta\|_2 \leq \lambda$, where
$\lambda$ is some number less than the lowest
eigenvalue
of $H_0$.  Then $H=H_0+\Delta$ is also positive definite.
\end{lemma}

\startproof
$H$ is positive definite if and only of
$Hv \cdot v>0$ for every nonzero unit vector $v$.
Let $v$ be such a vector. Writing $v$ in an
orthonormal basis of eigenvectors we see that
$H_0v \cdot v>\lambda$. Hence
$$
Hv \cdot v = (H_0 v+\Delta v) \cdot v  \geq
H_0v \cdot v - |\Delta v \cdot v|>
\lambda - \|\Delta v\| \geq
\lambda - \|\Delta\|_2 \geq 0.
$$
This completes the proof.
\endproof

Let $H$ be any of the Hessians we are considering -- e.g. $H=H_3$.
Define
\begin{equation}
\label{thirdbound}
F=\sqrt{\sum_{|J|=3} M_{J}^2}, \hskip 30 pt
M_{J}=\sup_{P \in B_0} |\partial_J {\cal E\/}(P)|.
\end{equation}
The sum is taken over all multi-indices $J$ of
weight $3$.

Recall that $\epsilon_0=2^{-18}$.

\begin{lemma}
\label{smallL2}
Let $\lambda$ be a strict lower bound
to the smallest eigenvalue of $H(\T)$.
If $\sqrt 7 \epsilon_0 F \leq \lambda$.
Then $H$ is positive
definite throughout $B_0$.
\end{lemma}

\startproof
Let $H_0=H(\T)$.
Let $\Delta=H-H_0$.  Clearly
$H=H_0+\Delta$.

Let $\gamma$ be the unit
speed line segment connecting $P$ to $P_0$ in $\R^7$.
Note that $\gamma \subset B_0$ and
$\gamma$ has length 
$L \leq \sqrt 7 \epsilon_0$.
We set $H_{L}=H$ and we let
$H_t$ be the Hessian of $E$ at the point of
$\gamma$ that is $t$ units from $H_0$.

We have
\begin{equation}
\Delta=\int_0^{L} D_t(H_t)\ dt.
\end{equation}
Here $D_t$ is the unit directional
derivative of $H_t$ along $\gamma$.

Let $(H_t)_{ij}$ denote the $ij$th entry of
$H_t$.  Let $(\gamma_1,...,\gamma_7)$ be the
components of the unit vector in the direction of $\gamma$.
Using the fact that $\sum_k \gamma_k^2=1$ and
the Cauchy-Schwarz inequality, and the
fact that mixed partials commute, we have
\begin{equation}
(D_tH_t)_{ij}^2=
\bigg(\sum_{k=1}^7 \gamma_k 
\frac{\partial}{\partial x_k}
\frac{\partial^2 H_t}{\partial x_i \partial x_j}\bigg)^2 \leq
\sum_{k=1}^7 
\bigg(\frac{\partial^3 H_t}{\partial x_i \partial x_j \partial x_k}\bigg)^2.
\end{equation}
Summing this inequality over $i$ and $j$ we get
\begin{equation}
\|D_tH_t\|_2^2 \leq
\sum_{i,j,k}
\bigg(\frac{\partial^3 H_t}{\partial x_i \partial x_j \partial x_k}\bigg)^2 \leq F^2.
\end{equation}
Hence
\begin{equation}
\label{small2}
\|\Delta\|_2 \leq \int_0^L \|D_t(H_t)\|_2\ dt \leq
L F \leq \sqrt 7 \epsilon_0 F<\lambda.
\end{equation}
This lemma now follows immediately from
Lemma \ref{smallL2}.
\endproof

Referring to Equation \ref{biggie}, and
with respect to the neighborhood $B_0$, define
$M_8({\cal E\/})$ and $\mu_j({\cal E\/})$ for $j=4,5,6,7$. 
Let $J$ be any multi-index of weight $3$. 
Using the fact that
$$
\mu_j(\partial_J {\cal E\/}) \leq \mu_{j+3}({\cal E\/}), \hskip 30 pt
M_{5}(\partial_J {\cal E\/}) \leq M_{8}({\cal E\/}),
$$
we see that Equation \ref{taylor} gives us the
bound
\begin{equation}
\label{BOUND}
M_{J} \leq |\partial_J{\cal E\/}(P_0)| + 
\sum_{j=1}^4
\frac{(7\epsilon_0)^j}{j!}\mu_{j+3}({\cal E\/})
+\frac{(7\epsilon_0)^5}{5!}M_8({\cal E\/}).
\end{equation}

\section{The Biggest Term}
\label{combinatorics}

In this section we will prove that the last term
in Equation \ref{BOUND} is at most $7^5/5!$,
which is in turn less than $141$.  The tiny
choice of $\epsilon_0=2^{-18}$ allows for very
crude calculations.

\begin{equation}
f_k(a,b)=\bigg(4-\|\Sigma^{-1}(a,b)-(0,0,1)\|^2\bigg)^k=
4^k\bigg(\frac{a^2+b^2}{1+a^2+b^2}\bigg)^k.
\end{equation}

$$
g_k(a,b,c,d)=\bigg(4-\|\Sigma^{-1}(a,b)-\Sigma^{-1}(c,d)\|^2\bigg)^k=$$
\begin{equation}
4^k\bigg(\frac{1 + 2 a c + 2 b d + (a^2+b^2)(c^2+d^2) }
{(1 + a^2 + b^2) (1 + c^2 + d^2)}\bigg)^k
\end{equation}
Note that
\begin{equation}
\label{shortcut}
f_k(a,b)=\lim_{c^2+d^2 \to \infty} g_k(a,b,c,d).
\end{equation}
We have
$$
{\cal E\/}_k(x_1,...,x_7)=
f_k(x_1,0)+f_k(x_2,x_3)+f_k(x_4,x_5)+f_k(x_6,x_7)+$$
$$
g_k(x_1,0,x_2,x_3)+
g_k(x_1,0,x_4,x_5)+
g_k(x_1,0,x_6,x_7)+$$
$$
g_k(x_2,x_3,x_4,x_5)+
g_k(x_2,x_3,x_6,x_7)+
g_k(x_4,x_5,x_6,x_7).
$$
Each variable appears in at most $4$
terms, $3$ of which appear in a
$g$-function and $1$ of which appears
in an $f$-function. Hence
\begin{equation}
\label{bound1}
M_8({\cal E\/}_k) \leq M_8(f_k)+3M_8(g_k) \leq 4M_8(g_k).
\end{equation}
The last inequality is a consequence of Equation \ref{shortcut}
and we use it so that we can concentrate on just one
of the two functions above.

\begin{lemma}
\label{easy}
When $r,s,D$ are non-negative integers and $r+s \leq 2D$,
$$\bigg{|}\frac{x^ry^s}{(1+x^2+y^2)^{D}}\bigg{|}<1.$$
\end{lemma}

\startproof
The quantity factors into expressions of the form
$|x^{\alpha}y^{\beta}/(1+x^2+y^2)|$ where
$\alpha+\beta \leq 2$. 
Such quantities are bounded above by $1$.
\endproof

For any polynomial $\Pi$, let $|\Pi|$ denote
the sum of the absolute values of the
coefficience of $\Pi$.
For each $8$th derivative $D_Ig_k$, we have
\begin{equation}
D_Ig_k=\frac{\Pi(a,b,c,d)}{(1+a^2+a^2)^{k+8}(1+c^2+d^2)^{k+8}},
\end{equation}
Where $\Pi_I$ is a polynomial of maximum
$(a,b)$ degree at most $2k+16$ and maximum
$(c,d)$ degree at most $2k+16$.
Lemma \ref{easy} then gives 
\begin{equation}
\sup_{(a,b,c,d) \in \R^4} |D_Ig_k(a,b,c,d)| \leq |\Pi_I|.
\end{equation}

Define
\begin{equation}
\Omega_k=\max_{j=1,...,k} M_8(g_k).
\end{equation}

We compute in Mathematica that

\begin{equation}
\label{bound2}
\Omega_6  \leq
\sup_{k=1,...,6} \sup_I |\Pi_I|=13400293856913653760
\end{equation}
The max is achieved when $k=6$ and $I=(8,0,0,0)$ or
$(0,8,0,0)$, etc.
We also compute that
\begin{equation}
\label{badbound}
   \Omega_{10}  \leq
\sup_{k=1,...,10} \sup_I |\Pi_I|= 162516942801336639946752000
\end{equation}
The max is achieved when $k=10$ and $I=(8,0,0,0)$ or $(0,8,0,0)$, etc.

The quantity
\begin{equation}
4(\Omega_{10}+130 \Omega_6) < 2^{90}
\end{equation}
serves as a grand upper bound to
$M_8(\cal E)$ for all the energies we consider.

This leads to
\begin{equation}
\frac{\bigg(7\times 2^{-18}\bigg)^5}{5!} \times M_8({\cal E\/}) <
7^5/5!<141.
\end{equation}

\section{The Remaining Terms}
\label{exception}

In this section we show that the sum of all the
terms in Equation \ref{BOUND} is less than $750$.

Let
\begin{equation}
\mu^*_{j+3,k}=\frac{(7\epsilon_0)^{j}}{(j)!} \mu_{j+3}({\cal E\/}_{G_k})
\end{equation}
We now estimate
$\mu^*_{j,k}$ for $j=4,5,6,7$ and
$k=2,3,4,5,6,10$.

The same considerations as above show that
\begin{equation}
\mu^*_{j+3,k} \leq \frac{(7\epsilon_0)^{j}}{(j)!}
\bigg(\mu_{j+3}(f_k)+ 3\mu_{j+3}(g_k)\bigg).
\end{equation}
Here we are evaluating the $(j+3)$rd partial
 at all points which
arise in the TBP configuration and then
taking the maximum.  For instance,
for $g_k$, one choice would be
$(a,b,c,d)=(1,0,0,1/\sqrt 3)$.

Again using Mathematica, we compute
\begin{equation}
\sup_{j=5,6,7}\ \sup_{k \leq 10} \mu^*_{j,k}<1/130,
\hskip 30 pt
\sup_{k \leq 6} \mu^*_{4,k}<4,
\hskip 30 pt
\sup_{k \leq 10} \mu^*_{4,k}<60.
\end{equation}
Therefore for any of the eneries we have, we get the bounds
\begin{equation}
\mu^*_4({\cal E\/})<60+4 \times 130=580, \hskip 30 pt
\mu^*_j({\cal E\/})<(1/130) + 130 \times (1/130)<2.
\end{equation}
The sum of all the terms in Equation \ref{BOUND} is
at most
\begin{equation}
580+2+2+2+141<750.
\end{equation}

\section{The End of the Proof}

For any real vector $V=(V_1,...,V_{343})$ define
\begin{equation}
\overline V=(|V_1|+750,...,|V_{343}|+750).
\end{equation}

Given one of our energies, let
$V$ denote the vector of
third partials of ${\cal E\/}$, evaluated at $\T$
and ordered (say) lexicographically.
In view of the bounds in the previous section,
we have
$F \leq \|\overline V\|$. 
Hence
\begin{equation}
\sqrt 7\epsilon_0 F \leq  \sqrt 7\epsilon_0 \|\overline V\|.
\end{equation}

Letting $V_k$ denote the vector corresponding to
$G_k$, we compute
\begin{equation}
\sup_{k=1,...,6} 7\epsilon^0\|V_k\|<1, \hskip 30 pt
\sup_{k=1,...,10} 7\epsilon_0\|V_k\|<60.
\end{equation}

This gives us the following bounds:
\begin{itemize}
\item $7\epsilon_0 F_k<1$ for $k=3,4,5,6$.
\item $7\epsilon_0 F_5^{\flat}<1+25=26$.
\item $7\epsilon_0 F_{10}^{\#}<60+13+68=141$.
\item $7\epsilon_0 F_{10}^{\#}<60+28+102=190$.
\end{itemize}
In all cases, the bound is less than the corresponding
lowest eigenvalue.  This completes our proof that
the Hessian in all cases is positive definite throughout
$B_0$.  Our proofs of the Big Theorem and the Small
Theorem are done.

\newpage
\part{Symmetrization}
\chapter{Preliminaries}

This chapter is something of a grab bag.
The purpose is to gather together the
various ingredients which go into
the proofs of the Symmetrization Lemma.

\section{Monotonicity Lemma}

\begin{lemma}[Monotonicity]
Let $\lambda_1,...,\lambda_n>0$.
If $\beta/\alpha>1$. Then
$$\sum_{i=1}^N \lambda_i^{\beta}-N \geq \frac{\beta}{\alpha}
 \bigg(\sum_{i=1}^N\lambda_i^{\alpha}-N\bigg).$$
\end{lemma}

\startproof
Replacing $\lambda_i$ by $\lambda_i^{\alpha}$ we can
assume without loss of generality that $1=\alpha<\beta$.
Now we have a constrained optimization 
problem.  We are fixing $\sum x_i$ and
we want to minimize $\sum x_i^\beta$ when $\beta>1$.
An easy exercise in Lagrange multipliers shows
that the min occurs when $\lambda_1=...=\lambda_N$.
So it suffices to consider the case $N=1$.

Write $\epsilon=\lambda-1$.  Suppose first that
$\epsilon>0$.  Then
$$\lambda^{\beta}-1=\int_{x=1}^{1+\epsilon} \beta x^{\beta-1} dx>
\int_{x=1}^{1+\epsilon} \beta dx = \beta \epsilon.$$
On the other hand, if $\epsilon<0$ then we have
$1+\epsilon \in (0,1)$ and
$$\lambda^{\beta}-1=-\int_{1+\epsilon}^1 \beta x^{\beta-1}dx>
-\int_{1+\epsilon}^1 \beta dx=\beta\epsilon.$$
That does it.
\endproof

\section{A Three  Step Process}

Our symmetrization is a retraction
from {\bf SMALL\/} to
{\bf K4\/}.  Here we describe it as a $3$
step process.   It is convenient to set $\Omega={\bf SMALL\/}$.
We start with the configuration $X$ having points
$(p_1,p_2,p_3,p_4) \in \Omega$.

\begin{enumerate}

\item We let $(p_1',p_2',p_3',p_4')$ be the
configuration which is obtained by
rotating $X$ about the origin so that
$p_0'$ and $p_2'$ lie on the same horizontal
line, with $p_0'$ lying on the right.
This slight rotation does not change the
energy of the configuration, but it does
slightly change the domain.  While $X \in \Omega$, the
new configuration $X'$ lies in a slightly
modified domain $\Omega'$ which we will describe in
\S \ref{domain2}.  

\item Given a configuration 
$X'=(p_0',p_1',p_2',p_3') \in \Omega'$, there is
a unique configuration 
$X''=(p_0'',p_1'',p_2'',p_3'')$, invariant under
under reflection in the $y$-axis, such that
 $p_j'$ and $p_j''$ lie
on the same horizontal line for $j=0,1,2,3$
and $\|p_0''-p_2''\|=\|p_0'-p_2'\|$.
There is a slightly different domain $\Omega''$
which contains $X''$.  Again, we will
describe $\Omega''$ below.  This is a $6$-dimensional
domain.

\item Given a configuration 
$X''=(p_0'',p_1'',p_2'',p_3'') \in \Omega''$ there
is a unique configuration
$X'''=(p_0''',p_1''',p_2''',p_3''') \in {\bf K4\/}$
such that $p_j''$ and $p_j'''$ lie on the same
vertical line for $j=0,1,2,3$.
The configuration $X'''$ coincides with the
configuration $X^*$ defined in the
Symmetrization Lemma.

\end{enumerate}

To prove the Symmetrization Lemma, it suffices to prove
the following two results.

\begin{lemma}
\label{symm1}
$R_s(X'') \leq R_s(X')$ for all $X' \in \Omega'$ and
$s>2$, with equality iff $X'=X''$.
\end{lemma}

\begin{lemma}
\label{symm2}
$R_s(X''') \leq R_s(X'')$ for all $X'' \in \Omega'$ and
$s \in [12,15+1/2]$, with equality iff  $X''=X'''$.
\end{lemma}

\section{Base and Bows}
\label{base0}

There are $10$ bonds (i.e. distances between pairs of points)
in a $5$ point configuration.  Ultimately, when we
perform the symmetrization operations we want to see that
a sum of $10$ terms decreases.  We find it convenient to
break this sum into $3$ pieces and (for the most part)
study the pieces separately.  We will fix some exponent
$s$ for the power law.  
Given a finite list $V_1,...,V_k$ of vectors,
we define
\begin{equation}
R_s(V_1,...,V_k)=\sum_{i=1}^k \|V_i-V_{i-1}\|^{-s}.
\end{equation}
The indices are taken mod $k$. 

We define the {\it base energy\/} of the configuration
$X$ to be 
\begin{equation}
A_s(X)=R_s(\widehat p_0,\widehat p_1,\widehat p_2,\widehat p_3),
\hskip 30 pt
\widehat p_k=\Sigma^{-1}(p_k).
\end{equation}
Here $\Sigma^{-1}$ is stereographic projection,
as in Equation \ref{inversestereo}.

We define the {\it bow energies\/} to be the two terms
\begin{equation}
B_{02,s}(X)=R_s(\widehat p_0,\widehat p_2,(0,0,1)), \hskip 30 pt
B_{13,s}(X)=R_s(\widehat p_1,\widehat p_3,(0,0,1)).
\end{equation}
We have
\begin{equation}
\label{basebow}
R_s(X)=A_s(X)+B_{02,s}(X)+B_{13,s}(X).
\end{equation}

\begin{center}
\resizebox{!}{1.7in}{\includegraphics{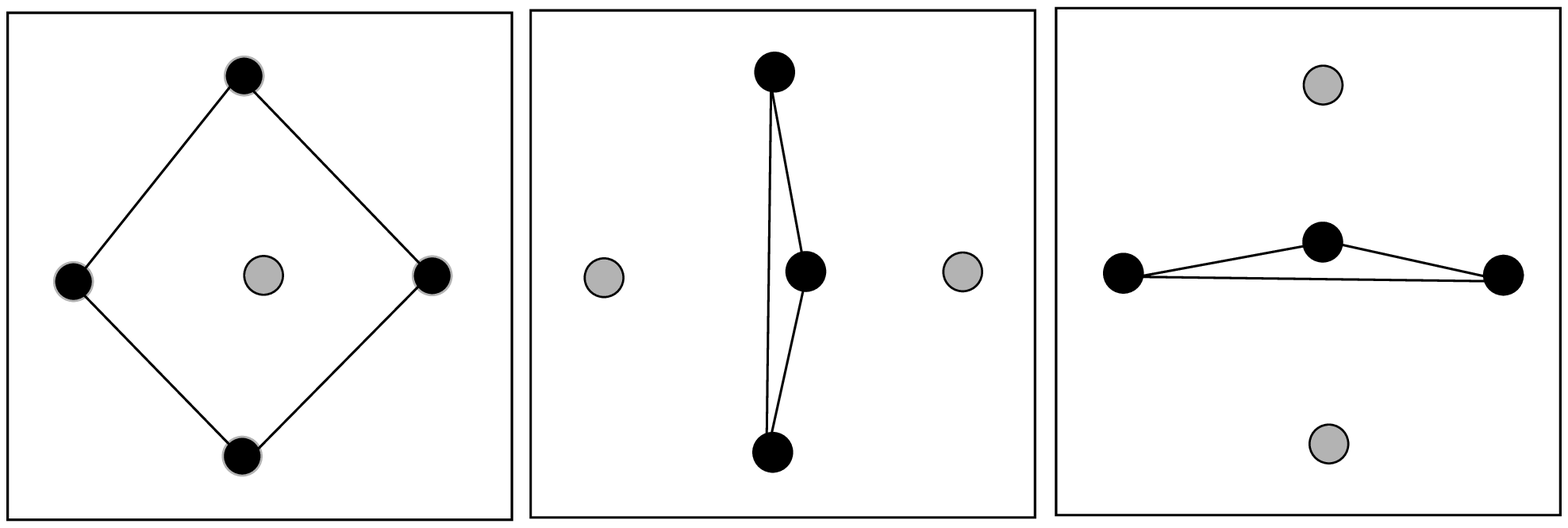}}
\newline
{\bf Figure 12.1:\/} Base and Bows
\end{center}

Our construction above refers to Figure 12.1.
The left hand side of Figure 12.1 illustrates
the bonds that go into the base energy.  The
two figures on the right in Figure 12.1 illustrate
the bonds that go into the two bow energies.
If we interpret Figure 12.1 as displaying points
in the plane, then the 
central dot is really the point at infinity.

\section{Proof Strategy}
\label{proofstrategy}

To prove Lemma \ref{symm1} we will establish
the following stronger results.  These
results are meant to hold for all $s \geq 2$
and all $X' \in \Omega'$.
(Again, we define $\Omega'$ in \S \ref{domain2} below.)
\begin{enumerate}
\item $A_{s}(X'') \leq A_{s}(X')$.
\item $B_{02,s}(X'') \leq A_{02,s}(X')$ with equality
iff $X'=X''$.  
\item $B_{13,s}(X'') \leq A_{13,s}(X')$ with equality
iff $X'=X''$.  
\end{enumerate}
Inequalities 2 and 3 are pretty
easy; they just involve pairs of points in the plane.
Inequality 1 breaks into two 
statements, each involving a triple of points in the plane.
This inequality is considerably harder to prove, and
in fact I do not know a good geometric proof.  However,
notice that the left hand side of the equation
is just the sum of the same term $4$ times.
It follows from the Monotonicity Lemma that
the truth of Inequality 1 at $s=2$ implies
the truth of Inequality 1 for $s>2$.  We
deal with Inequality 1 at $s=2$ with a
direct calculation.

It seems that Inequalities 1-3 also work in the
case of Lemma \ref{symm2}, except that for
Inequality 1 we need to take $s$ fairly large.
I found failures for exponents as high as $s=9$.
So, in principle, we could prove Lemma
\ref{symm2} using a similar strategy, except that
we would pick, say $s=12$, for Inequality 1.
The problem is that the high exponent leads to
enormous polynomials which are too large to
analyze directly.  In the remote future, perhaps such
a calculation would be feasible.

When trying to understand
the failure of Inequality 1 for small
exponents, I noticed the following:
When Inequality 1 fails, Inequalities 2 and 3
hold by a wide margin.  This suggested
a proof strategy.
To prove Lemma \ref{symm2} we will exhibit
constants $C_0$ and $C_1$, which depend on
the configuration and the exponent, such that
\begin{enumerate}
\item $A_s(X'')-A_s(X''') \geq -C_0-C_1$
for all $s \geq 2$,
\item $B_{02,s}(X'')-B_{02,s}(X''') \geq C_0$
for all $s \in [12,15+25/512]$. 
\item $B_{13,s}(X'')-B_{13,s}(X''') \geq C_1$
for all $s \in [12,15+25/512]$.
\end{enumerate}
We get equality in the second and third cases iff $X''=X'''$.
Adding these up we get that $R_s(X'') \geq R_s(X''')$
with equality iff $X''=X'''$.

These inequalities are supposed to hold for
all $X'' \in \Omega''$, a domain we define
in \S \ref{domain2} below.  Moreover, we
get equality in both Inequality 2 and Inequality 3 iff
$X''=X'''$. 
\newline
\newline
{\bf Remark:\/}
Inequalities 2 and 3 are rather delicate, though
part of the delicacy comes from the fact that
I simply squeezed the estimates as hard as I needed
to squeeze them until I got them to work. 
One couldl take the upper bound to be
$15+1/4$ in the Inequality 3 and still have it work.
\newline

In \S \ref{domain2} we will
define the domains $\Omega'$ and $\Omega''$ and
explain how we paramatrize (supersets which cover)
$\Omega'$ and $\Omega''$ by affine cubes.
In \S \ref{symm1proof} we will prove
Lemma \ref{symm1} and
in \S \ref{symm2proof} we will prove
Lemma \ref{symm2}.

\newpage
\chapter{The Domains}
\label{domain2}

\section{Basic Description}

In this section, we define the domains
mentioned above and prove that
they have the desired properties.

Recall that $\Omega$ is the domain of
points $p_0,p_1,p_2,p_3$ such that
\begin{enumerate}
\item $\|p_0\| \geq \|p_k\|$ for $k=1,2,3$.
\item $512 p_0 \in [433,498] \times [0,0]$.
\item $512 p_1 \in [-16,16] \times [-464,-349]$.
\item $512 p_2 \in [-498,-400] \times [0,24]$.
\item $512 p_3 \in [-16,16] \times [349,464]$.
\end{enumerate}
Let $\Omega'$ denote the domain of points
$p_0,p_1,p_2,p_3$ such that $p'_0$ and
$p_2'$ are on the same horizontal line, and
\begin{enumerate}
\item $\|p'_0\| \geq \|p'_k\|$ for $k=1,2,3$.
\item $512p'_0 \in [432,498] \times [0,16]$.
\item $512 p'_1 \in [-16,32] \times [-465,-348]$.
\item $512 p'_2 \in [-498,-400] \times [0,16]$.
\item $512 p'_3 \in [-32,16] \times [348,465]$.
\end{enumerate}

\begin{lemma}
\label{domain1}
If $X \in \Omega$ then $X' \in \Omega'$.
\end{lemma}

\startproof
Rotation about the origin does not change the norms
of the points, so $X'$ satisfies Property 1.
Let $\theta_0 \geq 0$ denote the angle we must rotate in order
to arrange that $p_0'$ and $p_2'$ must lie on the same
horizontal line.  This angle maximized when $p_0$ is
an endpoint of its segment of constraint and $p_2$ is
one of the two upper vertices of rectangle of constaint.
Not thinking too hard which of the $4$ possibilities
actually realizes the max, we check for all
$4$ pairs $(p_0,p_2)$ that
$$\bigg(R_{1/34}(p_0)\bigg)_2>
\bigg(R_{1/34}(p_2)\bigg)_2.$$
Here $(p)_2$ denotes the second coordinate of $p$
and $R_{\theta}$ denotes counterclockwise rotation
by $\theta$.
From this we conclude that $\theta_0<1/34$.
Now we need to see that if we rotate by less
than $1/34$ radians, the points belong to
$\Omega'$.
\newline
\newline
{\bf Case 1:\/}
Since $p_0$ lies in the positive $x$-axis and
$1-\cos(1/34)<1/512$, we have
$p_{01}' \in [p_{01}-1/512,p_{01}]$.
Since $\sin(1/34)<1/32$ we have
$p'_{02} \in [0,16/512]$.  These give
the constraints on $p_0'$ in Property 2.
\newline
\newline
{\bf Case 2:\/}
When we apply $R_{\theta}$, the first coordinate
of $p_1$ increases, and does so by at most
$\sin(1/34)<1/32$.  This gives the constraint
on $p_{11}'$ in Property 3.
Let $\theta_1$ denote the maximum angle that
$p_1$ can make with the $y$-axis.  
This angle is maximized when
$$
p_1=\bigg(\frac{\pm 16}{512},\frac{349}{512}\bigg).
$$
Since $\tan(1/21)>16/349$ we conclude that
$\theta_1<1/21$.
When we rotate, $p_{12}$ can change by at most
$$\cos\bigg(\frac{1}{21}+\frac{1}{34}\bigg)-
\cos\bigg(\frac{1}{21}\bigg)<\frac{1}{512}.$$
This gives the constraints on $p'_{12}$
in Property 3.
\newline
\newline
{\bf Case 3:\/} 
We have $p_{22}'=p_{02}'$ and 
$p_{21} \leq p_{21}' \leq \|p_0\|.$
This gives the constraints on $p_2'$ in Property 4.
\newline
\newline
{\bf Case 4:\/} This follows from Case 2 and symmetry.
\endproof

Let $\Omega''$ denote the domain of configurations
$X''$ which are invariant under reflection in
the $y$-axis and
\begin{enumerate}
\item $512 p''_{01} \in [416,498]$
\item $512 p''_{02} \in [0,16]$.
\item $512 p''_{12} \in [-465,-348]$.
\item $512 p''_{32} \in [348,465]$.
\end{enumerate}
Note that we also have
$$
p_{21}''=-p_{01}'', \hskip 30 pt
p_{22}''=p_{02}'', \hskip 30 pt
p_{11}''=p_{31}''=0.
$$
The domain $\Omega''$ is $4$ dimensional.
It follows almost immediately from the
definition of the map that $X' \in \Omega'$
implies $X'' \in \Omega''$.  The lower bound
of $416/512$ for $p_{01}'$ comes from
averaging the lower bound $(432/512)$ for $p_{01}$ with
the lower bound $(400/512)$ for $-p_{21}$.

Both $\Omega'$ and $\Omega''$ are polyhedral subsets
of $\R^8$, respectively $7$ and $6$ dimensional.
It turns out that we will really only need to
make computations on certain projections of these
domains.  That is, we can forget about the location
of one of the points at stage and concentrate
on just $3$ points at a time.  We now describe
these projections.
 
Let $\Omega'_1$ denote the domain of triples $(p'_0,p'_2,p_3')$
arising from configurations in $\Omega'$.  In other
words, we are just forgetting the location of $p'_1$.
Formally, $\Omega'_1$ is the image of $\Omega'$ under
the forgetful map  $(\R^2)^4 \to (\R^2)^3$.
We define $\Omega_3'$ similarly, with
the roles of $p'_1$ and $p'_3$ reversed.  Likewise, we
define $\Omega''_0$ to be the set of triples
$(p''_2,p''_1,p''_3)$ coming from quadruplies in $\Omega''$.
Finally, we define $\Omega''_2$ by switching the roles
played by $p''_0$ and $p''_2$.  The
spaces $\Omega'_1$ and $\Omega'_3$ are each $5$
dimensional and the spaces $\Omega''_0$ and $\Omega''_2$
are each $4$ dimensional.

\section{Affine Cubical Coverings}
\label{cube}

In this section we describe how to parametrize
supersets which contain the domains
$\Omega'_1$, $\Omega'_3$, $\Omega''_0$, and
$\Omega''_2$.

Let $S$ be a subset of Euclidean space, $\R^N$.
We say that an {\it affine cubical covering\/} of $S$
if a finite collection $(Q_1,\phi_1),...,(Q_m,\phi_m)$
where \begin{itemize}
\item $Q_k=[0,1]^n$ for each $k=1,...,m$.
\item $\phi_k: Q_k \to \R^N$ is an affine map for each $k=1,...m$.
\item $S$ is contained in the union $\bigcup_{k=1}^m \phi_k(Q_k)$.
\end{itemize}
In all $4$ cases we consider, we will have $N=6$ and $m=2$.
When $S=\Omega_1'$ or $\Omega_3'$ we will have $n=5$.
When $S=\Omega''_0$ or $S=\Omega''_2$ we
will have $m=4$.  Our maps will have the additional
virtue that they will be defined over $\Q$.  This
will make it so that our final calculations are
integer calculations.

Suppose we have some polynomial $f: \R^N \to \R$ and we
want to see that the restriction of $f$ to $S$ is non-negative.  
Our strategy is to find a cubical parametrization
of $S$ and to show that each of the functions
\begin{equation}
f_k=f \circ \phi_k: Q_k \to \R
\end{equation}
is non-negative.
This gives us a collection
$f_1,...,f_k$ which we want show are positive 
(respectively non-negative) on the unit cube.
For this purpose we use either the positive
dominance criterion or the positive dominance
algorithm.

\section{The Coverings for Lemma 1}

In this section we describe the cubical coverings for
$\Omega_1'$ and $\Omega_3'$.  In each case, we
will use $2$ affine cubes.  Thus, we have to
describe $4$ affine cubes all in all.

For rational numbers $r_1=a_1/b_1$ and $r_2=a_2/b_2$ we define
\begin{equation}
[r_1,r_2,t]=r_1(1-t)+r_2t.
\end{equation}
Thus function linearly interpolates between
$r_1$ and $r_2$ as $t$ interpolates
between $0$ and $1$.
We will always use rationals of the form $p/512$, so
we write $r_j=a_j/512$ and we store the information
for the above function as the pair of integers
$(a_1,a_2)$.

An array of the form
\begin{equation}
\left[\matrix{\ell_{11},\ell_{12} \cr
\cdots \cr
\ell_{n1},\ell_{n2}}\right]
\end{equation}
denotes the map
\begin{equation}
\label{linfun}
f(t_1,...,t_n)=([\ell_{11}/512,\ell_{12}/512,t_1],...,[\ell_{n1}/512,\ell_{n2}/512,t_n]).
\end{equation}

Here is our cubical covering for
$\Omega_1'$.
\begin{equation}
\left[\matrix{498&432\cr  0&16\cr  -32&16\cr 348&465\cr 0&64}\right] 
\hskip 30 pt 
\left[\matrix{432&416\cr 0&16\cr -32&16\cr 348&465\cr 0&64}\right]
\end{equation}
Now we describe the meaning of the variables.
We let the first coordinate of the map
in Equation \ref{linfun} be $a$ and the second $b$, and so on.
Then
\begin{itemize}
\item $p_{01}'=a+e$ and $p_{01}''=a$.
\item $p_{02}'=p_{02}''=p_{22}'=p_{22}''=b$.
\item $p_{21}'=-a+e$ and $p_{21}''=-a$.
\item $p_{31}'=c$ and $p_{31}''=0$.
\item $p_{32}'=p_{32}''=d$.
\end{itemize}
Here $p'_0=(p_{01},p'_{02})$, etc.
Thus, for example,
$$p_{01}'=\frac{498(1-t_1)+432 t_1 + 64 t_5}{512},
\hskip 30 pt t_1,t_5 \in [0,1].$$

We describe why the cubical covering for 
$\Omega_1'$ really is a covering. 
Here are the relevant conditions for $\Omega_1'$.
\begin{enumerate}
\item $512p'_0 \in [432,498] \times [0,16]$.
\item $512 p'_2 \in [-498,-400] \times [0,16]$.
\item $512 p'_3 \in [-32,16] \times [348,465]$.
\item $p_{22}'=p_{02}'$.
\item $|p_{21}| \leq p_{01}$. 
\end{enumerate}
We ignore the other conditions because we only
care about getting a covering.

A comparison between the terms in our matrix and
the conditions for $\Omega'_1$ reveals that the
numbers are exactly the same for all the variables
except $p_{01}'$ and $p_{21}'$.  Our parameterization
requires an explanation in this case.  Figure 13.1
shows a plot of the values of $(p_{01}',p_{21}')$ we
pick up with our parametrization.  The thickly
drawn rectangle is the domain of interest to us
without the condition $|p_{21}| \leq p_{01}$.
When we add this condition we remove the
darkly shaded triangular region.  The lightly
shaded regions show the projections of our
two affine cubes.  The matrix entry $64$ above
is somewhat arbitrary, and making this number
smaller would cause the shaded images to retract
a bit.  The labels in Figure 13.1 are scaled up by
$512$.

\begin{center}
\resizebox{!}{2.5in}{\includegraphics{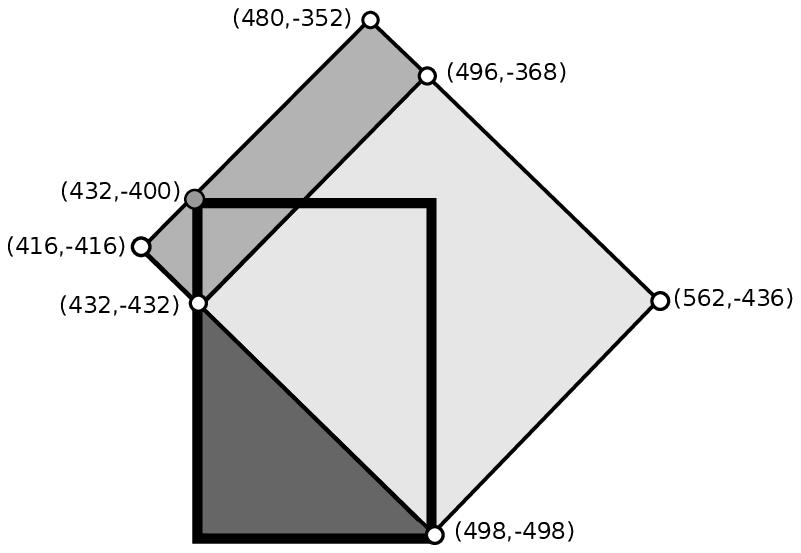}}
\newline
{\bf Figure 13.1:\/} Covering $\Omega_1'$ with affine cubes.
\end{center}

\noindent
{\bf Remark:\/}
It would have been nicer if we could have used a
cubical covering with just one affine cube.  We tried
this, but the resulting polynomial we got was not
positive dominant.  It is very likely that the
positive dominance algorithm would work on this
polynomial, but the polynomial is too big for the
P.D.A. to run quickly enough.  That is why we split
things up into two cubes.
\newline

Our cubical covering for $\Omega_3'$ is almost the same.
We put in bold the numbers which have changed.
\begin{equation}
\left[\matrix{498&432\cr  0&16\cr  -{\bf 16\/}&{\bf 32\/}\cr {\bf -348\/}&{\bf -465\/}\cr 0&64}\right] 
\hskip 30 pt 
\left[\matrix{432&416\cr 0&16\cr -{\bf 16\/}&{\bf 32\/}\cr {\bf -348\/}&{\bf -465\/}\cr 0&64}\right]
\end{equation}
The variables have the same meanings, except that we use the point
$p_1$ in place of the point $p_3$.  We get a covering in this case
for the same reason as in the previous case.

\section{The Coverings for Lemma 2}
\label{cube2}

In this case we only need to deal with
$\Omega_2''$.  The result for $\Omega_0''$
follows from reflection symmetry (in the $y$-axis.) 
We cover $\Omega_2''$ using $4$ affine cubes.
The cubes come in pairs, and the two cubes
in a pair differ only in the sign of a certain
entry. The sign determines whether we are
parametrizing configurations where $p_{12}''+p_{32}''$
is positive or negative.  Here are the $4$ cubes:

\begin{equation}
\left[\matrix{416&498\cr  0&16\cr  348&465\cr 0&\pm 24}\right],
\hskip 30 pt
\left[\matrix{416&498\cr  0&16\cr  364&449\cr 0&\pm 64}\right].
\end{equation}

Here are the meanings of the variables.
\begin{itemize}
\item $p_{01}''=p_{01}'''=a$.
\item $p_{02}''=b$ and $p_{02}'''=0$.
\item $p_{11}''=p_{31}''=p_{11}'''=p_{31}'''=0$.
\item $p_{12}''=-c+d$ and $p_{12}'''=-c$;
\item $p_{32}''=c+d$ and $p_{12}'''=c$;
\end{itemize}

Now we explain why this really gives a covering.
The conditions for $\Omega_2''$ are:
\begin{enumerate}
\item $512 p''_{01} \in [416,498]$
\item $512 p''_{02} \in [0,16]$.
\item $512 p''_{12} \in [-465,-348]$.
\item $512 p''_{32} \in [348,465]$.
\end{enumerate}
The only point about the covering 
that is not straightforward is why all possible
coordinate pairs $(p_{12}'',p_{32}'')$ are covered.
Figure 13.2 shows the picture of how this works.
The thick square is the domain of pairs we want
to cover and the shaded region shows that we
actually do cover.  Again, the labels are
scaled up by $512$.

\begin{center}
\resizebox{!}{4in}{\includegraphics{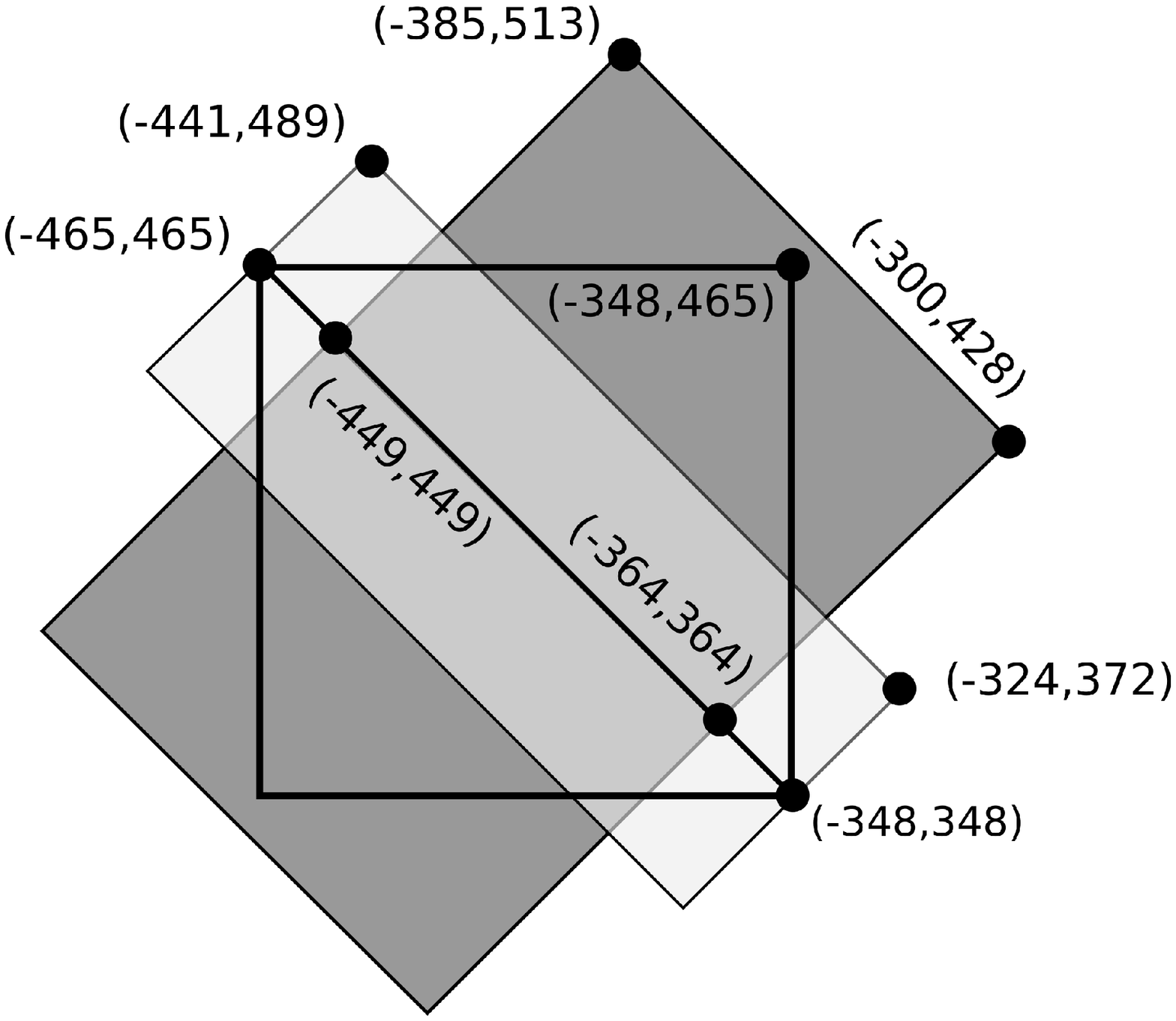}}
\newline
{\bf Figure 13.2:\/} Covering $\Omega''$ with affine cubes.
\end{center}

We have not labeled every point in the picture because
the whole picture is invariant under the map $(x,y) \to (-y,-x)$.
So, the remaining labels can be deduced from symmetry.

\newpage

\chapter{The First Symmetrization}
\label{symm1proof}

In this chapter we prove Theorem \ref{symm1} by
establishing Inequalities 1-3 listed
in \S \ref{proofstrategy} concerning
the configurations $X'$ and $X''$.

\section{Inequality 1}
\label{proofineq1}

We want to show that
\begin{equation}
A_s(X'') \leq A_s(X')
\end{equation}
for any $X' \in \Omega'$ and any $s \geq 2$. 
It follows from the Monotonicity Principle
that the case $s=2$ suffices.
We will split the base energy into two pieces
and handle each one separately.  We write,
for $k=1,3$,
\begin{equation}
A_{k,2}(X')=\|\Sigma^{-1}(p_0)-\Sigma^{-1}(p_{4-k})\|^{-2}+
\|\Sigma^{-1}(p_2)-\Sigma^{-1}(p_{4-k})\|^{-2}.
\end{equation}
In other words, we are omitting the interactions
with $p_k'$ when we define $A_{k,2}$.  The domain
for $A_{k,2}$ is $\Omega'_k$, described in
\S \ref{domain2}.
It clearly suffices to show that
\begin{equation}
A_{k,2}(X'') \leq A_{k,2}(X'), \hskip 30 pt k=1,3.
\end{equation}

What we now say works for either $k=1$ or $k=3$.
Let $(\phi_1,Q_1)$ and $(\phi_2,Q_2)$ be
the cubical covering of $\Omega'_k$ described
in \S \ref{cube}.  Again, what we say works
for either of these cubes.  We consider
the rational function
\begin{equation}
\label{mainfunction}
\frac{P(t_1,...,t_5)}{Q(t_1,...,t_5)}=(A_{k,2}(X')-A_{k,2}(X'')) \circ \phi.
\end{equation}
The polynomials are in $\Z[t_1,...,t_5]$,
the expression is reduced, and
$Q$ is chosen to be positive at
$(1/2,1/2,1/2,1/2,1/2)$.  This determines
$P$ and $Q$ uniquely.

We would like to show directly that this rational
function is positive on $[0,1]^5$, but we don't know how to
do this directly.   We take a different approach.

The two most interesting variables are $t_3$ and
$t_5$.  The first of these variables controls
the $x$-coordinate of $p'_3$ and the second of
these controls the deviation from the points
$p_0'$ and $p_2'$ being symmetrically placed about
the $y$-axis.  Setting $t_3=t_5=0$
is our first symmetrization operation.
We consider the Hessian with respect to these
two variables:

\begin{equation}
H_P=\left[\matrix{\partial_{33}P & \partial_{35}P \cr
\partial_{53}P &  \partial_{55} P} \right].
\end{equation}
Here $\partial_{33}P=\partial^2 P/\partial_{t_3}^2$, etc.
This is another polynomial on $[0,1]^5$.
Here is the calculation we can make.

\begin{theorem}
\label{comp1}
The polynomials $\det(H_P)$ and ${\rm trace\/}(H_P)$
are weak positive dominant on $[0,1]^5$.  Hence,
these functions are positive on $(0,1)^5$.
\end{theorem}

All in all, Theorem \ref{comp1} involves $8=2 \times 2 \times 2$
polynomials. There are $2$ domains, $2$ affine cubes, and
$2$ quantities to test.  Everything in sight is an
integer polynonial.  I did the calculation first in
Mathematica, and then reprogrammed the whole thing in
Java and did it again.

As is well known, a symmetric
$2 \times 2$ matrix is positive definite -- i.e.
has positive eigenvalues -- if and only if
its determinant and trace are positive.
Therefore, for all $4$ choices, we have:

\begin{corollary}
The Hessian $H_P$ is positive definite on $(0,1)^5$.
\end{corollary} 

Just so that we can make choices, we consider the
case $k=1$.  The choice $t_3=2/3$ corresponds
to $p_3''=p_3'$ and the choice $t_5=0$
corresponds to $p_0'=p_0''$ and $p_2'=p_2''$.
Let $\Pi$ be the codimension $2$ plane corresponding
to $t_3=2/3$ and $t_5=0$.

\begin{lemma}
$P$ and $\partial_3 P$ and $\partial_5 P$ all
vanish identically on $\Pi$.
\end{lemma}

\startproof
$P$ vanishes by definition.  $\partial_3 P$ and
$\partial_5 P$ vanish by symmetry, but just to
be sure we calculated symbolically that they do
in fact vanish on $\Pi$.
\endproof

Now we come to the main point of all these calculations.

\begin{lemma}
\label{linesegment}
$P>0$ on $(0,1)^5- \Pi$.
\end{lemma}

\startproof
One way to see this is to restrict $P$ to
a straight line segment $L$ that joins an arbitrary point
$\zeta \in (0,1)^5-\Pi$ (corresponding to a configuration
$X' \not = X''$) to the point $\zeta_0 \in \Pi$ corresponding
to $X''$.  Note that only the variables $t_3$ and
$t_5$ change along $L$.
Let $f$ be the restriction of $P$ to $L$.
Since $H_P$ is positive definite we
see that $f$ is a convex function.
Moreover, our calculation above shows that $f(\zeta_0)=0$
and $f'(\zeta_0)=0$.  But then we must have
$f'>0$ at almost all points of $L$.  Integrating,
we get $f(\zeta)>0$.
\endproof

Remember that $P$ is the numerator of a rational function $P/Q$.

\begin{lemma}
\label{Qpos}
$Q>0$ on $(0,1)^5-\Pi$.
\end{lemma}

\startproof
We sample one point in the connected
domain $(0,1)^5-\Pi$ to check that $Q>0$ somewhere in this set.
If $Q$ ever vanishes, then the fact that $P>0$ would cause
our well-defined rational function to blow up.  This does not happen
because the rational function is well defined.
Hence $Q>0$ everywhere on $(0,1)^5-\Pi$.
\endproof

Combining the last two lemmas, we see that the function
in Equation \ref{mainfunction} is positive on $(0,1)^5-\Pi$.
But then, by definition
$$A_{k,2}(X'')<A_{1,2}(X')$$ when $X' \not = X''$.  This is
what we wanted to prove in case $k=1$.  

The case
$k=3$ has the same proof except that we take
$t_3=1/3$ rather than $t_3=2/3$,

This completes the proof of Inequality 1.

\section{Inequality 2}

In this section we prove Inequality 2. That is,
\begin{equation}
B_{02,s}(X'') \leq B_{02,s}(X'), \hskip 30 pt
\forall s \geq 2,
\end{equation}
with equality iff the relevant points coincide.

The relevant points here are
\begin{equation}
\label{hor1}
p'_0=(x+d,y), \hskip 15 pt
p'_2=(-x+d,y), \hskip 15 pt
p''_0=(x,y), \hskip 15 pt
p''_2=(-x,y).
\end{equation}

The way our proof works is that we further
break the $3$-term expression for the bow
energy into a sum of an expression with
one term and an expression with two terms.
We show that separately these sums decrease
when $s=2$, and then we appeal the Monotonicity
Lemma proved in the previous chapter to take
care of the case $s>2$.  Along the way, we
reveal a beautiful formula.  Here are the details.

We compute
$$
\|\Sigma^{-1}(p'_0)-\Sigma^{-1}(p'_2)\|^{-2}-
\|\Sigma^{-1}(p''_0)-\Sigma^{-1}(p''_2)\|^{-2}=
\frac{d^2}{x^2} \times (2+d^2-2x^2+2y^2).
$$
Since $2-2x^2+2y^2>0$ we see that this expression
is always positive.  This shows that one of the
three terms in the bow energy decreases, when $s=2$.
When $s>2$ we get the same energy decrease, by the
Monotonicity Lemma from the previous chapter.

For the other two terms, we note the following
beautiful formula.
$$
\bigg(\|\Sigma^{-1}(p'_0)-(0,0,1)\|^{-2}+
  \|\Sigma^{-1}(p'_2)-(0,0,1)\|^{-2}\bigg)-$$
\begin{equation}
\label{beauty}
  \bigg(\|\Sigma^{-1}(p''_0)-(0,0,1)\|^{-2}+
  \|\Sigma^{-1}(p''_2)-(0,0,1)\|^{-2}\bigg)= \frac{d^2}{2},
\end{equation}
This holds identically for all choices of $x,y,d$.
Thus, the sum of the other two terms decreases at
the exponent $s=2$.  But then the Monotonicity
Lemma gives the same result for all $s>2$.
with equality iff $d=0$.

\section{Results about Triangles}

Before we prove Inequality 3, we need a
technical lemma about triangles.
The first two results in this section are surely
known, but I don't have a reference.  Lacking a
reference, I will prove them from scratch.
The third result is probably too obscure to
be known, but we need it for Inequality 3.
The sole purpose of this section is to prove
the technical result at the end.  This result
feeds into the proof of Lemma \ref{rotate}
and is not used elsewhere.

Let $C$ be
the unit circle and let $a,b$ be distinct
points in $C$.  Let $\widehat C$ denote one
of the two arcs bounded by $a$ and $b$.
For convenience we will draw pictures in
the case when $\widehat C$ is more than
a semicircle, but the result works in
all cases. 

\begin{lemma}
The function $f(c)=\|a-c\|+\|b-c\|$ has a unique
critical point at the midpoint of $\widehat C$,
and this point is a maximum. In other words,
the function $c \to f(c)$ increases as $c$
moves from $a$ to the midpoint of $\widehat C$.
\end{lemma}

\startproof
Our proof refers to Figure 14.1 below.  By symmetry
it suffices to consider the case when $c$ is
closer to $a$.  We will show that $f(c)$ increases
as $c$ moves further away from $a$. Let $V_a$ and
$V_b$ be unit vectors which are based at $c$ and
point from $a$ to $c$ and $b$ to $c$ respectively.
Let $T_c$ denote the tangent vector to $C$ at $c$.
The desired monotonicity follows from
the claim that
\begin{equation}
T \cdot (V_a+V_b)>0.
\end{equation}

\begin{center}
\resizebox{!}{2.6in}{\includegraphics{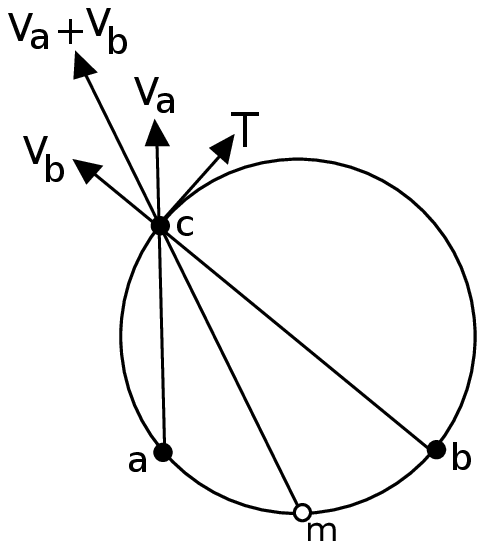}}
\newline
{\bf Figure 14.1:\/} The Circle Result
\end{center}

Let $m$ be the midpoint of the other arc of $C$
bounded by $a$ and $b$.
The vector $V_a+V_b$ makes an equal angle with both
$V_a$ and $V_b$.  But then, by a familiar theorem
from high school geometry, $V_a+V_b$ is parallel to
the ray connecting the point $m$ to $v$.  But this
ray obviously has positive dot product with $T$.
The result follows immediately.
\endproof

Here is a corollary of the preceding result.

\begin{lemma}
\label{monotone}
For any $s>0$ the function
$f(c)=\|a-c\|^{-s}+\|b-c\|^{-s}$ has a unique
critical point at the midpoint of $\widehat C$,
and this point is a minimum. In other words,
the function $c \to f(c)$ decreases as $c$
moves from $a$ to the midpoint of $\widehat C$.
\end{lemma}

\startproof
The key property we use is that the power law functions
are convex decreasing.  We use the same set up
as in the previous lemma.  Suppose we move $c$ to a
nearby location $c'$ further away from $a$.
There are constants $\epsilon_a>0$ and $\epsilon_b$ such that
\begin{equation}
\|a-c'\|=\|a-c\|+\epsilon_a, \hskip 30 pt
\|b-c'\|=\|b-c\|-\epsilon_b.
\end{equation}
By inspection we have $\epsilon_a>0$.  If
$\epsilon_b<0$ then both distances have increased
and the result of this lemma is obvious. 
So, we just have to worry about the case when
$\epsilon_b>0$.
By the previous result, we have
$\epsilon_a>\epsilon_b$.

We set $\epsilon=\epsilon_a$.  
It now follows from monotonicity that
\begin{equation}
\label{key00}
f(c')<(\|a-c\|+\epsilon)^{-s}+(\|b-c\|-\epsilon)^{-s}.
\end{equation}
We also have $\|a-c\|<\|b-c\|$.  Finally, we know that
the power function $R_s(x)=x^{-s}$ is convex decreasing.
This is to say that the derivative of $R_s$ is negative
and increasing.  Hence, when
$\epsilon$ is small, the first term on the right
side of Equation \ref{key00} decreases more
quickly than the second term increases.  This gives
$f(c')<f(c)$ for $\epsilon$ sufficiently small.
But then this is true independent of
$\epsilon$ and remains true until $\|a-c\|=\|b-c\|$.
\endproof

Now we come to the useful technical lemma.
As usual, let $\Sigma^{-1}$ be inverse stereographic
projection.
Call the pair of points
$((0,t),(0,u))$ a {\it good pair\/} if
$-t,u < \sqrt 3/3$.
The significance of the lower bound is
that the points
$\Sigma^{-1}(0,\pm 1/\sqrt 3)$ form an
equilateral triangle with $(0,0,1)$.

Let ${\cal E\/}_s(t,u)$ be the bow energy
associated to this pair.  That is
\begin{equation}
\label{term123}
{\cal E\/}_s(t,u)=
\matrix{
\|\Sigma^{-1}(0,-t)-(0,0,1)\|^{-s}+ \cr
\|\Sigma^{-1}(0,u)-(0,0,1)\|^{-s}+ \cr
\|\Sigma^{-1}(0,-t)-\Sigma^{-1}(0,u)\|^{-s}}
\end{equation}
For later reference, we will call the
three summands {\it term 1\/},
{\it term 2\/}, and {\it term 3\/}.

\begin{lemma}
\label{horincrease}
If $(u^*,t^*)$ is also a good pair with
$(-t^*) \geq (-t)$ and $u^* \geq u$ then
${\cal E\/}_s(t^*,u^*) \geq {\cal E\/}_s(t,u)$.
\end{lemma}

\startproof
Since we can switch the roles of $t$ and $u$, it
suffices to consider the case when $t=t^*$
and just $u$ moves.
Suppressing the dependence on $s$, we consider
the function
\begin{equation}
\phi(u^*)={\cal E\/}_s(t,u^*).
\end{equation}
We will look at the function geometrically.
By symmetry, this function has an extremum
when $U=\Sigma^{-1}(0,u^*)$ forms an isosceles
triangle with the points
$N=(0,0,1)$ and $T=\Sigma^{-1}(0-t)$.
Figure 14.2 shows the situation.

\begin{center}
\resizebox{!}{3.5in}{\includegraphics{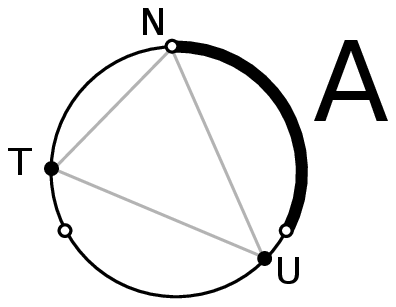}}
\newline
{\bf Figure 14.2:\/} Points on the a great circle
\end{center}

The thickened arc $A$ indicates the allowable locations
$U$ when $(t^*,u^*)$ is a good pair.  Significantly,
the apex of the isosceles triangle we have been
talking about lies outside $A$ when we have a
good pair.  By Lemma \ref{monotone}, the function
$\phi$ is monotone increasing on the arc $A$ as
we move away from the point $U$ to the point $N$.
This is what we wanted to prove.
\endproof

\section{Inequality 3}

In this section we prove Inequality 3.  That is,
\begin{equation}
\label{ver1}
B_{13,s}(X'') \leq B_{13,s}(X'), \hskip 30 pt
\forall s \geq 2,
\end{equation}
with equality iff the relevant points coincide.

The relevant points here are
$$
p'_1=(x_1,-y+d), \hskip 30 pt
p'_3=(x_3,y+d),$$
\begin{equation}
p''_1=(0,-y+d), \hskip 15 pt
p''_3=(0,y+d).
\end{equation}
To get from $(p_1',p_3')$ to $(p_1'',p_3'')$ we
are just sliding these points along horizontal
lines until they hit the $y$-axis.  Algebraically,
this sliding amounts to setting $x_1=x_3=0$.

Reflecting in the $y$-axis, we can assume without
loss of generality that $x_1 \geq 0$.  Suppose that
$x_3>0$.  What this means is that the points
$\widehat p_1$ and $\widehat p_3$ lie in the same
hemisphere on $S^2$.  We can increase their
chordal distance by reflecting one of the points
in the boundary of this hemisphere. In other words,
if we replace $x_3$ by $-x_3$ then we decrease
$A_{s,13}(X')$.  So, it suffices to consider
the case when $x_1 \geq 0$ and $x_3 \leq 0$.
Here is a critical special case.

\begin{lemma}
Equation \ref{ver1} holds when $x_1=0$ and $x_3<0$.
\end{lemma}

\startproof
We need to analyze what happens to the
$3$ terms in Equation \ref{term123}.
We have $p_1'=p_1''$ so 
$$\|\Sigma^{-1}(p_1')-(0,0,1)\|^{-s}=
\|\Sigma^{-1}(p_1'')-(0,0,1)\|^{-s}.$$
This takes care of term 1.

Now we deal with term 2.
Note that
$p_3''$ is closer to the origin than $p_3'$ so
$\Sigma^{-1}(p_3'')$ is farther from $(0,0,1)$ than
is $\Sigma^{-1}(p_3')$.  Hence 
\begin{equation}
\|\Sigma^{-1}(p_3'')-(0,0,1)\|^{-s}<
\|\Sigma^{-1}(p_3')-(0,0,1)\|^{-s}.
\end{equation}
Hence term 2 decreases.

Now we deal with term 3.
By the 
Monotonicity Lemma, it suffices to consider
the case when $s=2$.  Holding
$y_1$ and $y_3$ fixed and setting $x_3=x$, we compute
$$
\|\Sigma^{-1}(0,y_1),\Sigma^{-1}(x,y_3)\|^{-2}=
\frac{(1+y_1^2) u(x)}{4}, \hskip 10 pt
u(x)=\frac{x^2+(y_3^2+1)}{x^2+(y_1-y_3)^2}.
$$
To finish the proof it suffices to prove that
$u(x)$ takes its min at $x=0$.
Once we
know this, we see that term 3 decreases
when $s=2$. The Monotonicity Lemma then implies
that this happens for $s>2$.

To analyze $u$, we write
\begin{equation}
u(x)=\frac{x^2+C_1}{x^2+C_2}, \hskip 30 pt
C_1=y_3^2+1, \hskip 30 pt
C_2=(y_1-y_3)^2.
\end{equation}
This function has a unique extremum at $0$.
To find out if this extremum is a min or a max,
we take the second derivative:
\begin{equation}
\frac{d^2u}{dx^2}=
\frac{2(C_2-C_1)(C_2-3x^2)}{(C_2+x^2)^2}.
\end{equation}
We are hoping a min, so we want to verify that
the second derivative is positive.
The denominator is obviously positive.

To show the positivity of the first factor in
the numerator, we just
use $y_1<-\sqrt 3/3$ and
$y_3>\sqrt 3/3$.  This gives
$$y_1^2>1/3, \hskip 30 pt -2y_1y_3>2/3,
\hskip 30 pt
(y_1-y_3)^2>4/3.$$
Hence
\begin{equation}
C_2-C_1=-1+y_1^2-2y_1y_3>-1+(1/3)+(2/3)=0.
\end{equation}
\begin{equation}
C_2-3x^2>(2 \sqrt 3/3)^2-3(1/16)^2>(4/3)-(1/3)=1.
\end{equation}
This takes care of term 3.

In summary, term 1 in Equation \ref{term123}
does not change, and (unless we have equality in the
configurations) term 2 decreases and term 3 decreases.
\endproof

\noindent
{\bf Remark:\/} It is easy to make a sign error in
a calculation like the one just done.  So, to be sure, I plotted 
the graph of the function $u(x)$ for random choices
of $y_1$ and $y_3$ and say that it had a local min
at $x=0$.
\newline

The case when $x_1>0$ and $x_3=0$ has the
same treatment as the case just considered,
and indeed follows from this case and symmetry.
So, we just have to consider the case 
when $x_1>0$ and $x_3<0$.
In this case, the idea is to tilt the plane
and see what happens.  Since there are
no counterexamples when one of $x_1$ or $x_3$
equals $0$, the following result finishes our proof.

\begin{lemma}
\label{rotate}
If Equation \ref{ver1} has a counterexample
with $x_1>0$ and $x_3<0$ it also has a
counterexample when one of $x_1$ or $x_3$
equals $0$.
\end{lemma}

\startproof
Suppressing the exponent $s$ let
$B(a,b)$ be the bow energy associated to
a pair of points $(a,b \in \R^2$.

Imagine that we had a counterexample
to Equation \ref{ver1} with $x_1>0$ and $x_3<0$. 
This would mean that $$B(p_1',p_3')<B(p_1'',p_3'').$$
Let $q_1'$ and $q_3'$ denote the
points obtained by rotating the plane about the
origin until one of the points $p_1'$ or $p_3'$
hits the $y$-axis.
By rotational symmetry, we have
$$B(q_1',q_3')=B(p_1',q_3').$$
As we rotate, the vertical distance from
our points to the $x$-axis increases, because
they are ``rising up'' along circles centered
at the origin.  This means that
both $q_1''$ and $q_3''$ are further from
the origin respectively than
$p_1''$ and $p_3''$.

Note that both pairs $(p_1'',p_3'')$ and
$(q_1'',q_3'')$ are good pairs in
the sense of Lemma \ref{horincrease}.
Therefore, by Lemma \ref{horincrease},
$$B(q_1'',q_3'')> B(p_1'',p_3'').$$
Stringing together all our
inequalities, we have
$$B(q_1'',q_3'')>B(p_1'',p_3'')>B(p_1',p_3')=B(q_1',q_1'').$$
We have constructed a counterexample, namely the
pair $(q_1',q_3')$ in which one of the two points
lies on the $y$-axis.  That is what we wanted
to do.
\endproof

Having established Inequalities 1, 2, and 3, our proof
of Theorem \ref{symm1} is done.

\newpage

\chapter{The Second  Symmetrization}
\label{symm2proof}

In this chapter we prove Theorem \ref{symm2} by
establishing Inequalities 1-3 listed
in \S \ref{proofstrategy} concerning
the configurations $X''$ and $X'''$.
For reference, here are those inequalities again.
\begin{enumerate}
\item $A_s(X'')-A_s(X''') \geq -C_0-C_1$
for all $s \geq 2$,
\item $B_{02,s}(X'')-B_{02,s}(X''') \geq C_0$
for all $s \in [12,15+25/512]$. 
\item $B_{13,s}(X'')-B_{13,s}(X''') \geq C_1$
for all $s \in [12,15+25/512]$.
\end{enumerate}
The constants $C_0$ and $C_1$ depend on the
parameter $s$ and on the configuration $X''$.
As we mentioned in \S \ref{symm2} we will see that
we get equality in the second and third cases iff $X''=X'''$.
Adding up the inequalities,  we get that $R_s(X'') \geq R_s(X''')$
with equality iff $X''=X'''$.

Our argument for Inequality 3 is quite delicate.
At the end of the chapter we will explain a more
complicated but more robust argument.

\section{Inequality 1}

Figure 15.1 refers to the construction
in \S \ref{cube2}, and explains the
geometrical meaning of the variables
$a,b,c,d$ discussed there.

\begin{center}
\resizebox{!}{2.6in}{\includegraphics{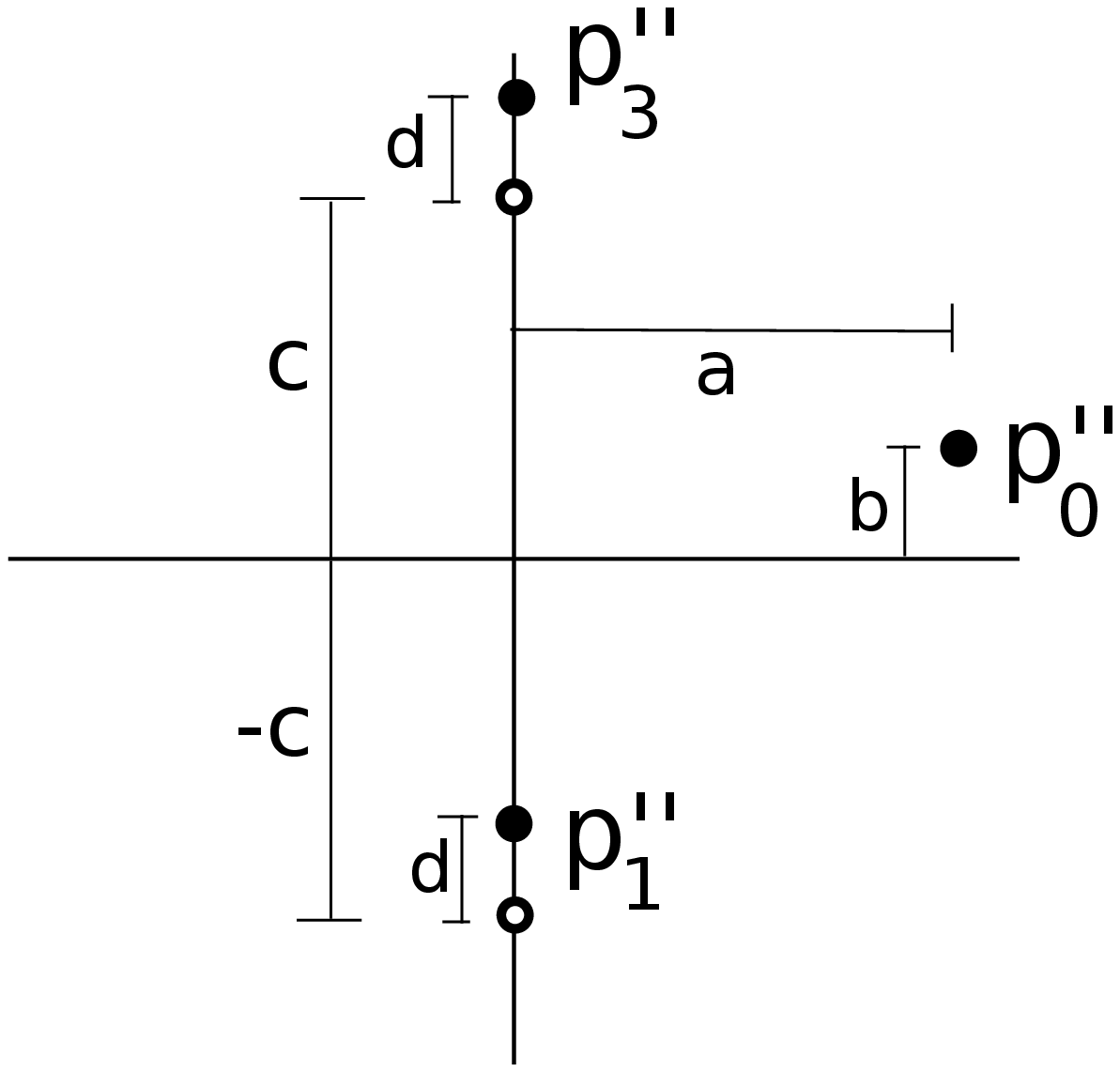}}
\newline
{\bf Figure 15.1:\/} The Meaning of the Cubical Parametrization
\end{center}

The variables in the polynomial produced by
the cubical parametrization are $t_1,t_2,t_3,t_4 \in [0,1]$.
The variablre $a,b,c,d$ are linear functions
in the variables $t_1,t_2,t_3,t_4$.   

In this chapter we are interested in studying
what happens when we hold $a$ and $c$ fixed and
vary $b$ and $d$.  The point is that our
symmetrization operation here only changes $b$ and $d$.
Suppressing the values of $a$ and $c$, we introduce
the notation

\begin{equation}
(j,k,s)_{b, d}=
\|\Sigma_{-1}(p_j'')-\Sigma^{-1}(p_k'')\|^{-s}.
\end{equation}
Here $\Sigma^{-1}$ is inverse stereographic projection.
When $b=d=0$ we write $(j,k,s)=(j,k,s)_{00}$ for
ease of notation.
We are interested in proving

\begin{equation}
(01,s)_{bd,} +
(03,s)_{bd} \geq 
(01,s) +
(03,s) \hskip 30 pt
\forall s \in [12,16],
\end{equation}
when the variables are in their intended range.
As we mentioned in \S \ref{proofstrategy}, this
calculation is too big for us.  We have to work
with a smaller exponent and claw our way back up.

We have already discussed how we need to incorporate
the bow energies into our estimate to make it work.
In practice we give a little ``positive kick'' to our
polynomial and then justify that the positive kick
is smaller than what we get from taking into account
the bow energies.  Here goes.

For each $k=1,2,3,4$, let $\phi_k$ be the map
in our cubical covering of $\Omega''_2$ described
in \S \ref{cube2}.  Define
\begin{equation}
\label{cubecompose}
\frac{P_k(t_1,t_2,t_3,t_4)}{Q_k(t_1,t_2,t_3,t_4)}=
\bigg((01,2)_{bd} +
(03,2)_{bd} -
(01,2) -
(03,2) +\frac{b^2}{94}+\frac{d^2}{13}\bigg) \circ \phi_k.
\end{equation}
The polynomials are in $\Z[t_1,t_2,t_3,t_4]$,
the expression is reduced, and
$Q$ is chosen to be positive at
$(1/2,1/2,1/2,1/2)$.  This determines
$P$ and $Q$ uniquely.

As in \S \ref{proofineq1}, we introduce the Hessians:
\begin{equation}
H_k=\left[\matrix{\partial_{22} P_k & \partial_{24} P_k \cr
\partial_{42} P_k & \partial_{44} P_k}\right], \hskip 30 pt
\partial_{ij} P_j=\frac{\partial^2 P_k}{\partial t_i \partial t_j}.
\end{equation}

\begin{theorem}
For each $k=1,2,3,4$, the functions
$\det H_k$ and ${\rm trace\/} H_k$ are weak positive
dominant on $[0,1]^4$.
\end{theorem}

\startproof
We do the calculation in Java and it works.
As a sanity check, we also compute the polynomials
in Mathematica and check that the two versions match.
\endproof

\noindent
{\bf Remark:\/} 
The integers $13$ and $94$ in our equation
are the largest integers which work.  
I went as far as checking that some of
the positive dominance fails with
$13+\frac{1}{4}$ in place of $13$.
The estimate of $94$ is extremely loose
for our purposes but the estimate of $13$
is quite tight.  Were we forced to use
an integer less than $13$ we would have
to really scramble to finish our proof.
See the discussion in \S \ref{scramble}.

\begin{corollary}
$P_k> 0$ on $(0,1)^4$.
\end{corollary}

\startproof
This has the same proof as Lemma
\ref{linesegment} once we verify the same
facts.  It follows from symmetry (and we do
a direct calculation) that $\partial_2 P_k$ and
$\partial_4 P_k$ vanish when $t_2=t_4=0$.
We choose any point $(t_1,t_2,t_3,t_4)
\in (0,1)^4$ and join it by a line segment $L$.
to $(t_1,0,t_3,0)$.  The two
functions $\det H_k$ and ${\rm trace\/} H_k$ 
are positive definition everywhere except
perhaps the endpoint of $L$.  Hence
$H_k$ is positive definite along $L$
except perhaps at the endpoint.  Hence, the
restriction of $P_k$ to $L$ is strictly
convex, except perhaps at the endpoint.
The vanishing of the first partials
combines with the convexity, to imply
that $(t_1,0,t_3,0)$ is the
global minimum of $P_k$ on $L$.
\endproof

\begin{corollary}
$$\Pi=(01,2)_{bd} +
(03,2)_{bd} -
(01,2) -
(03,2) +\frac{b^2}{94}+\frac{d^2}{13} \geq 0$$
on $\Omega_2''$.
\end{corollary}

\startproof
The function $P_k$ is positive on a connected
and open dense subset of $[0,1]^4$ and
the function $P_k/Q_k$ is well-defined
and finite on $[0,1]^4$.  Hence $Q_k$
does not vanish on a connected and open
dense subset of $[0,1]^4$.  But
we know that $Q_k>0$ for the
point $(1/2,1/2,1/2,1/2)$.  Hence
$Q_k>0$ on an open dense subset of
$[0,1]^4$.  Hence
$P_k/Q_k>0$ on an open dense subset of
$[0,1]^4$>  Hence $P_k/Q_k \geq 0$
on $[0,1]^4$. But the union of
our affine cubes contains
$\Omega_2''$.  Hence, by
definition, $\Pi \geq 0$ on $\Omega_2''$.
\endproof

Now we know that
\begin{equation}
\label{key1}
\bigg((01,2)_{bd} +
(03,2)_{bd}\bigg) -
\bigg((01,2) +
(03,2)\bigg)  \geq
-\bigg(\frac{b^2}{94} + \frac{d^2}{13}\bigg).
\end{equation}

Dividing both sides of Equation \ref{key1} by
$(01,2)$,

\begin{equation}
\label{key2}
\frac{(01,2)_{bd}}{(01,2)}+
\frac{(03,2)_{bd}}{(01,2)} - 2  \geq 
-\bigg(\frac{b^2}{94} + \frac{d^2}{13}\bigg) (01,-2).
\end{equation}

Applying the Monotonicity Lemma in the
case $N=2$ and $\alpha=2$ and $b=s$ we have

\begin{equation}
\label{key3}
\frac{(01,s)_{bd}}{(01,s)}+
\frac{(03,s)_{bd}}{(01,s)} - 2  \geq
-\frac{s}{2}\bigg(\frac{b^2}{94} + \frac{d^2}{13}\bigg) 
(01,-2).
\end{equation}

Clearing denominators, and using the fact
that $(01,s)=(03,s)$ we have
\begin{equation}
\label{key4}
(01,s)_{bd}+
(03,s)_{bd} - (01,s)-(03,s)  \geq
-\frac{s}{2}\bigg(\frac{b^2}{94} + \frac{d^2}{13}\bigg) 
 (01,s-2).
\end{equation}

The calculation above deals with half the bonds comprising
the base energy.  By reflection symmetry, we get the
same result for the other two bonds.  Hence
\begin{equation}
A_s(X'')-A_s(X''') \geq -s \bigg(\frac{b^2}{94} + \frac{d^2}{13}\bigg)
(01,s-2).
\end{equation}
We set
\begin{equation}
C_0=\frac{b^2 s}{94} (01,s-2), \hskip 30 pt
C_1=\frac{d^2 s}{13} (01,s-2).
\end{equation}
This gives us Inequality 1.

\section{Inequality 2}

We fix the variables $a,c,d$ and study the
dependence of our points on the variable $b$.
We define
\begin{equation}
(k,s)_{b}=\|\Sigma^{-1}(p''_k)-(0,0,1)\|^{-s}
\end{equation}
We also set $(k,s)=(k,s)_0$.
A direct calculation shows
\begin{equation}
(0,2)_{b}-(0,2)=\frac{b^2}{4}.
\end{equation}
Therefore
\begin{equation}
\frac{(0,2)_{b}}{(0,2)}-1 = \frac{b^2}{4} (0,-2).
\end{equation}
By the Monotonicity Lemma,
\begin{equation}
\frac{(0,s)_{b}}{(0,s)}-1  \geq
\frac{b^2s }{8} (0,-2).
\end{equation}
Hence
\begin{equation}
(0,s)_{b}-(0,s) \geq  \frac{b^2 s}{8} (0,s-2).
\end{equation}
Since $(0,s)_{b}=(2,s)_{b}$ for all $b$, we have
\begin{equation}
(2,s)_{b}-(2,s) \geq  \frac{b^2 s}{8} (0,s-2).
\end{equation}
Finally a direct calculation shows that
\begin{equation}
\label{direct}
(02,s)_{b}-(02,s)=
\frac{b^s}{16x^s} \times (2+2x^2+b^2)^{s/2} \geq 0.
\end{equation}
Adding these last $3$ equations together, we get
\begin{equation}
B_{02,s}(X'')-B_{02,s}(X''') \geq  
\frac{b^2 s}{4} (0,s-2).
\end{equation}

We want to show that the expression on the left hand side of
this equation is greater than or equal to $C_0$, with equality
if and only if $b=0$.  Comparing this equation
with Inequality 1, we see that it suffices to prove
\begin{equation}
\label{ineq2}
\frac{(0,s-2)}{(01,s-2)}>\frac{4}{94}.
\end{equation}
Setting $p_0'''=(x,0)$ and $p_3'''=(0,y)$, we compute that
\begin{equation}
\label{zeta1}
\frac{(0,s-2)}{(01,s-2)}=\bigg(\frac{x^2+y^2}{1+y^2}\bigg)^{s/2-1}.
\end{equation}
Since $x<1$, the right hand side decreases when we decrease $x$ and $y$
and increase $s$.  Therefore, the minimum occurs at
\begin{equation}
\label{minval}
x=\frac{416}{512} \hskip 30 pt
y=\frac{348}{512} \hskip 30 pt
s=15+\frac{25}{512}
\end{equation}
When we plug in these values and use Mathematica's ability to
compute numerical approximations to rational powers of
rational numbers, we get
\begin{equation}
0.1779791356057403219045637477452967293894273925...
\end{equation}
which is certainly greater than $4/94$.
This establishes Inequality 2.

\section{Improved Monotonicity}

We plan to carry out the same analysis for Inequality 3 that
we carried out for Inequality 2, but the estimate
$$\frac{(1,s)_{b}}{(1,s)}+\frac{(3,s)_{b}}{(3,s)}-2 
\geq \frac{s}{2}
\bigg(\frac{(1,2)_{b}}{(1,2)}+\frac{(3,2)_{b}}{(3,2)}-2\bigg)$$
coming from the Monotonicity Lemma is not good enough.
In this section we improve the estimate by roughly a
factor of $4$ within
a smaller range of exponents.

\begin{lemma}
\label{supermonotone}
For $s \geq 12$ we have
$$\frac{(1,s)_{b}}{(1,s)}+\frac{(3,s)_{b}}{(3,s)}-2 
\geq 2s \bigg(\frac{(1,2)_{b}}{(1,2)}+\frac{(3,2)_{b}}{(3,2)}-2\bigg)$$
\end{lemma}

\startproof
Suppose we know this result for $s=12$.  Then
for $s>12$ we have
$$\frac{(1,s)_{b}}{(1,s)}+\frac{(3,s)_{b}}{(3,s)}-2 
\geq 
\frac{s}{12}
\bigg(\frac{(1,12)_{b}}{(1,12)}+\frac{(3,12)_{b}}{(3,12)}-2\bigg) \geq $$

$$\frac{s}{12} \times 24 \times \bigg(
\frac{(1,2)_{b}}{(1,2)}+\frac{(3,2)_{b}}{(3,2)}-2\bigg) =
2s \bigg(\frac{(1,2)_{b}}{(1,2)}+\frac{(3,2)_{b}}{(3,2)}-2\bigg).$$
The first inequality is the Monotonicity Lemma.
So, we just have to prove our result when $s=12$.

Setting $p_1''=(0,-y+d)$ and $p_3''=(0,y+d)$
we compute that
$$
\frac{(1,12)_{b}}{(1,12)}+\frac{(3,12)_{b}}{(3,12)}-2 =
\frac{P(y,d)}{Q(y,d)}
\bigg(\frac{(1,2)_{b}}{(1,2)}+\frac{(3,2)_{b}}{(3,2)}-2\bigg),
$$
Where $P$ and $Q$ are entirely positive polynomials.  
The positivity is a piece of good luck.
Specifically, we have

$$P=6 + 15 d^2 + 20 d^4 + 15 d^6 + 6 d^8 + d^{10} + 
  90 y^2 + 300 d^2 y^2 + 420 d^4 y^2 + $$
$$ 270 d^6 y^2 +
66 d^8 y^2 +  300 y^4 + 1050 d^2 y^4 + 1260 d^4 y^4 + 495 d^6 y^4 + $$
\begin{equation}
 420 y^6 + 
  1260 d^2 y^6 + 924 d^4 y^6 + 270 y^8 + 495 d^2 y^8 + 66 y^{10}.
\end{equation}

\begin{equation}
Q=(1+y^2)^5
\end{equation}

Since $d$ does not appear in the denominator,
the ratio above is minimized when $d=0$.  When
$d=0$ the expression simplifies to
\begin{equation}
\frac{P}{Q}=\frac{6+66 y^2}{1+y^2}.
\end{equation}
This function is monotone increasing for $y>0$
and exceeds $24$ when we plug in $y=87/128$,
the minimum allowable value for points
in $\Omega''$.
\endproof

\section{Inequality 3}

We take $s \geq 12$.
A direct calculation shows
\begin{equation}
(1,2)_{d} + (3,2)_{d} - (1,2) - (3,2)=\frac{d^2}{2}.
\end{equation}

Using the fact that $(1,2)=(3,2)$, we have
\begin{equation}
\frac{(1,2)_{d}}{(1,2)}+\frac{(3,2)_{d}}{(1,2)}-2=
\frac{d^2}{2}(1,-2).
\end{equation}

When $s \geq 12$ Lemma \ref{supermonotone} gives

\begin{equation}
\frac{(1,s)_{d}}{(1,s)}+\frac{(3,s)_{d}}{(1,s)}-2 \geq
d^2 s \times (1,-2).
\end{equation}
Hence
\begin{equation}
\label{zoopX}
(1,s)_{d}+(3,s)_{d}-(1,s)-(3,s) \geq
d^2 s \times (1,s-2).
\end{equation}

Moreover, just as in Equation \ref{direct}, a direct
calculation shows
\begin{equation}
\label{direct2}
(13,s)_{d}-(13,s)=
\frac{d^s}{16y^s} \times (2-2y^2+d^2)^{s/2} \geq 0.
\end{equation}
This time we have to use the fact that $y^2<1$.
Adding the last two equations, we get

\begin{equation}
B_{13,s}(X'')-B_{13,s}(X''') \geq
d^2 s \times (1,s-2).
\end{equation}

We want to show that the expression on the left hand side of
this equation is greater than or equal to $C_1$, with equality
if and only if $d=0$.  Comparing this equation
with Inequality 1, we see that it suffices to prove
\begin{equation}
\label{ineq3}
\frac{(1,s-2)}{(01,s-2)}>\frac{1}{13}=0.0761239...
\end{equation}

Using the same parameters as in Inequality 2, we compute
\begin{equation}
\label{zeta2}
\frac{(1,s-2)}{(01,s-2)}=\bigg(\frac{x^2+y^2}{1+x^2}\bigg)^{s/2-1}.
\end{equation}
The minimum value occurs for the same values as in
Inequality 2, namely
\begin{equation}
\label{minval2}
x=\frac{416}{512} \hskip 30 pt
y=\frac{348}{512} \hskip 30 pt
s=15+25/512
\end{equation} 
When we plug in these values, and use Mathematica, we get
\begin{equation}
\label{squeaker}
0.0776538129762009654816474958983834370212062922715...,
\end{equation}
which exceeds $1/13$.
This establishes Inequality 3.

Our proof of Lemma \ref{symm2} is done.
This completes the proof of the Symmetrization Lemma.

\section{Discussion}
\label{scramble}

Our proof of Inequality 2 left a wide margin for error.
For instance, it works for all $s \in [2,16]$ and
indeed for $s$ considerably larger than $16$.
However, our proof of Inequality 3 seems to be quite a
close call.  This is somewhat of an illusion. As
we mentioned already, 
we just squeezed the estimates until they delivered
what we needed.   Here we sketch an alternate
proof of Inequality 3 that leaves more room for relaxation,
and incidentally covers all exponents in $[12,16]$.

We consider the affine cubes given by the arrays.
\begin{equation}
\left[\matrix{416&{\bf 470\/}\cr  0&16\cr  348&465\cr 0&\pm 24}\right],
\hskip 30 pt
\left[\matrix{416&{\bf 470\/}\cr  0&16\cr  364&449\cr 0&\pm 64}\right].
\end{equation}
The bolded numbers have changed from $498$ to $470$ and the
rest is the same as in \S \ref{cube2}.
We redefine the functions in Equation
\ref{cubecompose} but instead we use the functions
in \S \ref{cube2}.  
\begin{equation}
\label{cubecomposex}
\frac{P_k(t_1,t_2,t_3,t_4)}{Q_k(t_1,t_2,t_3,t_4)}=
\bigg((01,2)_{bd} +
(03,2)_{bd} -
(01,2) -
(03,2) +\frac{b^2}{94}+\frac{d^2}{{\bf 39\/}}\bigg) \circ \phi_k.
\end{equation}
All the positive dominance goes through, and so in
this region we could replace the constant $C_1$
from Inequality 1 with $C_1/3$.  Thus, we would have an easy
time eliminating the configurations covered by the
above affine cubes.

What remains is a subdomain of $\Omega''_2$ consisting
of points $(p_1,p_2,p_3,p_4)$ where 
$p_{01} \geq 470/512$.  Since we have
already taken care of the rest of $\Omega''_2$, for
the remaining points we can instead use the bound on
Equation \ref{zeta2} using the points
\begin{equation}
x=\frac{470}{512} \hskip 30 pt
y=\frac{348}{512} \hskip 30 pt
s=15+25/512
\end{equation} 
When we plug in these numbers we get a much better lower bound
of $0.105...$ in place of the bound in
Equation \ref{squeaker}.  This leaves considerably more breathing
room.  If we take $s=16$ in place of $s=15+25/512$ we
still get a bound which exceeds $1/13$.  Thus, this 
alternate argument completes the proof of the Symmetrization
Lemma for all $s \in [12,16]$.

I imagine that one could continue improving the Symmetrization
Lemma by tracking various bounds through a changing domain.  
For instance,
Consider the affine cubes
\begin{equation}
\left[\matrix{416&498\cr  0&16\cr  348&{\bf 512\/}\cr 0&\pm 24}\right],
\hskip 30 pt
\left[\matrix{416&498\cr  0&16\cr  364&{\bf 512\/}\cr 0&\pm 64}\right].
\end{equation}
These cubes cover a larger domain.  When we use these
cubes in place of the ones in \S \ref{cube2}, all the
arguments above go through word for word.  Thus, 
Lemma \ref{symm2} holds for configurations in which
$p_1$ and $p_3$ can go all the way up to the unit
circle.  (Unfortunately, I don't know how to do the
same for $p_0$ and $p_2$.)
This fact might helpful for an analysis of what
happens to the minimizer for very large power law
potentials.  See the discussion in
\S \ref{finalYY}.

\newpage

\part{Endgame}

\chapter{The Baby Energy Theorem}
\label{cluge}

\section{The Main Result}
\label{goal}

The goal of this chapter and the next
is to prove Lemma \ref{most0}.  In
this chapter we formulate and prove
a version of our Energy Theorem which
is adapted to the present situation.
The result here is much easier to
prove, and just comes down to
some elementary calculus. 
\newline
\newline
{\bf The Domain:\/}
Define
\begin{equation}
\Omega=[13,15+25/512] \times [43/64,1]^2.
\end{equation}
Every configuration in {\bf SMALL4\/} has
the 
form $(p_0,p_1,p_2,p_3)$ with
\begin{equation}
-p_2=p_0=(x,0), \hskip 30 pt
-p_1=p_3=(0,y), \hskip 30 pt 
(x,y) \in [43/64,1]^2.
\end{equation}
(These conditions define a somewhat larger set of
configurations than the ones in {\bf SMALL4\/}, but
we like the conditions above better because they
are simpler and more symmetric.)
We name such configurations, at the parameter $s$,
by the triple $(s,x,y)$.
The TBP corresponds to
the points $(s,1,1/\sqrt 3)$, though
this point does not lie in $\Omega$.
\newline
\newline
{\bf Blocks:\/}
We will use the hull notation developed in 
connection with our Energy Theorem.
See \S \ref{DYAD}.   So, if $X$ is a finite
set of points, $\langle X \rangle$ is the
convex hull of $X$.
We say that a {\it block\/} is 
the set of vertices of a 
rectangle solid, having the following form:
\begin{equation}
X = I \times Q \subset [0,16] \times [0,1]^2,
\end{equation}
where $I$ is the set of two endpoints of
a parameter interval and
$Q$ is the set of vertices of a
square.  We call $X$ {\it relevant\/} if
$X \subset \Omega$.  We work with the larger
domain just to have nice initial conditions for
our divide-and-conquer algorithm.

We define
\begin{equation}
|X|_1=|I|, \hskip 30 pt
|X|_2=|Q|.
\end{equation}
Here $|I|$ is the length of the
interval $\langle I \rangle$ and
$|Q|$ is the side length of the
square $\langle Q \rangle$.
\newline
\newline
{\bf Good Points:\/}
We call a point $(x_0,y_0) \in \R^2$ {\it good\/}
if it names a configuration such
that all distances between pairs of points
are at least $1.3$.  The point
$(1,1/\sqrt 3)$ certainly good because
all the corresponding distances are
at least $\sqrt 2=1.41...$.  For technical
reasons having to do with the rationality
of our calculation we will end up choosing
a good point extremely close to the
one representing the TBP and we will
work with it rather than the TBP.
For this reason, we state our result
more generally in terms of good points rather than
in terms of the TBP.
We fix some unspecified 
good point $(x_0,y_0) \in \R^2$ throughout
the chapter, and in the next chapter we will
make our choice.
\newline
\newline
{\bf The Main Result:\/}
Let $R_s(x,y)$ denote the Riesz $s$-energy
of the configuration $(x,y)$.  We set
$R(s,x,y)=R_s(x,y)$. 
Given a point $(s,x,y) \in \Omega$, we define
\begin{equation}
\Theta(s,x,y)=R(s,x,y)-R(s,x_0,y_0).
\end{equation}
Our main result is meant to hold with respect to
any good point.

\begin{theorem}[Baby Energy]
The following is true for any relevant block $X$
and any good point:
$$\min_{\langle X \rangle} \Theta \geq \min_X \Theta - 
(|X|_1^2/512+|X|_2^2).$$
\end{theorem}

\noindent
{\bf Remark:\/}
It is worth noting the lopsided nature of the
Baby Energy Theorem.  Thanks to the high exponent
in the power laws, our function is much more
linear in the $s$-direction.  Our divide-and-conquer
operation will exploit this situation.

\section{A General Estimate}
\label{genest}

Here we 
establish well-known facts about smooth functions and
their second derivatives. 

\begin{lemma}
\label{variation}
Suppose $F: [0,1] \to \R$ is a smooth function.
Then 
$$\min_{[0,1]} F \geq \min(F(0),F(1)) - \frac{1}{8} \max_{[0,1]} |F''(x)|.$$
\end{lemma}

\startproof
Replacing $F$ by $C_1(F-C_2)$ for suitable choices of
constants $C_1$ and $C_2$, we can assume without
loss of generality that
$\min(F(0),F(1))=0$ and
and $|F''| \leq 1$ on $[0,1]$. 
We want to show that
$F(x) \geq -1/8$ for all $x \in [0,1]$.

Let $a \in [0,1]$ be some global minimizer for $F$.
Replacing $F(x)$ by $F(1-x)$ if necessary,
we can  assume without loss of generality
that $a \in [1/2,1]$.
 We have $F'(a)=0$.
By Taylor's Theorem with remainder, there is some
$b \in [a,1]$ such that
$$0 \leq F(1)=F(a) + \frac{1}{2}F''(b)(1-a)^2$$
Since $(1-a)^2 \leq 1/4$ the last term
on the right lies in $[-1/8,1/8]$.
This forces $F(a) \geq -1/8$.
\endproof

\begin{lemma}
\label{secondD}
Let $X \subset \R^n$ denote the
set of vertices of a rectangular solid.
Let $d_1,...,d_n$ denote the side lengths of 
$\langle X \rangle$.  Let $F: \langle X \rangle \to \R$
be a smooth function. We have
$$\min_{\langle X \rangle} F \geq \min_{X} F - \frac{1}{8}
\sum_{i=1}^n(\max_{\langle X \rangle} |\partial^2 F/\partial x_i^2|)\ d_i^2.
$$
\end{lemma}

\startproof
The truth of the result does not change if we translate the
domain and/or range, and/or scale the coordinates
separately.  Thus, it suffices to prove the result
when $X=\{0,1\}^n$, the set of vertices of the
unit cube, and $F(0,...,0)=0$ and 
$\min_X F=0$.

We will prove the result by induction on $n$.
The warm-up lemma takes care of the case $n=1$.
Now we consider the $n$-dimensional case.
Consider some point $(x_1,...,x_n)$.
We can apply the $(n-1)$-dimensional
result to the two points
$(0,x_2,...,x_n)$ and $(1,x_2,...,x_n)$.
Remembering that $d_i=1$ for all $i$, we get
$$F(k,x_2,...,x_n) \geq  - \frac{1}{8}
\sum_{i=2}^n \max_{\langle X \rangle} |\partial^2 F/\partial x_i^2|
$$
for $k=0,1$.  At the same time, we can apply
the $1$-dimensional result to the line segment
connecting our two points.  This gives
$$F(x_1,...,x_n) \geq
\min(F(0,x_2,...,x_n),F(1,x_2,...,x_n)) -
\frac{1}{8} \max_{\langle X \rangle}  |\partial^2 F/\partial x_1^2|.
$$
Our result now follows from the triangle inequality.
\endproof

In view of Lemma \ref{secondD}, the Baby Energy Theorem
follows from these two estimates.
\begin{enumerate}
\item $|\Theta_{ss}(s,x,y)| \leq 1/64$ for all $(s,x,y) \in \Omega$.
\item $|\Theta_{xx}(s,x,y)| \leq 4$ and
$|\Theta_{yy}(s,x,y)| \leq 4$ for all $(s,x,y) \in \Omega$.
\end{enumerate}

\noindent
{\bf Notation:\/}
Here we are using the shorthand notation
\begin{equation}
\Theta_{ss}=\frac{\partial^2 \Theta}{\partial s^2}, \hskip 30 pt
\Theta_{xx}=\frac{\partial^2 \Theta}{\partial x^2}, \hskip 30 pt
\Theta_{yy}=\frac{\partial^2 \Theta}{\partial y^2}.
\end{equation}
In general, we will denote partial derivatives this way.
Note that
the partial derivative $R_s=\partial R/\partial s$ could
be confused with the Riesz energy function at the parameter $s$.
Fortunately, we will not need to consider $\partial R/\partial s$
directly and in the one section where we come close to it we
will change our notation.

\section{A Formula for the Energy}
\label{energyf}

A straightforward calculation shows that

\begin{equation}
\label{maineq}
R(s,x,y)=A(s,x)+A(s,y)+2B(s,x)+2B(s,y)+4C(s,x,y),
\end{equation}

\begin{equation}
\matrix{
A(s,x) & = &
\bigg(\frac{(1+x^2)^2}{16x^2}\bigg)^{s/2} \cr
B(s,x) & = &
\bigg(\frac{1+x^2}{4}\bigg)^{s/2} \cr
C(s,x) & = &
\bigg(\frac{(1+x^2)(1+y^2)}{4(x^2+y^2)}\bigg)^{s/2}}
\end{equation}

\section{The First Estimate}
\label{logbound}

\begin{lemma}
\label{powerbound}
Let $\psi(s,b)=b^{-s}$.
When $b \geq 1.3$ and $s \in [13,\infty)$ we have
$\psi_{ss}(s,b) \in (0,1/440)$.
\end{lemma}

\startproof
We compute
\begin{equation}
\psi_{ss}(s,b)=b^{-s}\log(b)^2>0.
\end{equation}
Fixing $b \geq 1.3$, this positive quantity is monotone
decreasing in $s$.  So, it suffices to prove
our result when $s=13$.

The equation $$\psi_{ssb}(13,b)=0$$
has its unique solution in $[1,\infty)$ at
$b=\exp(k/13)<1.3$.  Moreover, the function
$\psi_{ss}(13,b)$ tends to $0$ as $b \to \infty$.
Hence the restriction of $\psi_{ss}(13,*)$ to
$[1.3,\infty)$ takes on its max at $b=1.3$.
We compute $|\psi_ss(13,1.3)|<1/440$.
\endproof

Let $d$ denote the smallest distance
between a pair of points in a configuration
parametrized by a point in $\Omega \cup \{x_0,y_0\}$.
Finding $d$ is equivalent to finding the global
maximum $M$ of the functions $A,B,C$ and then
taking $d=M^{-1/2}$.  We will see in the next
section that we have $M<.59$.  
The global max in $[43/64,1]^2$ occurs at the vertex
$(43/64,43/64)$.
Since $M<.59$, we have
$M^{-1/2}>1.3$ on $\Omega$.  Also, by definition,
$M^{-1/2}(x_0,y_0) \geq 1.3$.
Hence $d \geq 1.3$.
But then, by Lemma \ref{powerbound},
$\Theta_{ss}(s,x,y)$ is the sum of $20$ terms,
all less than $1/440$ in absolute value.

Now, $10$ of the $20$ terms comprising $\Theta_{ss}(s,x,y)$ are
positive and $10$ are negative. 
Note that $4$ of
the configuration distances associated to $(x,y)$ are
greater or equal to $4$ of the configuration 
distances associated to 
$(x_0,y_0)$ and {\it vice versa\/}.  Hence
$4$ of the positive terms are
are less or equal to the absolute values of $4$ of
the negative terms, and {\it vice versa\/}.
Hence, $|\Theta_{xx}|$ is at most the maximum of
the following two quantities..
\begin{itemize}
\item The sum of the $6$ largest positive terms.
\item The sum of the absolute values of the $6$ largest negative
terms.
\end{itemize}
In short, 
$$|\Theta_{xx}| \leq \frac{6}{440}<\frac{1}{64}.$$
This establishes the estimate.

\section{A Table of Values}
\label{table}

We define
\begin{equation}
a(x)=A(2,x), \hskip 30 pt
b(x)=B(2,x), \hskip 30 pt
b(x,y)=C(2,x,y).
\end{equation}
These are all rational functions.

We call a function $f$ {\it easy\/} if
the two partial derivatives $f_x$ and $f_y$
are nonzero for all $(x,y) \in [43/64,1]^2$.
For an easy function, the max and min
values are achieved at the vertices.
A routine exercise in calculus shows that
$a,a_x,a_{xx},b,b_x,b_{xx},c,c_x$ are all easy.
The function $c_{xx}$ is not easy because
$$c_{xxx}=\frac{6x(x^2-y^2)(-1+y^4)}{(x^2+y^2)^4}$$
vanishes along the diagonal.
However $c_{xxy}$ does not vanish on our
domain and $c_{xxx}$ only vanishes on the
diagonal.  So, the max and min values of
$c_{xx}$ are also achieved at the vertices.
Computing at the vertices and rounding outward, we have
$$
\matrix{
a \in [.25,.3] &
a_x \in [-.33,0] &
a_{xx} \in [.5,1.97] \cr
b \in [.35,.5] &
b_x \in [.34,.5] &
b_{xx} \in [.5,.5] \cr
c \in [.5,{\bf .59\/}] &
c_x \in [-.33,0] &
c_{xx}=[0,.49]}
$$
The value in bold is the one used in the
previous section.

Let $F$ be any of $A,B,C$.
We write
\begin{equation}
F=f^u, \hskip 30 pt u=\frac{s}{2}.
\end{equation}
For instance $A=a^u$.

By the chain rule,
\begin{equation}
\label{monoX}
F_{xx}=u(u-1)f^{u-2}f_x^2 + uf^{u-1}f_{xx}
\end{equation}
This equation will let us get estimates
on $A_{xx}, B_{xx}, C_{xx}$ based on
the lookup table above. The following
technical lemma helps with this goal.

\begin{lemma}
\label{powerdecrease}
Let $g(u)=u d^{u-1}$ and
$h(u)=u(u-1)d^{u-2}$.  If
$d<.7$ the functions $g$ and $h$
are monotone decreasing as
functions of $u$ on $[6.5,\infty)$.
\end{lemma}

\startproof
Since $g=d/(u-1) \times h$, it suffices to prove our result for $h$.
Since $h>0$ the function $\log h$ is well defined.
Since $\log$ is a monotone increasing function
on $\R_+$, it suffices to prove that
$\log h(u)$ is decreasing.  We compute
$$\frac{d}{du} \log h(u) = \frac{1}{u}+\frac{1}{u-1}+\log(d).$$
Since $\log$ is monotone increasing and $d \leq .7$, we have
$$\frac{d}{du} \log h(u) \leq \frac{1}{6.5}+\frac{1}{5.5}+\log(.7)<-.02.$$
Hence $\log h$ is decreasing on $[6.5,\infty)$ as long as $d \leq .7$.
\endproof

\section{The Second Estimate}

\begin{lemma}
\label{EST}
Throught $\Omega$ we have
\begin{enumerate}
\item $0<A_{xx}<0.035$
\item $0<B_{xx}<0.47$
\item $0<C_{xx}<0.55$.
\end{enumerate}
\end{lemma}

\startproof
From our lookup table, we see that
$f_{xx} \geq 0$ for all $F \in \{A,B,C\}$ and
so all of $A_{xx},B_{xx},C_{xx}$ are non-negative.
Plugging in the values from our lookup table, and using
Lemma \ref{powerdecrease}, we have
\begin{equation}
A_{xx} \leq (13/2) (11/2)(.3)^{9/2}(.11)+ (13/2) (.3)^{11/2}(1.97)<0.035
\end{equation}
\begin{equation}
B_{xx} \leq (13/2) (11/2)(.5)^{u-2}(.25)+ (13/2) (.5)^{11/2}(.5)<0.47
\end{equation}
\begin{equation}
C_{xx} \leq (13/2) (11/2)(.59)^{9/2}(.11)+ (13/2) (.59)^{11/2}(.49)<0.55
\end{equation}
This completes the proof.
\endproof

Note that $\Theta_{xx}=R_{xx}$ because
the good point $(x_0,y_0)$ is not a function
of the variable $x$.
By Equation
\ref{maineq} we have
\begin{equation}
\label{rbound}
R_{xx}=A_{xx}+2B_{xx}+4C_{xx}
 \leq
.035 + 2(.467) + 4(.55)<3.2.
\end{equation}
We also have $R_{xx}>0$.  So,
$R_{xx} \in (0,3.2)$, which
certainly implies our
less precise estimate $|R_{xx}|<4$.
This establishes Estimate 2.  
The bound for $R_{yy}$ follows
from the bound for $R_{xx}$ and symmetry.

The proof of the
Baby Energy Theorem is done.

\newpage

\chapter{Divide and Conquer Again}

\section{The Goal}
\label{goal0}

In this chapter we
prove Lemma \ref{most0}.
The first order of business is to make a choice
of a good point.   Let
\begin{equation}
(x_0,y_0)=\bigg(1,\frac{37837}{65536}\bigg).
\end{equation}
We define the function $\Theta$ with respect to this point.
Note that $65536=2^{16}$.

The configuration $(x_0,y_0)$ is extremely close to
the TBP on account of the fact that
$$\bigg|\frac{37837}{65536}-\frac{1}{\sqrt 3}\bigg|<2^{-18}.$$
No configuration isometric to the
one represented by $Y$ lies in {\bf SMALL4\/}.
Hence, by the Dimension Reduction Lemma,
\begin{equation}
R_s(Y)>R_s(\T), \hskip 30 pt \forall s \in \bigg[13,15+\frac{25}{512}\bigg].
\end{equation}
We work with $(x_0,y_0)$ rather than $(1,1/\sqrt 3)$
because we want an entirely rational calculation.

Recall that the set
{\bf TINY4\/} consists of those points $(x,y) \in [55/64,56/64]$.
Lemma \ref{most0} is immediately implied by the
following two statements: 
\begin{enumerate}
\item $\Theta(s,x,y)>0$ if $(s,x,y) \in \Omega$ and
$s \leq 15+\frac{24}{512}$.
\item $\Theta(s,x,y)>0$ if $(s,x,y) \in \Omega$ and
$s \geq 15+\frac{24}{512}$ and $(x,y) \in {\bf TINY4\/}$.
\end{enumerate}

\section{The Ideal Calculation}

We first describe our calculation as if it were being
done on an ideal computer, without roundoff error.

Given a square $Q$, we let the vertices
of $Q$ be $Q_{ij}$ for $i,j \in \{0,1\}$.
We choose $Q_{00}$ to be the bottom left vertex -- i.e.,
the vertex with the smallest coordinates.
Given an interval $I$ we write $I=[s_0,s_1]$,
with $s_0<s_1$.
\newline
\newline
{\bf The Grading Step:\/}
Given a block $X=I \times Q \subset [0,16] \times [0,1]^2$
we perform the following evaluation,
\begin{enumerate}
\item If $I \subset [0,13]$ we pass $X$ because $X$ is irrelevant.
\item If $I \subset [15+25/512,16]$ we pass $X$ because $X$ is irrelevant.
\item If $Q$ is disjoint from $\Omega$ we pass $X$ because $X$
is irrelevant.  We test this by checking that either
$Q_{10} \leq 43/64$ or $Q_{01} \leq 43/64$.
\item If $s_0 \geq 15+24/512$ and $Q \subset {\bf TINY4\/}$ we pass $X$.
\item $s_0<13$ and $s_1>13$ we fail $X$ because we don't want to
make any computations which involve exponents less than $13$.
\item If $X$ has not been passed or failed, we
compute
$$\min_X \Theta-|X|_1^2/512-|X|_2^2.$$
If this quantity is positive we pass $X$.
Otherwise we fail $X$.
\end{enumerate}
To establish the two statements mentioned at the
end of the last section, it suffices to
find a partition of $[0,16] \times [0,1]^2$ into
blocks, all of which pass the grading step.
\newline
\newline
{\bf Subdivision:\/}
We will either divide a block $X$ into two
pieces or $4$, depending on the dimensions.
We call $X$ {\it fat\/} if 
\begin{equation}
16|X|_2>|X|_1,
\end{equation}
and otherwise thin.
If $X$ is fat we subdivide $X$ into $4$
equal equal pieces by dyadically subdividing $Q$.
If $X$ is thin we subdivide $X$ into $2$ pieces
by dyadically subdividing $I$.  We found by trial
and error that this scheme takes advantage of the
lopsided form of the Baby Energy Theorem and
produces a small partition.
\newline
\newline
{\bf The Main Algorithm:\/}
We run the same algorithm as described
in \S \ref{divideandconquer}, except
for the following differences.  
\begin{itemize}
\item
We start with the initial block $\{0,16\} \times \{0,1\}^2$.
\item Rather than use some kind of subdivision recommendation
coming from the Energy Theorem, we use the
fat/thin method just described.
\end{itemize}

\section{Floating Point Calculations}

When I run the algorithm on my $2016$ macbook pro, it
finishes in about $0.1$ seconds, takes $23149$ steps,
and produces a partition of size $15481$.
Figure 17.1 shows a slice of this partition
at a parameter $s=13+\epsilon$. Here $\epsilon$
is a tiny positive number we don't care about.
We are showing a region that is just
slightly larger than $[1/2,1]^2$.  This is
where all the detail is.  The grey blocks are
irrelevant and the white ones have passed the
calculation from the Baby Energy Theorem.

\begin{center}
\resizebox{!}{3in}{\includegraphics{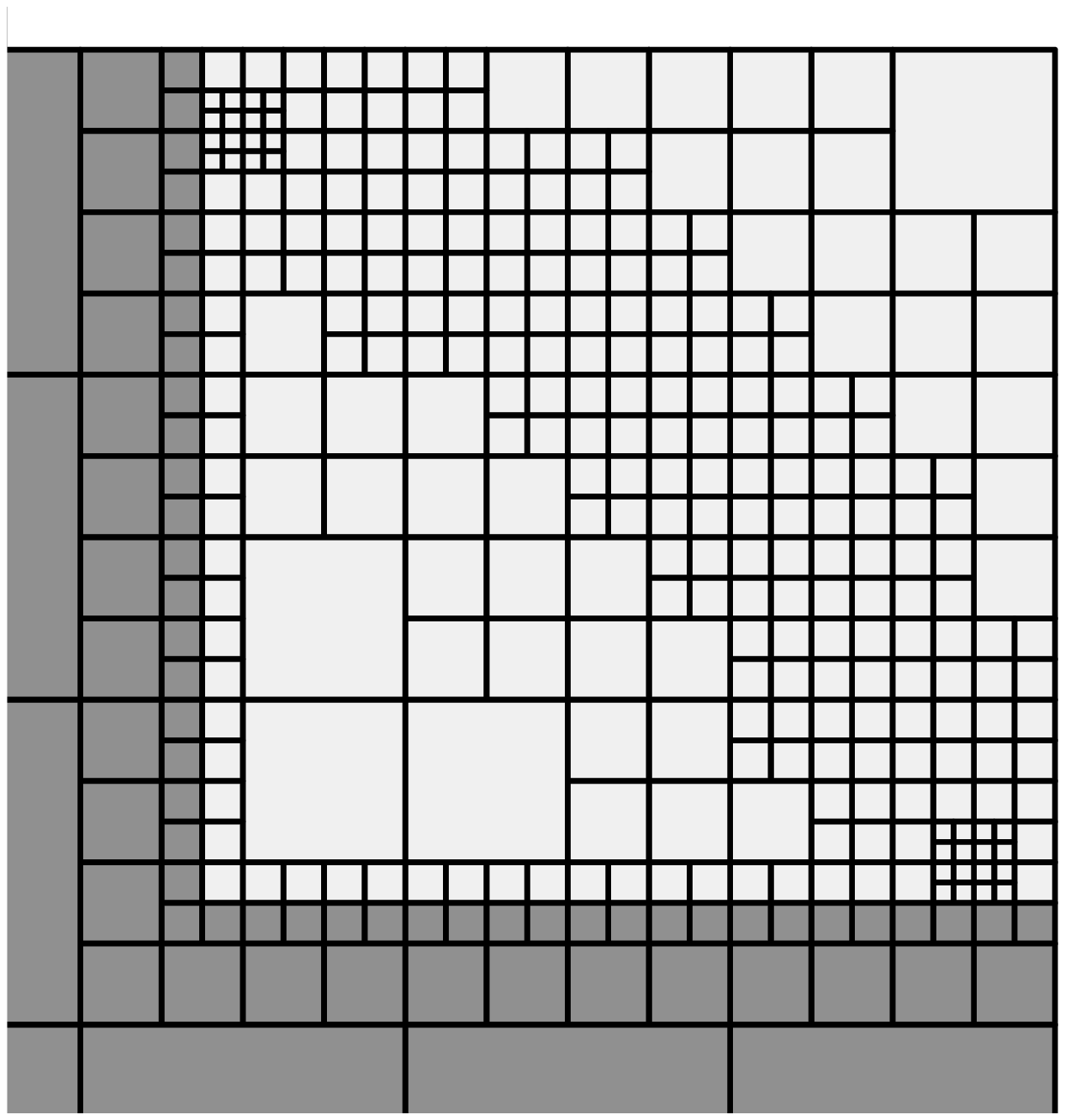}}
\newline
{\bf Figure 17.1:\/} Slice of the partition at $s=13+\epsilon$.
\end{center}

Figure 17.2 shows a slice of this partition at 
$s=15+24/512+\epsilon$.  The little squares concentrate
around the FPs which are very close to the TBP
in terms of energy.

\begin{center}
\resizebox{!}{5.2in}{\includegraphics{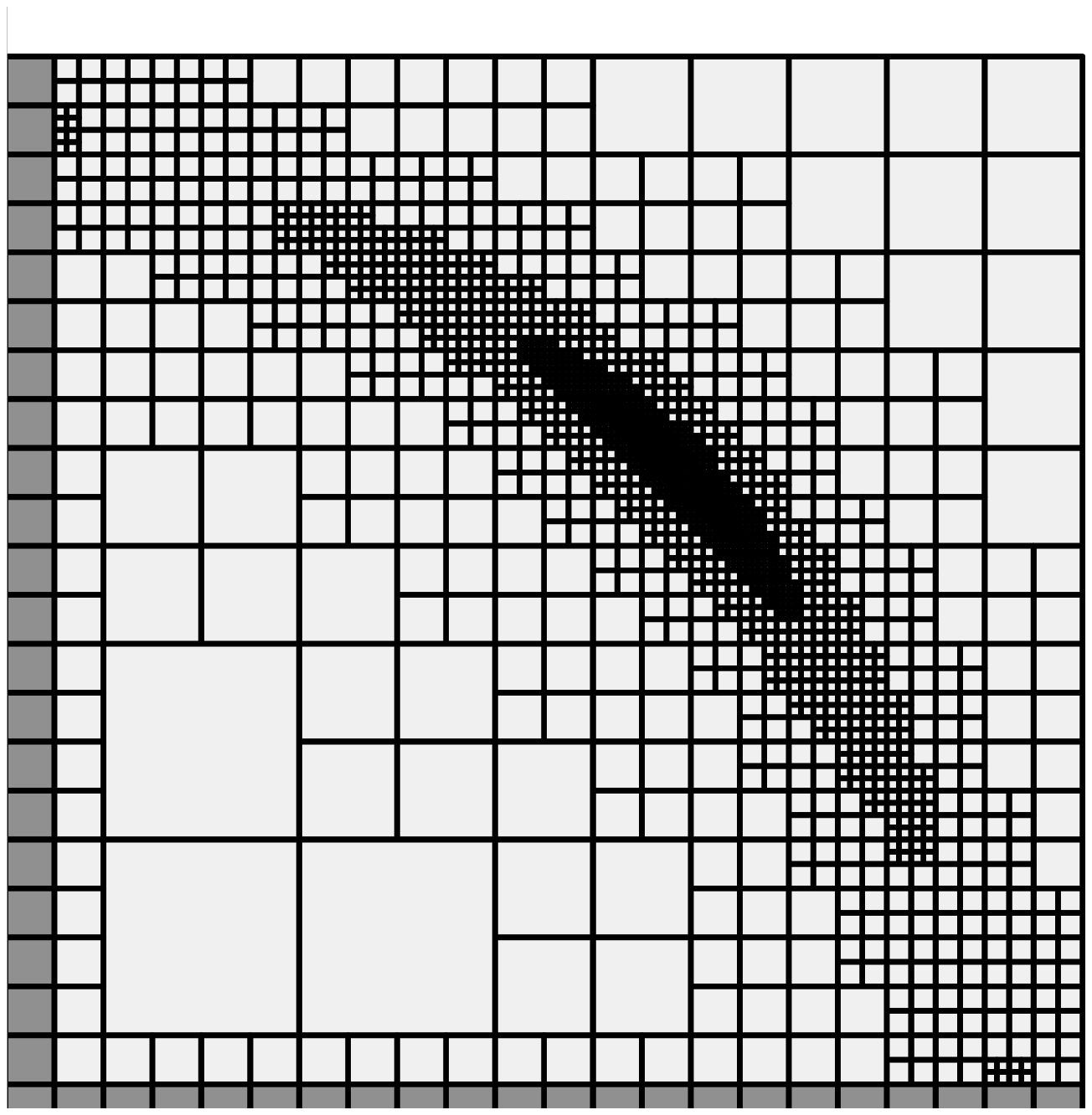}}
\newline
{\bf Figure 17.2:\/} Slice of the partition at $s=13+\epsilon$.
\end{center}

These pictures don't do the partition justice.
Using my program, the reader can slice the
partition at any parameter and zoom in and out
of the picture.

\section{Confronting the Power Function}

The floating point calculation just discussed
does not prove anything, on account
of floating point error.  The problem is that we need
to deal with expressions of the form $a^s$ where $a$ and
$s$ are rationals.  I don't know how well the floating
point errors are regulated for the power function.
Without this knowledge, it is difficult to set up an
interval arithmetic calculation.  Instead, I will
explain how to do the calculation
using exact integer arithmetic.

First of all, we represent our blocks as integers in
the following way:
\begin{equation}
2^{16}X=[S_0,S_1] \times [X_0,X_1] \times [Y_0,Y_1],
\hskip 30 pt S_i,X_j,Y_k \in \Z.
\end{equation}
In theory, the subdivision process could render
this representation invalid.  This would happen,
say, if $S_0+S_1$ is odd and we have to
subdivide $X$ along the interval $I=[S_0/2^{16},S_1/2^{16}]$.
However, our algorithm comes nowhere near doing enough
nested subdivisions for this to happen.

Referring to the functions in \S \ref{energyf}, in order
to evaluate the functions $F_k(s,x,y)$ for
$s,x,y \in \Q$ we need to be able to handle expressions
of the form
\begin{equation}
b^s, \hskip 30 pt b,s \in \Q.
\end{equation}
Now we explain how we do this.
\newline
\newline
{\bf A Specialized Power Function:\/}
We introduce the function
\begin{equation}
P(\sigma,b,s), \hskip 30 pt \sigma \in \{-1,1\}.
\end{equation}
Here is the definition:

\begin{enumerate}
\item We compute the floating point value $y=b^s$.
\item We compute the floating point value $2^{32}y$
and cast it to an integer $N$.  The casting process
rounds this value to a nearby integer -- presumably
one of the two nearest choices. 
\item $P(\sigma,b,s)=2^{-32}(N+2\sigma)$.  (Here $2\sigma$ is deliberate.)
\end{enumerate}
It seems possible that the value $P(\sigma,b,s)$ might
be computer-dependent, because of the way that $2^{32}y$
is rounded to an integer.  
\newline
\newline
{\bf Heuristic Discussion:\/}
Though this is irrelevant for the formal proof,
let me first justify
heuristically why one should expect
\begin{equation}
\label{rigor}
P(-1,b,s)<b^s<(+1,b,s),
\end{equation}
on any modern computer,
at least for smallish rationals $b$ and $s$.

Unless we had an unbelievably poor computer, we should expect
$$|N-2^{32} y|<1.5,$$
because the computer is supposed to choose an
integer within $1$ unit of $2^{32} y$.
At the same time, unless we had an extremely poor
computer, we should expect
$$|y-b^s|<.5 \times 2^{-32},$$
because otherwise the power function would be off in
roughly the $10$th digit.

Assuming these equations, we have
$$(N-2)-2^{32}y<-.5,$$
which implies
$$(N-2)-2^{32}b^s<(N-2)-2^{32}y+.5<0.$$
But this gives $P(-1,b,s)<b^s$. Similarly,
we should expect $P(1,b,s)>b^s$.
\newline
\newline
{\bf Formal Verification:\/}
The above heuristics do not play any role in our
proof, but we mentioned them in order to explain
the idea behind the construction.  What we actually
do is compute $P(\pm 1,b,s)$ and then formally
verify Equation \ref{rigor} for the relevant
numbers.  Here we explain how we check that
either $p<b^s$ or $b^s<p$ for eational numbers $p,b,s$ using
integer arithmetic. 

We define integers $\alpha,\beta,\gamma,\delta,u,v$ such
that
\begin{equation}
p=\frac{\alpha}{\beta}, \hskip 30 pt
b=\frac{\gamma}{\delta}, \hskip 30 pt
s=\frac{u}{v}.
\end{equation}
The sign of $p-b^s$ is the same as the sign of
\begin{equation}
\alpha^v \delta^u -\beta^v \gamma^u.
\end{equation}
We just expand out these integers and
check the sign.  We implement the verification
using Java's BigInteger class, which performs
exact integer arithmetic on integers whose size
is limited only by the capacity of the machine.
In theory, the needed verification could be
larger than the capacity of the machine, but
for our algorithm, the calculation runs to
completion, doing all verificaations in less
than about $2$ minutes.

\section{The Rational Calculation}

Now we describe the calculation we actually do.
We represent our blocks as in the preceding
section.  

For each of the functions $F_k$ in \S \ref{energyf}
we define $\underline F_k$ and $\overline F_k$ to
be the functions obtained by replacing 
the operation $b^s$ with $P(-1,b,s)$ and
$P(1,b,s)$ respectively.
For instance
\begin{equation}
\underline F_2(s,x)=P\bigg(-1,\frac{1+x^2}{4},s/2\bigg).
\end{equation}
We then define
$\underline R$ and $\overline R$ as in
Equation \ref{maineq} except that we replace
each occurance of a function $F$ with
$\underline F$ or $\overline F$.
For each instance we actually compute, we formally
verify all relevant instances of Equation \ref{rigor}.
(We set the calculation to abort with a conspicious
notification if any of these formal verifications fails.)
These verifications give us the inequalities
\begin{equation}
\label{rigor2}
\underline R(s,x,y)<R(s,x,y)<\overline R(s,x,y)
\end{equation}
in every instance that we need it.

We define
\begin{equation}
\underline \Theta(s,x,y)=\underline R(s,x,y) - \overline R(s,x_0,y_0).
\end{equation}
In all the cases we need, we have
\begin{equation}
\underline \Theta(s,x,y)<\Theta(s,x,y).
\end{equation}
Finally, in Step 5 of the Grading Algorithm, we pass
the block $X$ if we have
\begin{equation}
\min_X \underline \Theta - |X|_1^2/512 - |X|_2^2>0.
\end{equation}
For all the values we need, this equation guarantees
that
$$\min_X \Theta - |X|_1^2/512 - |X|_2^2>0,$$
which combines with the Baby Energy Theorem to show
that $\Theta>0$ on $\langle X \rangle$.

When we re-run our calculation using this rational
approach, the calculation runs to completion
in about $2$ minutes and $15$ seconds.  This is about
$1500$ times as long as the floating point
calculation.  The rational calculation takes
$23213$ steps and produces a partition of size
$15519$.  This partition is slightly larger
than the one produced by the floating
point calculation because the rational 
version of Grading Step 6 is slightly more stringent.

The fact that the rational version of our program
runs to completion constitutes a proof of the
two statements mentioned at the end of
\S \ref{goal0}, and thereby furnishes a proof
of Lemma \ref{most0}.

\section{Discussion}

It seems worth mentioning some other
approaches one could take when dealing
with the calculations described in this
chapter.  First of all, the algorithm
halts in about $23000$ steps, and
perhaps there are about a half a million
times we must invoke the power function.
This is a pretty small number in modern
computing terms, and so we could have
simply made a lookup table, where we
store the needed values of the power
function in a file and let the
computation access the file as needed.
This approach is somewhat like the
construction of our function $P(\sigma,b,s)$,
except that the values would be
frozen once and for all, and someone
could go back and verify that the
values stored in the file had the
claimed properties by any means available.

Another approach would be to observe that
we are only considering expressions
of the form $b^s$ where $b \in Q$ and
$s=a/2^k$.  That is, $s$ is a dyadic
rational. We could consider the identity
$$b^s=\sqrt{ ... \sqrt{b \times ... \times b}},$$
which only involves multiplication and the
square root function, both of which are
governed by the IEEE standards.  One
could use this as a basis for an
interval arithmetic calculation, though I
don't know how efficient it would be. 
Like our formal verifier, the success of
this approach would depend on $k$ being
fairly small.

Yet another approach would be to construct
polynomial over and under approximations
to the function $p(b,s)=b^s$, somewhat
along the lines of \S \ref{polyapx}.
This would be a more robust approach,
because it wouldn't depend on the denominator
of $s$ being fairly small.

\newpage

\chapter{The Tiny Piece}

\section{Overview}

Define
\begin{equation}
\Omega'=\bigg[15+\frac{24}{512},15+\frac{25}{512}\bigg] \times
\bigg[\frac{55}{64},\frac{56}{64}\bigg].
\end{equation}
Lemma \ref{tiny} concerns triples $(s,x,y) \in \Omega'$.
As in the previous chapters, we have our Riesz
energy function $R: \Omega' \to \R$, given by
$R(s,x,y)=R_s(x,y)$.

We consider the curves
\begin{equation}
\gamma(x_0,t)=(x_0+t-2 t^2,x_0-t-2t^2).
\end{equation}
Lemma \ref{tiny} says that 
\begin{equation}
\label{delicate}
R(s,\gamma(x_0,0)) \leq R(s,\gamma(x_0,t)),
\hskip 30 pt \forall (x,\gamma(x_0,t)) \in \Omega',
\end{equation}
with equality iff $t=0$.

We've mentioned already that 
Equation \ref{delicate} is delicate.  The result
does not work, for instance,
is we forget the quadratic term and
just use the more straightforward
retraction along lines of slope $-1$ in {\bf TINY4\/}.
The trouble is that the expression
\begin{equation}
\label{hess0}
R_{xx}+R_{yy}-2R_{xy}
\end{equation}
is negative at some points of $\Omega'$.
The straight line retraction would acually work
if we could use a 
region much smaller than $\Omega'$ but I didn't
think I could prove a much stronger version of
Lemma \ref{most0}.
The more complicated retraction has nicer
properties which allow it to work throughout $\Omega'$.

\section{The Proof Modulo Derivative Bounds}

When we twice differentiate Equation \ref{delicate} with
respect to $t$ we get an expression that
is more complicated than Equation \ref{hess0},
and the other terms turn out to be helpful.
What we will prove is that
$R_{tt}>0$ throughout $\Omega'$.  This means
that the function $t \to R(s,\gamma(x_0,t))$ is
convex as long as the image lies in $\Omega'$.
Since $R(s,\gamma(x_0,t))=R(s,\gamma(x_0,-t)$,
it follows from symmetry and convexity that
$R$ takes its minimum on such curves when $t=0$.
In short, Lemma \ref{tiny} follows directly
from the fact that $R_{tt}>0$ on $\Omega'$.

Using the Chain Rule (or, more honestly, Mathematica) we
compute that
$$
R_{tt}=(1-8t+16t^2)R_{xx}+
(1+8t+16t^2)R_{yy}+(-2+32t^2)R_{xy}-4(R_x+R_y).
$$
There is nothing special about the function $R$ in
this equation.  We would get the same result
for any map $t \to F(s,\gamma(x_0,t))$.

Now we simplify the expression that we need to consider.
It follows from Equation \ref{maineq} and
the results in Lemma \ref{EST} that
$R_{xx}>0$ and $R_{yy}>0$ in $\Omega'$.
We will prove below in Lemma \ref{Cbound} that
$R_{xy}>0$ on $\Omega'$ as well.

Since $R_{xx}>0$ and $R_{xy}>0$ and $R_{xy}>0$ on 
$\Omega'$, we see that
$$
R_{tt} \geq (1-8t)R_{xx}+(1+8t)R_{yy}-2R_{xy}-4(R_x+R_y).
$$
Since $2t=x-y$, see that
$R_{tt} \geq H$, where
\begin{equation}
H = (1-4x+4y)R_{xx}+(1+4x-4y)R_{yy}-2R_{xy}-4(R_x+R_y).
\end{equation}
We will complete the proof of Lemma \ref{tiny} by
showing that $H>0$ on $\Omega'$. 

 Since (as it turns out)
$H$ is bounded away from $0$ on $\Omega'$, we can
get away with just using information about the
first partial derivatives.   The goal of the next
three sections is to prove the folllowing estimates.
\begin{equation}
\label{mainbound}
|H_x|, |H_y| <16, \hskip 30 pt H_s<2.
\end{equation}
These are meant to hold in $\Omega'$.
We are not ruling out the possibility that
$H_s$ could be very negative, though a
further examination of our argument would rule this out.
Also, for what it is worth, a careful examination
of the proofs will show that really we prove these
bounds on the larger domain $[15,16] \times {\bf TINY4\/}$.

We set
\begin{equation}
\beta=15+\frac{25}{512}.
\end{equation}
The bound on $H_s$ tells us that
\begin{equation}
H(x,y,s)>H\bigg(\beta,x,y\bigg)-\frac{1}{256},
\hskip 30 pt \forall (s,x,y) \in \Omega'.
\end{equation}
So, to finish the proof of Lemma \ref{tiny} we just have
to show that
\begin{equation}
\label{finalcalc}
H\bigg(\beta,x,y\bigg)>\frac{1}{256},
\hskip 30 pt \forall x,y \in [55/64,56/64].
\end{equation}
This reduces our calculation down to a single parameter.

We compute, for the $17^2$ points
\begin{equation}
x^*,y^* \in \{55/64+i/1024,\ i=0,...,16\}
\end{equation}
that
\begin{equation}
H\bigg(\beta,x^*,y^*\bigg)>\frac{1}{64}+\frac{1}{256}.
\end{equation}
We do this calculation in Java using exact integer arithmetic,
exactly as we did the calculation for Lemma \ref{most0}.

Given any point $(\beta,x,y)$ there is some
point $(x^*,y^*)$ on our list such that
$|x-x^*| \leq 2048$ and $|y-y^*| \leq 2048$.
But then the bounds $|H_x|<16$ and $|H_y|<16$ give us
$$
H\bigg(\beta,x,y\bigg)>
H\bigg(\beta,x,y^*\bigg)-\frac{1}{128}>
H\bigg(\beta,x^*,y^*\bigg)-\frac{1}{64}>\frac{1}{256}.
$$
This establishes Equation \ref{finalcalc}.
Our proof of Lemma \ref{tiny} is done, modulo the
derivative bounds in Equation \ref{mainbound}.
The rest of the chapter is devoted to establishing
these bounds.

The bound on $|H_x|$ implies the bound on
$|H_y|$ by symmetry.  So, we just have to
deal with $|H_x|$ and $|H_s|$.
Before we launch into the tedious calculations,
we note that a precise estimate like
$|H_x|<16$ is not so important for our
overall proof.  Were we to have the weaker
estimates $|H_x|, |H_y|<32$ we
just need to compute at $4$ times as many
values. Given the speed of the computation above (less
than $30$ seconds)
we would have a feasible calculation with much
worse estimates.  We mention this because, even
though we tried very hard to avoid any
errors of arithmetic, we don't want to the reader
to think that the proof is so fragile that a
tiny error in arithmetic would destroy it.

\section{The First Bound}

We first prepare another table like the
one in \S \ref{table}.  This time we
compute on the smaller domain
$${\bf TINY4\/}=[55/64,56/64].$$  We define
the functions $a,b,c$ as in \S \ref{table}.
In addition to the functions considered in
\S \ref{table}, the functions
$a_{xxx},b_{xxx},c_{xy},c_{xxx},c_{xxy}$ are also easy
on {\bf TINY4\/}.
We have

$$
\matrix{
a \in [.254,.256] &
a_x \in [-.09,-.07] &
a_{xx} \in [.76,.83]& 
a_{xxx} \in [-3.3,-2.9] \cr
b \in [.43,.45] &
b_x \in [.42,.44] &
b_{xx} \in [.5,.5]&
b_{xxx} \in [0,0]. \cr
c \in [.50,.52] &
c_x \in [-.09,-.07] &
c_{xx} \in [.08,.11]&
c_{xxx} \in [-.012,.012] \cr
c_{xy} \in [.67,.71] &
c_{xxy} \in [-1.07,-.97]}
$$
Note that $b_{xx}=1/2$ identically
and $b_{xxx}=0$ identically.
One can deduce other bounds by symmetry.
For instance, the bounds on
$c_{y}$ are the same as the bounds on $c_x$.
Now we use these bounds to get derivative
estimates on the functions $A,B,C$ and $R$.
We already mentioned the inequality
$0<R_{xy}$ above.  This is part of our
next result.

\begin{lemma}
\label{Cbound}
$0<R_{xy}<.35$ in $\Omega'$.
\end{lemma}

\startproof
Since the functions $A$ and $B$ in 
Equation \ref{maineq} only involve one of
the variables, we have
$R_{xy}=4C_{xy}$.  We will
show that $0<C_{xy}<.087$ in $\Omega'$.
As in Lemma \ref{EST} we write
$u=s/2$ and $C=c^u$.

We compute
\begin{equation}
\label{posxy}
C_{xy}=u(u-1)c^{u-2}c_xc_y + u c^{u-1}c_{xy}
\end{equation}
Since $c_x<0$ and $c_y<0$ and $c_{xy}>0$, we have
$C_{xy}>0$.
Given the bounds in our table, and
Lemma \ref{powerdecrease}, we have
$$C_{xy}<(15/2)\times (14/2) \times (.52)^{11/2} \times (.09)^2 +
(15/2) \times (.52)^{13/2} (.71)<.087.$$
This completes the proof.
\endproof

Now we revisit Lemma \ref{EST}, using the
values from the smaller domain to get better
estimates.

\begin{lemma}
\label{EST2}
Throught $\Omega'$ we have
\begin{enumerate}
\item $0<A_{xx}<0.002$
\item $0<B_{xx}<0.138$
\item $0<C_{xx}<0.023$.
\end{enumerate}
\end{lemma}

\startproof 
The positivity
follows from Lemma \ref{EST}.  For the new upper bounds,
we proceed just as in Lemma \ref{EST} except that now
we have different values to plug in.    Using
Lemma \ref{powerdecrease} and our lookup table, we have
\begin{equation}
A_{xx} \leq (15/2)(13/2)(.256)^{11/2}(.09)^2 (15/2) (.256)^{13/2}(.83)<.002
\end{equation}

\begin{equation}
B_{xx} \leq (15/2)(13/2)(.45)^{11/2}(.44)^2+ 15/2 (.45)^{13/2}(.5)<.138
\end{equation}

\begin{equation}
C_{xx} \leq (15/2)(13/2)(.52)^{11/2}(.09)^2+ 15/2 (.45)^{13/2}(.11)<.023
\end{equation}
\endproof

\begin{corollary}
$|R_{xx}|,|R_{xx}|<.37$ on $\Omega'$.
\end{corollary}

\startproof
From
\ref{maineq} we get
\begin{equation}
|R_{xx}| \leq |A_{xx}|+2|B_{xx}|+4|C_{xx}| \leq
.002 + 2(.138) + 4(.023)<.37.
\end{equation}
\endproof

\begin{lemma}
The following is true in $\Omega'$.
\begin{enumerate}
\item $|A_{xxx}|<.011$.
\item $|B_{xxx}|<.77$.
\item $|C_{xxx}|<.052$.
\item $|C_{xxy}|<.31$.
\end{enumerate}
\end{lemma}

\startproof
As usual we set $u=s/2$ and $F=f^u$ for
each $F=A,B,C$.  We have
$$
F_{xxx}=u(u-1)(u-2)f^{u-3}f_x^3+
3u(u-1)f^{u-2}f_xf_{xx} + u f^{u-1}f_{xxx}.
$$
Hence
$$
|F_{xxx}| \leq $$
$$ u(u-1)(u-2)f^{u-3}|f_x|^3 +$$
$$3u(u-1)f^{u-2}|f_x||f_{xx}| +$$
\begin{equation}
 u f^{u-1}|f_{xxx}|.
\end{equation}
Plugging in the max values from our charts, we
get the bounds in Items 1-3.

A similar calculation gives
$$
|F_{xxy}| \leq $$
$$u(u-1)(u-2)f^{u-3}|f_x|^2|f_y| +$$
$$
2u(u-1)f^{u-2}|f_x||f_{xy}| +$$
$$
u(u-1)f^{u-2}|f_y||f_{xx}| +$$
\begin{equation}
 u f^{u-1}|f_{xxx}|.
\end{equation}
Plugging in the max values from our charts, we
get the bound in Item 4.
\endproof

\begin{corollary}
\label{triple}
$|R_{xxx}|, |R_{yyy}|<3.3$ and
$|R_{xxy}|, |R_{xyy}|<1.24$
on $\Omega'$.
\end{corollary}

\startproof
Using the bounds in the preceding lemma, we have
$$|R_{xxx}| \leq 2|A_{xxx}| + 4 |B_{xxx}| + 4|C_{xxx}|<
2(.011)+4(.76)+4(.052)<3.3$$
Similarly
$$|R_{xxy}|=4|C_{xxy}|<1.24.$$
The other two cases follow from symmetry.
\endproof

Now for the bound on $|H_x|$.
We first remind the reader of the equation for $H$.
 \begin{equation}
H = (1-4x+4y)R_{xx}+(1+4x-4y)R_{yy}-2R_{xy}-4(R_x+R_y).
\end{equation}

We compute
$$H_x= 
(1+4x-4y)R_{xxx}-4R_{xx}+
(1-4x+4y)R_{xyy}+4R_{yy}-
2R_{xxy}-
4R_{xx}-4R_{xy}.
$$
Noting that $|4x-4y| \leq 1/32$ on $\Omega'$, and using
the bounds in this section, we get
$$|H_x|<\frac{33}{32}(3.3)+4 (.37)+\frac{33}{32}(1.24)+4(.37)+2(1.24)+
4(.37)+4(.35)<16.$$
This completes the proof.

\section{The Second Bound}

For most of the calculations we use the
parameter $u=s/2$ in our calculations.
Note that $F_s=2F_u$ for any function $F$.
One principle we use repeatedly is that
$|\log(x)|$ is monotone decreasing when $x \in (0,1)$.

\begin{lemma}
$A_{xu}, C_{xu}>0$ and
$|B_{xu}|<.015$ on $\Omega$.
\end{lemma}

\startproof
Letting $F=f^u$ be any of our functions, we have
$F_x=u f^{u-1} f_x$.
Differentiating with respect to $u$, and noting
that $f_{xu}=0$ because $f$ does not depend on $u$,
we have
\begin{equation}
F_{xu}=f^{u-1}f_x \times (1+u\log f).
\end{equation}
In all cases we have
$\log f \in [-1.35,-.65]$ and this combines
$u \in [15/2,8]$ to force the
$1+u\log f<0$.

Since $a_x,c_x<0$ we see
that $A_{xu},C_{xu}>0$.
Since $u \leq 8$, we have
$$|1+u \log b|<8 |\log(.42)|-1<6.$$ 
Since $u>15/2$, Lemma \ref{powerdecrease}
and the bounds in the lookup table give
$$|B_{xu}|<(6) (.45)^{13/2}(.44)<.015.$$
This completes the proof.
\endproof

\begin{corollary}
\label{rbound1}
$-4(R_{xs}+R_{ys})<1$.
\end{corollary}

\startproof
$R_{xu}$ is the sum of $10$ terms.
$6$ of these terms are positive and
the other $4$ are greater than $-.015$.
Hence $R_{xu}>-.06$.  By symmetry
$R_{yu}>-.06$.
Hence $$
-4(R_{xs}+R_{ys})=-8(R_{xu}+R_{yu})<8 \times .12<1$$
This completes the proof.
\endproof

\begin{lemma}
$A_{xxu},B_{xxu},C_{xxu}<0$.
\end{lemma}

\startproof
Here is Equation \ref{monoX} again:
$$
F_{xx}=u(u-1)f^{u-2}f_x^2 + uf^{u-1}f_{xx}.$$
Differentiating with respect to $u$, we have
\begin{equation}
\label{partial2}
F_{xxu}=f^{u-2}f_x^2 \times (2u-1+u(u-1)\log f) +
f^{u-1}f_{xx} \times (1+u \log f).
\end{equation}
The same analysis as in the preceding lemma
shows that 
$$(2u-1)+u(u-1)\log f<0, \hskip 30 pt
(1+u \log f)<0.$$
Again, what is going on is that $\log f$ is
not too far from $-1$, and the coefficient
in front of $\log f$ is much larger than the
constant term.

Hence, in all cases, the first term on the
right hand side of Equation \ref{partial2} is
negative.  Since $f_{xx}>0$, the second term is
also negative.
\endproof

\begin{corollary}
\label{negative}
$(1+4x-4y)R_{xxs}+(1+4x-4y)R_{yys}<0$ on $\Omega'$.
\end{corollary}

\startproof
$R_{xxu}$ is the sum of $10$ terms all negative.
Hence $R_{xxu}<0$.  By symmetry $R_{yyu}<0$ as well.
But then the same goes for $R_{xxs}$ and $R_{yys}$.
More over $|4x-4y| \leq 1/128$ on $\Omega'$, so
the coefficients in front of our negative
quantities are positive.
\endproof

\noindent
{\bf Remark:\/}
Of course, we have also proved that $R_{ssx},R_{yys}<0$.
This fact will be useful in the next chapter.
\newline

\begin{lemma}
$C_{xyu} \in [-.053,0)$ on $\Omega'$.
\end{lemma}

\startproof
Treating Equation \ref{posxy} just as we
treated Equation \ref{monoX}, we have
\begin{equation}
C_{xyu}=c^{u-2}c_xc_y \times (2u-1+u(u-1)\log c) +
c^{u-1}c_{xy} \times (1+u \log c).
\end{equation}
From the bounds in the lookup table, we get
$\log(c) \in [-.70,-.65]$. From this fact, 
we get the bounds
$$(2u-1+u(u-1)\log c) \in (-25.2,0),
\hskip 30 pt (1-u \log c) \in (-4.6,0).$$
Since $c_xc_y>0$ and $c_{xy}>0$ we see that
$C_{xyu}<0$.
Finally, since $u>15/2$, Lemma \ref{powerdecrease} combines
with the values in the table to give
$$
|C_{xyu}|<(25.2) (.52)^{11/2} (.09)^2 +(4.6) (.52)^{13/2} (.71)<.053.
$$
This completes the proof.
\endproof

\begin{corollary}
$-2R_{xys}<1$ on $\Omega'$.
\end{corollary}

\startproof
We have $|R_{xyu}|=4|C_{xyu}|<.22$.
Hence $|2R_{xys}|=4 \times .22<1$.
\endproof

Our bound $H_s<2$ follows from adding the bounds
in the $3$ corollaries.  Our proof of Lemma
\ref{tiny} is now complete.  We mention one more
corollary that we will use in the next chapter.

\begin{corollary}
\label{neg2}
$R_{xys}<0$ on $\Omega'$.
\end{corollary}

\startproof
We have $R_{xys}=2R_{xyu}=8C_{xyu}<0$.
\endproof

\newpage

\chapter{The Final Battle}
\label{finalXX}

\section{Discussion}
\label{finalYY}

Before getting into the details of our proof,
it seems worth pausing to give some geometric
intuition about the competition between the
TBP and the FPs.  Suppose we start with the 
{\it equatorial\/} FP, the one
which has its apex at the north pole and
$4$ points at the equator.  Near the equator the sphere
is well approximated by an osculating
cylinder.  If we think of the points as living on the
cylinder instead of the sphere, then pushing
the $4$ points down dramatically increases the
distance to the apex without changing any of
the other distances.  When the exponent in the
power law is very high, this motion gives a
tremendous savings in the energy.
When we 
do the motion on the sphere, the tiny decrease
in some of the distances is a much smaller
effect relatively speaking.  This mismatch
between the dramatic energy decrease in some
bonds and a small energy increase in others is
enough for an FP to overtake
the TBP as the energy minimizer.  The question
is really just when this happens.

One interesting thing to observe is
the TBP beats any given FP in the long run.
That is, if $X$ is any TBP, we have $R_s(X)>R_s(\T)$
for all $s$ sufficiently large.  In all cases other
than the equatorial FP, the minimum distance
between a pair of points
in $X$ is less than the minimum distance between
pairs of points in the TBP.  In the long run, these
minimum distances dominate in the energy calculation.
(For the equatorial FP, a direct calculation
bears out that the TBP always has lower energy.)
Put another way, the amount
we need to push the points downward from the
equatorial FP described above
decreases as the power law increases.  As
$s \to \infty$, the winning FP converges to the
equatorial FP.

Finally, here are some throughts on how to continue
our analysis into the realm of very large exponents.
In playing around with the Symmetrization Lemma from
Part 3, I noticed that as the power law increases, the
domain in which our proof of the Symmetrization Lemma
seems to work stretches out to include points
near the axes and closer to the unit circle.  
Thus, the ``domain of proof'' for the Symmetrization Lemma
somehow follows the domain of interest in the minimization
problem. See the discussion at the end of
\S \ref{scramble}.  Perhaps this means that a generalization of
the Symmetrization Lemma would prove that an FP is the
global minimizer w.r.t. $R_s$ when $s>15+25/512$.

\section{Proof Overview}

Now we turn to the proof of the Main Theorem. Recall  that
\begin{equation}
\alpha=15+\frac{24}{512}, \hskip 30 pt
\beta = 15 + \frac{25}{512}.
\end{equation}
We now know that
\begin{itemize}
\item The TBP is the unique minimizer w.r.t
$R_s$ when $s \in (0,\alpha]$.
\item When $s \in [\alpha,\beta]$.
the minimizer w.r.t $R_s$ is either
the TBP or an FP.  In the latter case, the FP
lies in {\bf TINY4\/}.
\end{itemize}
It remains only to compare the TBP with the FPs and
see what happens. 
We divide the problem into $4$ regimes:
\begin{enumerate}
\item The critical interval,
$[\alpha,\beta]$.
We will prove a monotonicity result which allows
us to detect the phase-transition constant
$\shin \in (\alpha,\beta)$.
\item The small interval, 
$[\beta,16]$.  We show that the TBP
is never the minimizer w.r.t $R_s$ when $s$ is
in this interval. We do this with a $6$-term
Taylor series expansion of the energy about 
the point $(\beta,445/512,445/512)$.
\item The long interval $[16,24]$.  We do this
with a carefully engineered brute-force calculation
of the energy at $64$ specific points.
\item The end $[24,\infty)$.  Here
we give an analytic argument that is
inspired by the geometric picture
discussed above.
\end{enumerate}

In this chapter we will use Mathematica's ability
to do very high precision approximations of 
logs of rational numbers, and rational powers
of rational numbers.  We will list our calculations
with about $40$ digits of accuracy, and the reader
will see that usually we only need a few digits of
accuracy, and in all cases at most $8$ digits.
The calculations are so concrete
 that a reader who does not trust Mathematica can
re-do the calculations using some other system.
One could also use methods such as those in
\S \ref{polyapx} to do the evaluations.

\section{The Critical Interval}

Let $\Omega'$ be the region from the previous chapter.
Let $\Omega'' \subset \Omega'$ be the rectangle of
points $(s,x,x)$.  These points parametrize the FPs in
{\bf TINY4\/}.  

Let $\Theta: \Omega' \to \R$ be the function 
considered in the proof of Lemma \ref{most},
except that this time we define
$\Theta$ w.r.t the good point $(x_0,y_0)=(1,1/\sqrt 3)$
corresponding to the TBP.  In other words,
\begin{equation}
\Theta(s,x,x)=R(s,x,x)-R(s,1,1/\sqrt 3).
\end{equation}
We consider the partial derivative.
\begin{equation}
\Psi=\Theta_s.
\end{equation}

\begin{lemma}[Negative Differential]
$\Psi<0$ on $\Omega'$.
\end{lemma}

\startproof
In \S \ref{cluge} we proved that
\begin{equation}
|\Psi_s|<2^{-6}.
\end{equation}  Given this bound, and the
fact that every
parameter associated to
$\Omega'$ is within $2^{-9}$ of
$s^*$, it suffices to prove that
\begin{equation}
\label{shortmono}
\Psi(\alpha,x,x)<-2^{-15}=.0000305..., \hskip 30 pt
\forall\ x \in [55/64,56/64].
\end{equation}

We have already shown that
$R_{xxs}<0$ and $R_{yys}<0$ on $\Omega'$.
See the remark following Corollary \ref{negative}.
Moreover, Lemma
\ref{neg2} says that $R_{xys}<0$.  But now we observe that
\begin{equation}
\Psi_{xx}=R_{xxs}+R_{yys}+2R_{xys}<0.
\end{equation}
In other words, the single variable function 
\begin{equation}
\zeta(x)= \Psi(\alpha,x,x)
\end{equation} is
concave.  That is, $\zeta''(x)<0$.

We compute explicitly in Mathematica that
$$\zeta'(55/64)=-0.00134860003628414037588921953585756185...$$
$$\zeta'(56/64)=-0.00625293301806721131316180538457106515...$$
Our point here is that Mathematica can approximate the quantities
involved to arbitrarily high precision.  All we care about is
that these quantities are both negative, and for this we just need
$3$ digits of accuracy.

Since $\zeta'$ is negative at both endpoint and
$\zeta''>0$ we see that
\begin{equation}
\zeta'(x)<0, \hskip 30 pt \forall x \in [55/64,56/64].
\end{equation}
Hence $\zeta$ takes its maximum in our interval at
$55/64$.  Again using Mathematica, we compute that
$$\psi(55/64)=-0.000112582444773451000945188641412...<-2^{-15}.$$
Since this is the maximum value on our interval, 
Equation \ref{shortmono} is true.
\endproof

Now we reach the point where we define the
phase transition constant $\shin$.  We don't
get an explicit value, of course, but we
get the existence from the Negative Differential Lemma.

\begin{lemma}
There exists a constant 
$\shin \in (\alpha,\beta)$ such that
\begin{enumerate}
\item For $s \in (0,\shin)$ the TBP is the unique
global minimizer w.r.t $R_s$.
\item For $s=\shin$, both the TBP and some FP are
minimizers w.r.t $R_s$.
\item For $s \in (\shin,\beta]$,
the minimizer w.r.t $R_s$ is an FP.
\end{enumerate}
\end{lemma}

\startproof
Combining the Big Theorem, the Small Theorem, the
Symmetrization Lemma, and
Lemma \ref{most0}, we can say that the TBP is
the unique global minimizer w.r.t $R_s$ for all
$s \in (0,\alpha]$.

On the other hand, we compute in Mathematica that
\begin{equation}
\Theta\bigg(\beta,\frac{445}{512},\frac{445}{512}\bigg)
<-.00000079529630852573464048109...
\end{equation}
So, the TBP is not the minimizer w.r.t. $R_s$ when
$s=\beta$.

These two facts combine with Lemma \ref{tiny} 
to say that there exists a
constant $\shin \in (\alpha,\beta)$
which has Properties 1 and 2.  But for $s \in (\shin,\beta]$
we have
\begin{equation}
\Theta\bigg(s,\frac{445}{512},\frac{445}{512}\bigg)=
\int_{\shin}^s \frac{d\Theta}{d\mu}d\mu<0,
\end{equation}
by the Negative Differential Lemma.
So, for all such $s$, some FP beats the TBP.
This establishes Property 3.
\endproof

\noindent
{\bf Remark:\/} Of course, we found the value $445/512$ by
trial and error.  No other dyadic rational with denominator
less than $1024$ works.

\section{The Small Interval}
\label{midrange}

Our goal in this section is to prove that the
TBP is not a minimizer w.r.t $R_s$ when
$s \in [\beta,16]$.
We change our notation for partial derivatives
to one which is more convenient for the present purposes.
We write
\begin{equation}
\partial_s F=\partial F/\partial s, \hskip 30 pt
\partial^2_s F=\partial^2F/\partial s^2, ...
\end{equation}
for partial derivatives, rather than $F_s$, $F_{ss}$, etc,
as we have been doing. 
We have $\Psi=\partial \Theta/\partial s$, as above.

We start this section by getting good estimates
on $\Psi$ restricted to the ray
\begin{equation}
\bigg[\beta,\infty\bigg) \times \{(x^*,x^*)\},
\hskip 30 pt x^*=\frac{445}{512}.
\end{equation}
The idea is to expand things out in a Taylor series
about the point $(\beta,x^*,x^*)$.

Let $\{b_{ij}\}$ denote the set of
$10$ distances which arise in connection with
the TBP and let $\{b^*_{ij}\}$ denote the
set of $10$ distances which arise in connection
with the FP corresponding to $(x^*,x^*)$.
We have
\begin{equation}
\label{exact0}
\partial_s^k \Psi(\beta,x^*,x^*)=
\sum_{i<j} (-1)^{k+1}(\log b_{ij}^*)^{k+1} (b_{ij}^*)^{-\beta}-
\sum_{i<j} (-1)^{k+1}(\log b_{ij})^{k+1} b_{ij}^{-\beta}.
\end{equation}

Note that $b_{ij}>1$ and $b^*_{ij}>1$ for all $i,j$ and so
the two summands have opposite signs, and the signs
alternate with $k$.  Both summands decrease in
absolute value with $s$.  Therefore,
\begin{equation}
|\partial_s^k \Psi(s,x^*,x^*)| \leq
\max \bigg(
\bigg|\sum_{i<j} (\log b_{ij}^*)^{k+1} (b_{ij}^*)^{-\beta}\bigg|,
\bigg|\sum_{i<j} (\log b_{ij})^{k+1} b_{ij}^{-\beta}\bigg|\bigg).
\end{equation}
Call this bound $M_k$.

We set
\begin{equation}
C_k=\frac{1}{k!} \partial_k \Psi(\beta,x_0,y_0),
 \hskip 30 pt t=s-\beta.
\end{equation}
From Taylor's Theorem with Remainder, we have
\begin{equation}
\Psi(s,x^*,y^*) \leq C_0+C_1 t + C_2 t^2+C_3 t^3 + C_4 t^4 + (M_5/5!)\ t^5,
\end{equation}

Using Mathematica, we compute
$$
\matrix{
C_0 & = & -0.0001406516656256745696597771130279382483753... \cr
C_1 & = & +0.0001070798482749202926537278987632143748257... \cr
C_3 & = & -0.0000289210547274380613542311596565183497424... \cr
C_4 & = & +0.0000042894047365987515852310727996101178313... \cr
C_5 & = & -0.0000003935379322252784711450416897633466139... \cr
M_5/5! & = &+0.0000006949816058762526443113985497276969081...}
$$
From these values we see that
\begin{equation}
-|C_0|^{-1}\Psi(s) \geq
1-.77 t + .2 t^2 - .031 t^3 +.003 t^4 -.005 t^5.
\end{equation}
The polynomial on the right is positive dominant in
the sense of \S \ref{posdom0}.  Hence, it is positive for
$t \in [0,1]$.  Therefore,
\begin{equation}
\frac{\partial \Theta}{\partial s}=\Psi(s,x^*,x^*)<0, \hskip 30 pt
s \in [\beta,\beta+1].
\end{equation}
But we also know that $\Theta(\beta,x^*,x^*)<0$.
Hence
$\Theta(s,x^*,x^*)<0$ for all
$s \in [\beta,16]$.
Hdnce the TBP is not a
minimizer w.r.t. $R_s$ when $s \in [\beta,16]$.

\section{The Long Interval}

Now we prove that the TBP is not a minimizer
w.r.t. $R_s$ when $s \in [16,24]$. As we mentioned
in the beginning of the chapter, the winning FP
keeps changing.  So, we have to adjust our
domain.  Now we work in the domain
$[16,24] \times [7/8,1)^2$.  Referring to
the argument in \S \ref{logbound}, the smallest
distance that arises in connection with
configurations corresponding to points in
this domain is
$$\bigg(C(7/8,7/8)\bigg)^{-1/2}=\frac{112}{113} \sqrt 2>7/5.$$
This means that each individual summand in the $20$
terms which go into $\partial^2_s\Theta$ is at most
$(7/5)^{-s} \log(7/5)^2$ in absolute value.
Reasoning as in \S \ref{logbound}, we have
\begin{equation}
|\partial_s^2 \Theta(s,x,x)|<6 \times \log(7/5)^2 \times
(7/5)^{-s}<(7/5)^{-s}.
\end{equation}

Define
\begin{equation}
s_k=16+k/4.
\end{equation}
When $k=32$ we have $s_k=24$.
We compute explicitly in Mathematica that
\begin{equation}
-(7/5)^s \times \Theta\bigg(s,1-\frac{2}{s},1-\frac{2}{s}\bigg)-\frac{1}{64}>0,
\end{equation}
\begin{equation}
-(7/5)^s \times \Theta\bigg(s,1-\frac{2}{s-(1/4)},1-\frac{2}{s-(1/4)}\bigg)-\frac{1}{64}>0,
\end{equation}
for the $33$ values $s=s_0,...,s_{32}$.
In both cases, the min is greater than $.005$ and so this
calculation only requires $3$ digits of accuracy.

Consider the interval $[s_k,s_{k+1}]$ for $k \in \{0,...,31\}$.
The calculation above shows that
$$
\Theta\bigg(s_k,1-\frac{2}{s_k},1-\frac{2}{s_k}\bigg),
\Theta\bigg(s_k,1-\frac{2}{s_{k+1}},1-\frac{2}{s_{k+1}}\bigg)<$$
$$
-\frac{1}{64}(7/5)^{-s}<-\frac{1}{128}(7/5)^{-s}=
-\frac{1}{8} \times (7/5)^{-s} \times \bigg(\frac{1}{4}\bigg)^2<$$

\begin{equation}
-\frac{1}{8} \max_{s \in [s_{k},s_{k+1}]} 
\bigg|\partial_s^2 \Theta\bigg(s,1-\frac{2}{s},
1-\frac{2}{s}\bigg)\bigg| \times
 |s_{k+1}-s_k|^2.
\end{equation}
But then, by Lemma \ref{variation}, applied
to $-\Theta$,
\begin{equation}
\Theta\bigg(s,1-\frac{2}{s},1-\frac{2}{s}\bigg)<0,
\hskip 30 pt \forall s \in [s_k,s_{k+1}].
\end{equation}
We get this result for all $k=0,...,31$, and
so we conclude that
\begin{equation}
\Theta\bigg(s,1-\frac{2}{s},1-\frac{2}{s}\bigg)<0,
\hskip 30 pt
\forall s \in [16,24].
\end{equation}
This proves that the TBP is not a minimizer w.r.t 
$R_s$ for any $s \in [16,24]$.

\section{The End}

Letting $x_s=1-2/s$, we compute
\begin{equation}
\label{abound}
A(s,x_s)=\bigg(\frac{2-2s+s^2}{4(s-2)^2s^2}\bigg)^{s/2}<
\bigg(\frac{1}{4}-\frac{1}{20s}\bigg)^{s/2}.
\end{equation}
Our estimate holds for all $s \geq 22.0453...$
We prove this by setting the two quantities inside
the brackets equal to each other, solving for $s$,
and noting that all real roots are less than $22.0453...$.
This shows that the difference between the two quantities
does not change sign after this value.  A single evaluation
at $s=24$ confirms that the inequality goes in the
correct direction.

Next, we compute that
\begin{equation}
B(s,x_s)=
\bigg(\frac{1}{2}-\frac{1}{s}+\frac{1}{s^2}\bigg)^{s/2} \leq \bigg(\frac{1}{2}-\frac{23}{24 s}\bigg)^{s/2}.
\end{equation}
Using the same method, we establish that this estimate
holds for all $s \geq 24$.

Finally, we compute
\begin{equation}
C(s,x_s,x_s)=\bigg(\frac{2-2s+s^2}{2(s-2)^2s^2}\bigg)^{s/2}
\leq \bigg(\frac{1}{2}-\frac{1}{10 s}\bigg)^{s/2}
\end{equation}
This inequality follows from
Equation \ref{abound} and holds in the same range.

We take $s \geq 24$ so that all the inequalities above
hold. 
From Equation \ref{maineq} we have
\begin{equation}
R(s,x_s,x_s)=2^{-s/2}(a+b+c),
\end{equation}
where
\begin{equation}
\label{brute}
a=2\bigg(\frac{1}{2}+\frac{1}{10s}\bigg)^{s/2}<0.01.
\end{equation}

\begin{equation}
\label{subtle1}
b=4\bigg(1-\frac{23}{12s}\bigg)^{s/2} \leq
4\exp(-23/24)<1.54.
\end{equation}

\begin{equation}
\label{subtle2}
c=4\bigg(1+\frac{1}{5s}\bigg)^{s/2} \leq
4\exp(1/10)<4.43.
\end{equation}

Equation \ref{brute} is true by a mile.
Equations \ref{subtle1} and
\ref{subtle2} follow from the well known
fact that
\begin{equation}
\lim_{t \to \infty} \mu(t)=\exp(k), \hskip 30 pt
\mu(t)=
\bigg(1+\frac{k}{t}\bigg)^t,
\end{equation}
and that $\mu(t)$ is monotone increasing
on $(|k|,\infty)$.

Adding up these estimates, we get
$$
R_s(s,x_s,x_s)<5.98 \times 2^{-s/2}<$$
\begin{equation}
6 \times 2^{-s/2}+
3 \times 3^{-s/2} + 4^{-s/2}=R_s(\T).
\end{equation}
This proves that the TBP is not a minimizer
w.r.t. $R_s$ for $s \in [24,\infty)$.  

Our proof of the Main Theorem is done.

\newpage
\chapter*{References}

[{\bf A\/}] A. N. Andreev,
{\it An extremal property of the icosahedron\/}
East J Approx {\bf 2\/} (1996) no. 4 pp. 459-462
\newline
\newline
[{\bf BBCGKS\/}] Brandon Ballinger, Grigoriy Blekherman, Henry Cohn, Noah Giansiracusa, Elizabeth Kelly, Achill Schurmann, \newline
{\it Experimental Study of Energy-Minimizing Point Configurations on Spheres\/}, 
arXiv: math/0611451v3, 7 Oct 2008
\newline
\newline
[{\bf BDHSS\/}] P. G. Boyvalenkov, P. D. Dragnev, D. P. Hardin, E. B. Saff, M. M. Stoyanova,
{\it Universal Lower Bounds and Potential Energy of Spherical Codes\/}, 
Constructive Approximation 2016 (to appear)
\newline
\newline
[{\bf BHS\/}], S. V. Bondarenko, D. P. Hardin, E.B. Saff, {\it Mesh Ratios for Best Packings and Limits of Minimal Energy Configurations\/}, 
\newline
\newline
[{\bf C\/}] Harvey Cohn, {\it Stability Configurations of Electrons on a Sphere\/},
Mathematical Tables and Other Aids to Computation, Vol 10, No 55,
July 1956, pp 117-120.
\newline
\newline
[{\bf CK\/}] Henry Cohn and Abhinav Kumar, {\it Universally 
Optimal Distributions of Points on Spheres\/}, J.A.M.S. {\bf 20\/} (2007) 99-147
\newline
\newline
[{\bf CCD\/}] online website: \newline
http://www-wales.ch.cam.ac.uk/$\sim$ wales/CCD/Thomson/table.html
\newline
\newline
[{\bf DLT\/}] P. D. Dragnev, D. A. Legg, and D. W. Townsend,
{\it Discrete Logarithmic Energy on the Sphere\/}, Pacific Journal of Mathematics,
Volume 207, Number 2 (2002) pp 345--357
\newline
\newline
[{\bf F\"o\/}], F\"oppl {\it Stabile Anordnungen von Electron in Atom\/},
J. fur die Reine Agnew Math. {\bf 141\/}, 1912, pp 251-301.
\newline
\newline
[{\bf HS\/}], Xiaorong Hou and Junwei Shao,
{\it Spherical Distribution of 5 Points with Maximal Distance Sum\/}, 
arXiv:0906.0937v1 [cs.DM] 4 Jun 2009
\newline
\newline
[{\bf I\/}] IEEE Standard for Binary Floating-Point Arithmetic
(IEEE Std 754-1985)
Institute of Electrical and Electronics Engineers, July 26, 1985
\newline
\newline
[{\bf KY\/}], A. V. Kolushov and V. A. Yudin, {\it Extremal Dispositions of Points on the Sphere\/}, Anal. Math {\bf 23\/} (1997) 143-146
\newline
\newline
[{\bf MKS\/}], T. W. Melnyk, O. Knop, W.R. Smith, {\it Extremal arrangements of point and and unit charges on the sphere: equilibrium configurations revisited\/}, Canadian Journal of Chemistry 55.10 (1977) pp 1745-1761
\newline
\newline
[{\bf RSZ\/}] E. A. Rakhmanoff, E. B. Saff, and Y. M. Zhou,
{\it Electrons on the Sphere\/}, \newline
  Computational Methods and Function Theory,
R. M. Ali, St. Ruscheweyh, and E. B. Saff, Eds. (1995) pp 111-127
\newline
\newline
[{\bf S1\/}] R. E. Schwartz, {\it The $5$ Electron Case of Thomson's Problem\/},
Journal of Experimental Math, 2013.
\newline
\newline
[{\bf S2\/}] R. E. Schwartz, {\it The Projective Heat Map\/},
A.M.S. Research Monograph, 2017.
\newline
\newline
[{\bf S3\/}] R. E. Schwartz, {\it Lengthening a Tetrahedron\/},
Geometriae Dedicata, 2014.
\newline
\newline
[{\bf SK\/}] E. B. Saff and A. B. J. Kuijlaars,
{\it Distributing many points on a Sphere\/}, 
Math. Intelligencer, Volume 19, Number 1, December 1997 pp 5-11
\newline
\newline
\noindent
[{\bf Th\/}] J. J. Thomson, {\it On the Structure of the Atom: an Investigation of the
Stability of the Periods of Oscillation of a number of Corpuscles arranged at equal intervals around the
Circumference of a Circle with Application of the results to the Theory of Atomic Structure\/}.
Philosophical magazine, Series 6, Volume 7, Number 39, pp 237-265, March 1904.
\newline
\newline
[{\bf T\/}] A. Tumanov, {\it Minimal Bi-Quadratic energy of $5$ particles on $2$-sphere\/}, Indiana Univ. Math Journal, {\bf 62\/} (2013) pp 1717-1731.
\newline
\newline
[{\bf W\/}] S. Wolfram, {\it The Mathematica Book\/}, 4th ed. Wolfram Media/Cambridge
University Press, Champaign/Cambridge (1999)
\newline
\newline
[{\bf Y\/}], V. A. Yudin, {\it Minimum potential energy of a point system of charges\/}
(Russian) Diskret. Mat. {\bf 4\/} (1992), 115-121, translation in Discrete Math Appl. {\bf 3\/} (1993) 75-81

\end{document}